\documentclass[leqno,oneeqnum,onefignum,onetabnum]{siamart171218}

\usepackage{amsmath,amssymb,mathtools}
\usepackage[sc]{mathpazo}
\usepackage{algorithm,algorithmic}
\usepackage[sort,compress]{cite}
\usepackage{graphicx}
\usepackage{multirow}
\usepackage{makecell}
\usepackage{textcomp}
\usepackage{rotating}
\usepackage{siunitx}
\usepackage{paralist}
\usepackage{nicefrac}
\usepackage{booktabs}
\usepackage{lscape}
\usepackage{geometry}
\usepackage[hyperpageref]{backref}
\usepackage{tikz,pgfplots}
\usepackage{enumitem}
\usepackage{xr}
\newtheorem{remark}{Remark}

\graphicspath{{figures/}{tables/}}

\definecolor{red}{RGB}{200,16,46}

\makeatletter
\@ifclassloaded{standalone}{}{
\hypersetup{
linkcolor=red,
urlcolor=red,
citecolor=red}}
\makeatother%

\newcommand\sw[1]{\texttt{#1}}

\newcommand\demons{\sw{Demons}}
\newcommand\nirep{\sw{NIREP}}
\newcommand\itk{\sw{ITK}}
\newcommand\irtk{\sw{IRTK}}
\newcommand\ants{\sw{ANTs}}
\newcommand\dartel{\sw{DARTEL}}
\newcommand\elastix{\sw{elastix}}
\newcommand\tao{\sw{TAO}}
\newcommand\petsc{\sw{PETSc}}
\newcommand\accfft{\sw{AccFFT}}
\newcommand\niftyreg{\sw{NiftyReg}}
\newcommand\fair{\sw{FAIR}}
\newcommand\pyca{\sw{PyCA}}
\newcommand\defmet{\sw{deformetrica}}
\newcommand\fftw{\sw{FFTW}}
\newcommand\claire{\sw{CLAIRE}}

\newcommand{\figref}[1]{Fig.~\ref{#1}}
\newcommand{\tabref}[1]{Tab.~\ref{#1}}
\newcommand{\secref}[1]{\S\ref{#1}}
\newcommand{\algref}[1]{Alg.~\ref{#1}}
\newcommand{\runref}[1]{run~\##1}

\newcounter{runidnum}
\newcommand{\runid}{{\color{gray}\stepcounter{runidnum}\#\therunidnum}}
\newcommand{\resetrunid}{\setcounter{runidnum}{0}}
\DeclareMathOperator*{\minopt}{minimize}

\newcolumntype{R}{>{\columncolor{gray!20}}r}
\newcolumntype{L}{>{\columncolor{gray!20}}l}
\newcolumntype{C}{>{\columncolor{gray!20}}c}
\newcommand{\mcol}[2]{\multicolumn{#1}{r}{#2}}

\newcommand{\algadjust}{\centering\small\renewcommand\arraystretch{1.2}}

\newcommand{\vect}[1]{\ensuremath{\boldsymbol{#1}}}             
\newcommand{\mat}[1]{\ensuremath{\boldsymbol{#1}}}              
\newcommand{\diag}[1]{\ensuremath{\operatorname{diag}(#1)}}

\newcommand{\fun}[1]{\ensuremath{\mathcal{#1}}}                   
\newcommand{\dop}[1]{\ensuremath{\mathcal{#1}}}                   
\newcommand{\ns}[1]{\ensuremath{\mathbf{#1}}}                   
\newcommand{\fs}[1]{\ensuremath{\mathcal{#1}}}                  
\newcommand{\idiv}{\ensuremath{\nabla\cdot}}                    
\newcommand{\igrad}{\ensuremath{\nabla}}                        
\newcommand{\ilap}{\rotatebox[origin=c]{180}{$\nabla$}}         
\newcommand{\icurl}{\ensuremath{\nabla \times}}                 

\newcommand{\id}{\ensuremath{\operatorname{id}}}

\newcommand{\half}[1]{\frac{#1}{2}}
\renewcommand{\d}[1]{\mathop{}\!\mathrm{d}#1}
\newcommand{\dt}{\d{t}}
\newcommand{\dx}{\d{\vect{x}}}

\newcommand{\p} {\partial}
\newcommand{\iom}[1]{\int_{\Omega}#1\dx}

\newcommand{\iut}[1]{\int_0^1#1\dt}

\newcommand{\defeq}{\ensuremath{\mathrel{\mathop:}=}}
\newcommand{\eqdef}{\ensuremath{=\mathrel{\mathop:}}}
\newcommand{\T}{\ensuremath{\mathsf{T}}}

\newcommand{\di}[1]{\ensuremath{\mathbf{#1}}}                    
\newcommand{\dig}[1]{\ensuremath{\boldsymbol{#1}}}               

\newcommand{\bipa}{\begin{inparaenum}[(\itshape i\upshape)]}
\newcommand{\eipa}{\end{inparaenum}}
\newcommand{\bipasub}{\begin{inparaenum}[(\itshape a\upshape)]}
\newcommand{\eipasub}{\end{inparaenum}}

\newcommand{\ipoint}[1]{\textit{\textbf{#1}}.}

\newcommand{\iquote}[1]{``\emph{#1}''}

\reversemarginpar
\newcommand{\mmargin}[1]{{\marginpar{\em\tiny #1}}}\renewcommand{\mmargin}[1]{}

\usepackage{placeins}

\def\stitle{CLAIRE: A SOLVER FOR DIFFEOMORPHIC IMAGE REGISTRATION}
\title{CLAIRE: A distributed-memory solver for constrained large deformation diffeomorphic image registration\thanks{This material is based upon work supported by NIH award 5R01NS042645-14; by NSF awards CCF-1817048 and CCF-1725743; by the U.S. Department of Energy, Office of Science, Office of Advanced Scientific Computing Research, Applied Mathematics program under Award Number DE-SC0019393; by the U.S. Air Force Office of Scientific Research award FA9550-17-1-0190; and by the Simons Foundation award 586055. Any opinions, findings, and conclusions or recommendations expressed herein are those of the authors and do not necessarily reflect the views of the AFOSR, DOE, NIH, NSF, and Simons Foundation. Computing time on the Texas Advanced Computing Centers Stampede system was provided by an allocation from TACC and the NSF. This work was completed in part with resources provided by the Research Computing Data Core at the University of Houston.}}
\author{Andreas~Mang\thanks{Department of Mathematics, University of Houston, Houston, Texas 77204-5008, USA, \href{mailto:andreas@math.uh.edu}{andreas@math.uh.edu}} \and Amir~Gholami\thanks{Department of Electrical Engineering and Computer Sciences, University of California, Berkeley, CA 94720-1770, USA, \href{amirgh@berkeley.edu}{amirgh@berkeley.edu}} \and Christos~Davatzikos\thanks{Center for Biomedical Image Computing and Analytics, Department of Radiology, University of Pennsylvania, Philadephia, PA 19104-2643, USA, \href{christos.davatzikos@uphs.upenn.edu}{christos.davatzikos@uphs.upenn.edu}} \and George~Biros\thanks{Oden Institute for Computational Engineering and Sciences, University of Texas at Austin, Austin, TX 78712-1229, USA, \href{gbiros@acm.org}{gbiros@acm.org}}}

\begin{document}
\maketitle

\begin{abstract}
With this work we release \claire, a distributed-memory implementation of an effective solver for constrained large deformation diffeomorphic image registration problems in three dimensions. We consider an optimal control formulation. We invert for a stationary velocity field that parameterizes the deformation map. Our solver is based on a globalized, preconditioned, inexact reduced space Gauss--Newton--Krylov scheme.

We exploit state-of-the-art techniques in scientific computing to develop an effective solver that scales to thousand of distributed memory nodes on high-end clusters. We present the formulation, discuss algorithmic features, describe the software package, and introduce an improved preconditioner for the reduced space Hessian to speed up the convergence of our solver.

We test registration performance on synthetic and real data. We demonstrate registration accuracy on several neuroimaging datasets. We compare the performance of our scheme against different flavors of the \demons\ algorithm for diffeomorphic image registration. We study convergence of our preconditioner and our overall algorithm. We report scalability results on state-of-the-art supercomputing platforms. We demonstrate that we can solve registration problems for clinically relevant data sizes in two to four minutes on a standard compute node with 20 cores, attaining excellent data fidelity. With the present work we achieve a speedup of (on average) 5$\times$ with a peak performance of up to 17$\times$ compared to our former work.
\end{abstract}

\newcommand{\slugmaster}{\slugger{sisc}{xxxx}{xx}{x}{x--x}} 



\pagestyle{myheadings}
\thispagestyle{plain}
\markboth{ANDREAS MANG ET AL.}
{
\stitle
}

\section{Introduction}
\label{s:introduction}

\begin{figure}
\centering
\includegraphics[width=0.99\textwidth]
{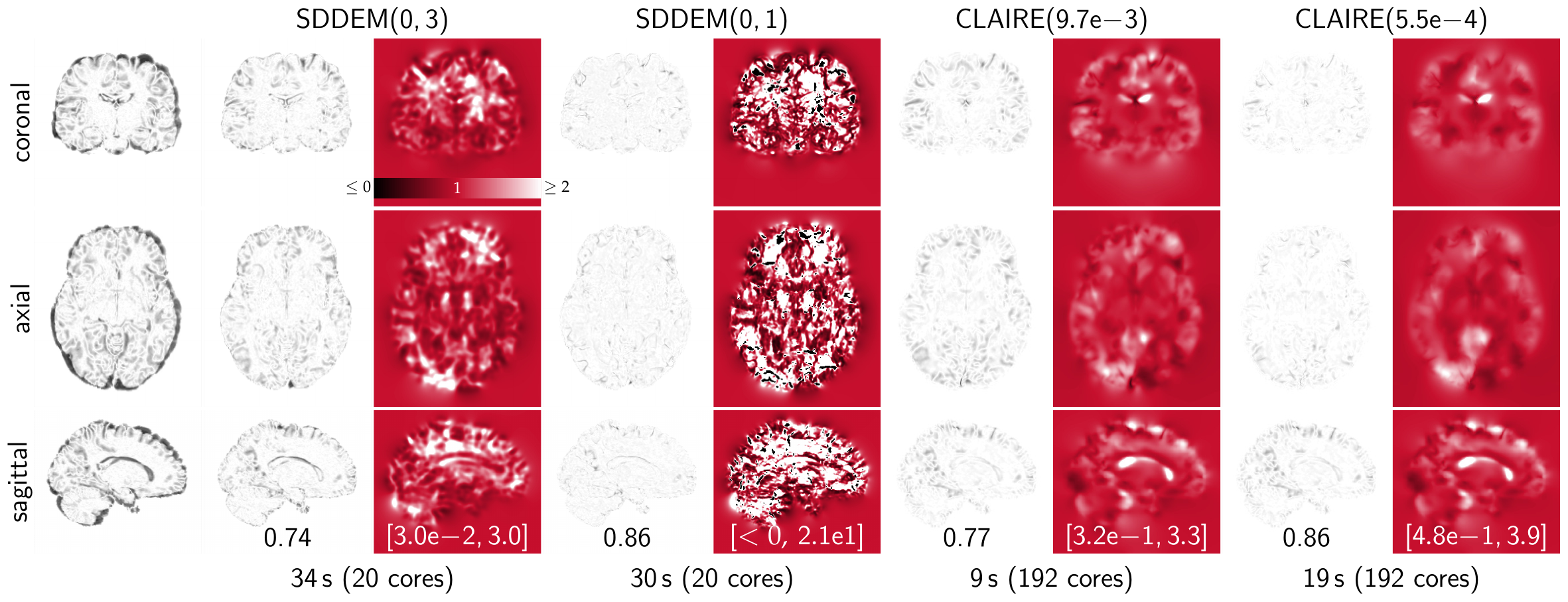}
\caption{We compare results for \claire{} and the diffeomorphic \demons{} algorithm. We consider the first two volumes of the \nirep{} dataset. We report results for the symmetric diffeomorphic \demons{} algorithm (SDDEM) with regularization parameters ($\sigma_d,\sigma_u$) determined by an exhaustive search. We report results for \claire{} for different choices for the regularization parameter for the velocity ($\beta_v=\num{3.7e-3}$ and $\beta_v=\num{5.5e-4}$; determined by a binary search). We show the original mismatch on the left. For each variant of the considered algorithms we show the mismatch after registration and a map for the determinant of the deformation gradient. We report values of the Dice score of the union of all available gray matter labels below the mismatch. We also report the extremal values for the determinant of the deformation gradient. We execute the \demons{} algorithm on one node of the RCDC's Opuntia server (Intel ten-core Xeon E5-2680v2 at \SI{2.8}{\giga\hertz} with \SI{64}{\giga\byte} memory; 2 sockets for a total of 20 cores; \cite{cacds-web}) using 20 threads. We use a grid continuation scheme with 15, 10, and 5 iterations per level, respectively. If we execute \claire{} on the same system, the runtime is \SI{103}{\second}  and \SI{202}{\second}, respectively. If we increase the number of iterations of SDDEM to 150, 100, 50 per level, we obtain a dice score of $0.75$ and $0.86$ with a runtime of \SI{322}{\second} and \SI{297}{\second}, respectively. The results for \claire{} are for 16 nodes with 12 MPI tasks per node on TACC's Lonestar 5 system (2-socket Xeon E5-2690 v3 (Haswell) with 12 cores/socket, \SI{64}{\giga\byte} memory per node; \cite{tacc-web}). We execute \claire{} at full resolution using a parameter continuation scheme in $\beta_v$. Detailed results for these runs can be found in the supplementary materials, in particular \tabref{t:nirep-regquality-h1sdiv-3d3e-3}, \tabref{t:nirep-regquality-h1sdiv-5d5e-4}, and \tabref{t:demons-all-nirep-data}.\label{f:claire-exres}}
\end{figure}

Deformable registration is a key technology in medical imaging. It is about computing a map $\vect{y}$ that establishes a \emph{meaningful} spatial correspondence between two (or more) images $m_R$ (the \emph{reference} (fixed) \emph{image}) and $m_T$ (the \emph{template} (deformable or moving) \emph{image}; image to be registered) of the same scene~\cite{Modersitzki:2004a,Fischer:2008a}. Numerous approaches for formulating and solving image registration problems have appeared in the past; we refer to~\cite{Modersitzki:2004a,Modersitzki:2009a,Fischer:2008a,Hajnal:2001a,Sotiras:2013a} for lucid overviews. Image registration is typically formulated as a variational optimization problem that consists of a data fidelity term and a Tikhonov regularization functional to over-come ill-posedness~\cite{Fischer:2008a,Engl:1996a}. In many applications, a key concern is that $\vect{y}$ is a \emph{diffeomorphism}, i.e., the map $\vect{y}$ is differentiable, a bijection, and has a differentiable inverse. A prominent strategy to ensure regularity of $\vect{y}$ is to introduce a pseudo-time variable $t \geq 0$ and invert for a smooth, time-dependent velocity field $\vect{v}$ that parameterizes the map $\vect{y}$~\cite{Beg:2005a,Dupuis:1998a,Miller:2001a,Vercauteren:2009a}; existence of a diffeomorphism $\vect{y}$ can be guaranteed if $\vect{v}$ is adequately smooth~\cite{Beg:2005a,Chen:2011a,Dupuis:1998a,Trouve:1998a}. There exists a large body of literature of diffeomorphic registration parameterized by velocity fields $\vect{v}$ that, in many cases, focuses on theoretical considerations~\cite{Dupuis:1998a,Miller:2001a,Younes:2009a,Younes:2007a,Younes:2010a}. There is much less work on the design of efficient solvers; examples are~\cite{Azencott:2010a,Avants:2011a,Ashburner:2007a,Ashburner:2011a,Beg:2005a,Crum:2005a,Zhang:2015a,Hernandez:2009a,Vercauteren:2009a,Polzin:2016a}. Most existing solvers use first order methods for numerical optimization and/or are based on heuristics that do not guarantee convergence. Due to computational costs, early termination results in compromised registration quality. Our intention in this work is to deploy an efficient solver for diffeomorphic image registration problems that \bipa\item uses state-of-the art algorithms, \item is scalable to thousands of cores, \item requires minimal parameter tuning, \item and produces high-fidelity results with guaranteed regularity on a discrete level.\eipa
We showcase exemplary results for CLAIRE for a neuroimaging dataset in \figref{f:claire-exres}. We compare \claire{} to different variants of the \demons{} algorithm.

\begin{table}
\caption{Notation and symbols.\label{t:notation-and-symbols}}
\centering\scriptsize
\begin{tabular}[t]{lllll}\toprule
Symbol                          & Description                                                             &\quad & Acronym         & Description                                                                                              \\\midrule
$\Omega$                        & spatial domain; $\Omega\defeq(0,2\pi)^3\subset\ns{R}^3$                 &      & \claire         & {\bf c}onstrained {\bf la}rge deformation diffeomorphic {\bf i}mage {\bf re}gistration \cite{claire-web} \\
$\vect{x}$                      & spatial coordinate; $\vect{x}\defeq(x_1,x_2,x_3)^\T\in\ns{R}^3$         &      & CFL             & Courant--Friedrichs--Lewy (condition)                                                                    \\
$t$                             & pseudo-time variable; $t \in [0,1]$                                     &      & CHEB($k$)       & Chebyshev (iteration) with fixed iteration number $k\in\ns{N}$ \cite{Golub:1961a,Gutknecht:2002a}        \\
$m_R(\vect{x})$                 & reference image                                                         &      & FFT             & fast Fourier transform                                                                                   \\
$m_T(\vect{x})$                 & template image (image to be registered)                                 &      & GPL             & GNU General Public License                                                                               \\
$\vect{v}(\vect{x})$            & stationary velocity field                                               &      & HPC             & high performance computing                                                                               \\
$\vect{y}(\vect{x})$            & deformation map                                                         &      & KKT             & Karush--Kuhn--Tucker                                                                                     \\
$m(\vect{x},t)$                 & state variable (transported intensities)                                &      & LDDMM           & large deformation diffeomorphic metric mapping~\cite{Beg:2005a}                                          \\
$m_1(\vect{x})$                 & final state; $m_1(\vect{x}) \defeq m(\vect{x},t=1)$                     &      & matvec          & matrix vector product                                                                                    \\
$\lambda(\vect{x},t)$           & adjoint variable                                                        &      & MPI             & Message Passing Interface                                                                                \\
$\tilde{m}(\vect{x},t)$         & incremental state variable                                              &      & \petsc          & Portable Extensible Toolkit for Scientific Computation~\cite{Balay:2016a,petsc-web}                      \\
$\tilde{\lambda}(\vect{x},t)$   & incremental adjoint variable                                            &      & PCG             & preconditioned conjugate gradient (method)~\cite{Hestenes:1952a}                                         \\
$\fun{L}$                       & Lagrangian functional                                                   &      & PCG($\epsilon$) & PCG method with relative tolerance $\epsilon\in(0,1)$                                                    \\
$\vect{g}$                      & (reduced) gradient                                                      &      & RK2             & 2nd order Runge--Kutta method                                                                            \\
$\dop{H}$                       & (reduced) Hessian operator                                              &      & (S)DDEM         & (symmetric) diffeomorphic demons~\cite{Vercauteren:2007a,Vercauteren:2009a}                              \\
$\p_i$                          & partial derivative with respect $x_i$                                   &      & (S)LDDEM        & (symmetric) log-domain diffeomorphic demons~\cite{Vercauteren:2008a}                                     \\
$\igrad$                        & gradient operator; $\igrad\defeq(\p_1,\p_2,\p_3)^\T$                    &      & \tao            & Toolkit for Advanced Optimization~\cite{Munson:2017a}                                                    \\
$\idiv$                         & divergence operator                                                     &      & \\
$\ilap$                         & Laplacian operator (vectorial and scalar)                               &      & \\
$\vect{n}_x$                    & number of grid points; $\vect{n}_x = (n_1,n_2,n_3)^\T$                  &      & \\
$n_t$                           & number of cells in temporal grid                                        &      & \\
$n$                             & number of unknowns; $n = 3\cdot \prod_{i=1}^3n_i$                       &      & \\
\bottomrule
\end{tabular}
\end{table}

\subsection{Outline of the Method}
\label{s:method-outline}

We summarize our notation and commonly used acronyms in \tabref{t:notation-and-symbols}. We use an optimal control formulation. The task is to find a smooth velocity field $\vect{v}$ (the ``control variable'') such that the distance between two images (or densities) is minimized, subject to a regularization norm for $\vect{v}$ and a deformation model given by a hyperbolic PDE constraint. More precisely, given two functions $m_R(\vect{x})$ (\emph{reference image}) and $m_T(\vect{x})$ (\emph{template image}) compactly supported on an open set $\Omega\subset\ns{R}^3$ with boundary $\p\Omega$, we solve for a \emph{stationary} velocity field $\vect{v}(\vect{x})$ as follows:
\begin{subequations}\label{e:ip}
\begin{align}\label{e:ip:objective}
\minopt_{\vect{v}, m}
& \quad \half{1}\int_{\Omega}(m_1(\vect{x}) - m_R(\vect{x}))^2\d \vect{x} + \fun{S}(\vect{v})\\
\begin{aligned}
\label{e:ip:transport}
\text{subject to}\\\\
\end{aligned}
&\quad
\begin{aligned}
\p_t m +  \vect{v} \cdot \igrad m & = 0 && \text{in}\;\Omega\times(0,1] \\
m & = m_T && \text{in}\;\Omega\times\{0\}
\end{aligned}
\end{align}
\end{subequations}

\noindent with periodic boundary conditions on $\p\Omega$. Here, $m(\vect{x},t)$ (the `'state variable`') corresponds to the transported intensities of $m_T(\vect{x})$ subject to the velocity field $\vect{v}(\vect{x})$; in our formulation, $m_1(\vect{x}) \defeq m(\vect{x},t=1)$---i.e., the solution of~\eqref{e:ip:transport} at $t=1$---is equivalent to $m_T(\vect{y}(\vect{x}))$ for all $\vect{x}$ in $\Omega$. The first part of the functional in~\eqref{e:ip:objective} measures the discrepancy between $m_1$ and $m_R$. The regularization functional $\fun{S}$ is a Sobolev norm that, if chosen appropriately, ensures that $\vect{v}$ gives rise to a diffeomorphism $\vect{y}$~\cite{Beg:2005a,Dupuis:1998a,Hart:2009a,Trouve:1998a}. We augment the formulation in~\eqref{e:ip} by constraints on the divergence of $\vect{v}$ to control volume change. A more explicit version of our formulation can be found in~\secref{s:formulation}.

Problem~\eqref{e:ip} is ill-posed and involves ill-conditioned operators. We use the method of Lagrange multipliers to solve the constrained optimization problem~\eqref{e:ip}. Our solver is based on an \emph{optimize-then-discretize} approach; we first derive the optimality conditions and then discretize in space using a pseudospectral discretization with a Fourier basis. We use a globalized, inexact, preconditioned  Gauss--Newton--Krylov method to solve for the first order optimality conditions. The hyperbolic transport equations that appear in our formulation are integrated in time using a semi-Lagrangian method. Our solver uses MPI for distributed-memory parallelism and can be scaled up to thousands of cores.

\subsection{Contributions}
\label{s:contributions}

We follow up on our former work on constrained diffeomorphic image registration~\cite{Mang:2015a,Mang:2016a,Mang:2017b,Mang:2016c}.  We focus on registration performance, implementation aspects and the deployment of our solver, and introduce additional algorithmic improvements. Our contributions are the following:
\begin{itemize}[align=left,leftmargin=1.8em,itemindent=0pt,labelsep=0pt,labelwidth=1.2em]
\item We present several algorithmic improvements compared to our past work. Most notably, we implement an improved preconditioner for the reduced space Hessian (originally described in~\cite{Mang:2017b} for the two-dimensional case). We empirically evaluate several variants of this preconditioner.
\item We evaluate registration quality and compare our new, improved solver to different variants of the diffeomorphic \demons{} algorithm~\cite{Vercauteren:2008a,Vercauteren:2009a}.
\item We study strong scaling performance of our improved solver.
\item We make our software termed \claire{}~\cite{claire-web} (which stands for {\bf c}onstrained {\bf la}rge deformation diffeomorphic {\bf i}mage {\bf re}gistration) available under GPL license. The code can be downloaded here:
\begin{center}
\url{https://github.com/andreasmang/claire}.
\end{center}
\noindent The URL for the deployment page is \url{https://andreasmang.github.io/claire}.
\end{itemize}

\subsection{Limitations and Unresolved Issues}
\label{s:limitations}

Several limitations and unresolved issues remain: $\bullet$ We assume similar intensity statistics for the reference image $m_R$ and the template image $m_T$. This is a common assumption in many deformable image registration algorithms~\cite{Beg:2005a,Hart:2009a,Lee:2010a,Museyko:2009a,Vialard:2012a}. To enable the registration of images with a more complicated intensity relationship, more involved distance measures need to be considered~\cite{Modersitzki:2004a,Sotiras:2013a}. $\bullet$ Our formulation is not symmetric, i.e., not invariant to a permutation of the reference and template image. The extension of our scheme to the symmetric case is mathematically straightforward~\cite{Avants:2008a,Lorenzi:2013b,Vercauteren:2008a} but its efficient implementation is nontrivial. This will be the subject of future work. $\bullet$ We invert for a stationary velocity field $\vect{v}(\vect{x})$ (i.e., the velocity does not change in time). Stationary paths on the manifold of diffeomorphisms are the group exponentials (i.e., one-parameter subgroups that do not depend on any metric); they do not cover the entire space of diffeomorphisms. The definition of a metric may be desirable in certain applications~\cite{Beg:2005a,Miller:2004a,Zhang:2015a} and, in general, requires nonstationary velocities. Developing an effective, parallel solver for nonstationary $\vect{v}$ requires more work.

\subsection{Related Work}
\label{s:related-work}

With this work we follow up on our prior work on constrained diffeomorphic image registration~\cite{Mang:2015a,Mang:2016a,Mang:2016c,Mang:2017a,Mang:2017c}. We release \claire, a software package for velocity-based diffeomorphic image registration. For excellent reviews on image registration see~\cite{Modersitzki:2004a,Hajnal:2001a,Sotiras:2013a}. In diffeomorphic registration, we formally require that $\det\igrad\vect{y}$ does not vanish or change sign. An intuitive approach to safeguard against nondiffeomorphic $\vect{y}$ is to add hard and/or soft constraints on $\det\igrad\vect{y}$ to the variational problem~\cite{Burger:2013a,Haber:2007a,Rohlfing:2003a,Sdika:2008a}. An alternative strategy is to introduce a pseudo-time variable $t$ and invert for a smooth velocity field $\vect{v}$ that parameterizes $\vect{y}$~\cite{Beg:2005a,Dupuis:1998a,Miller:2001a,Vercauteren:2009a}; existence of a diffeomorphism $\vect{y}$ can be guaranteed if $\vect{v}$ is adequately smooth~\cite{Beg:2005a,Chen:2011a,Dupuis:1998a,Trouve:1998a}. Our approach falls into this category. We use a PDE-constrained optimal control formulation; we refer to~\cite{Biegler:2003a,Borzi:2012a,Gunzburger:2003a,Hinze:2009a,Lions:1971a} for insight into theory and algorithmic developments in optimal control. In general, the solver has to be tailored to the structure of the control problem, which is dominated by the PDE constraints; examples for elliptic, parabolic, and hyperbolic PDEs can be found in~\cite{Adavani:2008b,Biros:2008a}, \cite{Adavani:2008a,Gholami:2016a,Mang:2012b,Stoll:2015a}, and~\cite{Benzi:2011a,Borzi:2002a,Lee:2010a,Herzog:2018a,Wilcox:2015a}, respectively. In our formulation, the PDE constraint is---in its simplest form---a hyperbolic transport equation (see \eqref{e:ip}). Our formulation has been introduced in~\cite{Mang:2015a,Mang:2016a,Hart:2009a}. A prototype implementation of our solver has been described in~\cite{Mang:2015a} and has been improved in~\cite{Mang:2017b}. We have extend our original solver~\cite{Mang:2015a} to the 3D setting in~\cite{Mang:2016c,Gholami:2017a}. The focus in~\cite{Mang:2016c,Gholami:2017a} is the scalability of our solver on HPC platforms. In~\cite{Scheufele:2019a} we presented an integrated formulation for registration and biophysical tumor growth simulations that has been successfully applied to segmentation of neuroimaging data~\cite{Mang:2017c,Gholami:2019a}.

Optimal control formulations that are related to ours have been described in~\cite{Borzi:2002a,Chen:2011a,Hart:2009a,Lee:2010a,Lee:2011a,Vialard:2012a,Herzog:2018a}. Related formulations for optimal mass transport are described in~\cite{Benzi:2011a,Haber:2015a,Rehman:2009a,Mang:2017a}. Our work differs from optimal mass transport in that intensities are constant along the characteristics (i.e., mass is not preserved). Our formulation shares numerous characteristics with traditional optical flow formulations~\cite{Horn:1981a,Kalmoun:2011a,Ruhnau:2007a}. The key difference is that we treat the transport equation for the image intensities as a hard constraint. PDE-constrained formulations for optical flow, which are equivalent to our formulation, are described in~\cite{Andreev:2015a,Barbu:2016a,Borzi:2002a,Chen:2011a}. Our work is closely related to the LDDMM approach~\cite{Avants:2008a,Avants:2011a,Beg:2005a,Dupuis:1998a,Trouve:1998a,Younes:2007a}, which builds upon the pioneering work in~\cite{Christensen:1996a}. LDDMM uses a nonstationary velocity but there exist variants that use stationary $\vect{v}$~\cite{Arsigny:2006a,Ashburner:2007a,Hernandez:2009a,Lorenzi:2013a,Lorenzi:2013b,Vercauteren:2009a}; they are more efficient. If we are only interested in registering two images, stationary $\vect{v}$ produce good results. Another strategy to reduce the size of the search space is geodesic shooting~\cite{Ashburner:2011a,Miller:2006a,Vialard:2012a,Younes:2007a,Zhang:2015b}; the control variable of the associated optimal control problem is an initial momentum/velocity at $t=0$.

Among the most popular, publicly available packages for diffeomorphic registration are \demons~\cite{Vercauteren:2008a,Vercauteren:2009a}, \ants~\cite{Avants:2011a}, \pyca~\cite{pyca-git}, \defmet~\cite{Bone:2018b,Fishbaugh:2017a,deformetrica-git} and \dartel~\cite{Ashburner:2007a}. Other popular packages for deformable registration are \irtk~\cite{Rueckert:1999a}, \elastix~\cite{Klein:2010a}, \niftyreg~\cite{Modat:2010a}, and \fair~\cite{Modersitzki:2009a}. The latter are, with the exception of \fair, based on (low-dimensional) parametric deformation models. Unlike existing approaches, \claire{} features explicit control on the determinant of the deformation gradient; we introduce hard constraints on the divergence of $\vect{v}$. Our formulation was originally proposed in~\cite{Mang:2016a}; a similar approach is described in~\cite{Borzi:2002a}. Other works that consider divergence-free $\vect{v}$ have been described in~\cite{Chen:2011a,Hinkle:2009a,Mansi:2011a,Ruhnau:2007a,Saddi:2008a}.

There exist few works on effective numerical methods. Despite the fact that first order methods for optimization have poor convergence rates for nonlinear, ill-posed problems, most works, with the exception of ours~\cite{Mang:2015a,Mang:2016a,Mang:2018a,Mang:2016c,Mang:2017a,Mang:2017b} and~\cite{Ashburner:2011a,Benzi:2011a,Hernandez:2014a,Simoncini:2012a,Vercauteren:2009a,Herzog:2018a}, use first order gradient descent-type approaches. We use a globalized Newton--Krylov method, instead. For these methods to be efficient, it is critical to design an effective preconditioner. (We refer the reader to~\cite{Benzi:2005a} for an overview on preconditioning of saddle point problems.) Preconditioners for problems similar to ours can be found in~\cite{Benzi:2011a,Simoncini:2012a,Herzog:2018a}. Another critical component is the PDE solver. In our case, the expensive PDE operators are hyperbolic transport equations. Several strategies to efficiently solve these equations have been considered in the past~\cite{Borzi:2002a,Mang:2015a,Mang:2017a,Mang:2016a,Polzin:2016a,Hart:2009a,Benzi:2011a,Simoncini:2012a,Beg:2005a,Chen:2011a,Mang:2017c,Mang:2016c}. We use a semi-Lagrangian scheme~\cite{Beg:2005a,Chen:2011a,Mang:2017c,Mang:2016c}.

Another key feature of \claire{} is that it can be executed in parallel~\cite{Mang:2016c,Gholami:2017a}. Examples for parallel solvers for PDE-constrained optimization problems can be found in~\cite{Akcelik:2002a,Akcelik:2006a,Biros:1999a,Biros:2005a,Biros:2005b,Biegler:2003a,Biegler:2007a,Schenk:2009a}. We refer the reader to~\cite{Eklund:2013a,Fluck:2011a,Shackleford:2013a,Shams:2010a} for surveys on parallel algorithms for image registration. Implementations, such as \demons~\cite{Vercauteren:2008a,Vercauteren:2009a}, \ants~\cite{Avants:2011a}, or \elastix~\cite{Klein:2010a}, which are largely based on kernels implemented in the \itk\ package~\cite{Johnson:2015a}, exploit multithreading for parallelism. GPU implementations of different variants of map-based, low-dimensional parametric approaches are described in~\cite{Shackleford:2010a,Modat:2010a,Shamonin:2014a}. A GPU implementation of a map-based nonparametric approach is described in~\cite{Koenig:2018a}. GPU implementations with formulations that are similar to ours are described in~\cite{Ha:2009a,Ha:2011a,Sommer:2011a,Rehman:2009a,ValeroLara:2014a,Bone:2018b}. The work that is most closely related to ours, is~\cite{Ha:2009a,Ha:2011a,ValeroLara:2014a}. In~\cite{Ha:2009a,Ha:2011a} a (multi-)GPU implementation of the approach described in~\cite{Joshi:2005a} is presented. The work in~\cite{ValeroLara:2014a} discusses a GPU implementation of \dartel~\cite{Ashburner:2007a}.

What sets our work apart are the numerics and our distributed-memory implementation: We use high-order numerical methods (second order time integration, cubic interpolation, and spectral differentiation). The linear solvers and the Gauss--Newton optimizer are built on top of \petsc~\cite{Balay:2016a} and \tao~\cite{Munson:2017a}. Our solver uses MPI for parallelism and has been deployed to HPC systems~\cite{Mang:2016a,Gholami:2017a}. This allows us to target applications of unprecedented scale (such as CLARITY imaging~\cite{Tomer:2014a}) without posing the need to downsample the data~\cite{Kutten:2017a}. We will see that we can solve problems with $\num{3221225472}$ unknowns in \SI{2}{\minute} on 22 compute nodes (256 MPI tasks) and in less than \SI{5}{\second} if we use 342 compute nodes (4096 MPI tasks). Exploiting parallelism also allows us to deliver runtimes that approach real-time capabilities.

\subsection{Outline}
\label{s:outline}

We present our approach for large deformation diffeomorphic image registration in~\secref{s:methods}, which comprises the formulation of the problem (see \secref{s:formulation}), a formal presentation of the optimality conditions (see \secref{s:optsys-and-newtonstep}), and a discussion of the numerical implementation (see \secref{s:numerics}). We present details about our software package in \secref{s:software}. Numerical experiments are reported in~\secref{s:experiments}. We conclude with~\secref{s:conclusions}. This publication is accompanied by supplementary materials. There, we report more detailed results and provide some background material.

\section{Methods}
\label{s:methods}

In what follows, we describe the main building blocks of our formulation, our solver, and its implementation, and introduce new features that distinguish this work from our former work~\cite{Mang:2015a,Mang:2016a,Mang:2016c,Gholami:2017a,Mang:2017a,Mang:2017b,Mang:2018a}.

\subsection{Formulation}
\label{s:formulation}

Given two images---the \emph{reference image} $m_R(\vect{x})$ and the \emph{template image} $m_T(\vect{x})$---com\-pact\-ly supported on $\Omega = (0,2\pi)^3\subset\ns{R}^3$, with boundary $\p\Omega$ and closure $\bar{\Omega}$, our aim is to compute a \emph{plausible} deformation map $\vect{y}(\vect{x})$ such that for all $\vect{x}\in\Omega$, $m_R(\vect{x}) \approx m_T(\vect{y}(\vect{x}))$~\cite{Modersitzki:2004a,Modersitzki:2009a,Fischer:2008a}. We consider a map $\vect{y}$ to be plausible if it is a diffeomorphism, i.e., an invertible map, which is continuously differentiable (a $C^1$-function) and maps $\Omega$ onto itself. In our formulation, we do not directly invert for $\vect{y}$; we introduce a pseudo-time variable $t\in[0,1]$ and invert for a \emph{stationary} velocity $\vect{v}(\vect{x})$, instead. In particular, we solve for $\vect{v}(\vect{x})$ and a mass source map $w(\vect{x})$ as follows~\cite{Mang:2016a}:
\begin{subequations}
\label{e:varopt}
\begin{align}
\label{e:varopt:objective}
\minopt_{m,\,\vect{v},\,w}&\quad
  \half{1} \iom{(m_1 - m_R)^2}
+ \half{\beta_v} \langle\dop{B}[\vect{v}],\dop{B}[\vect{v}]\rangle_{L^2(\Omega)^s}
+ \half{\beta_w} \iom{\igrad w\cdot \igrad w + w^2}\\
\begin{aligned}
\label{e:varopt:constraint}
\text{subject to}\\\\\\
\end{aligned} & \quad
\begin{aligned}
\p_t m + \igrad m \cdot \vect{v} & = 0   && {\rm in}\;\Omega\times(0,1]\\
                               m & = m_T && {\rm in}\;\Omega\times\{0\}\\
                  \idiv \vect{v} & = w   && {\rm in}\;\Omega
\end{aligned}
\end{align}
\end{subequations}

\noindent with periodic boundary conditions on $\p\Omega$, and $s>0$, $\beta_v > 0$, $\beta_w > 0$. The state variable $m(\vect{x},t)$ in~\eqref{e:varopt:constraint} represents the transported intensities of $m_T$ subjected to the velocity field $\vect{v}$; the solution of the first equation in~\eqref{e:varopt:constraint}, i.e., $m_1(\vect{x}) \defeq m(\vect{x}, t=1)$, is equivalent to $m_T(\vect{y}(\vect{x}))$, where $\vect{y}$ is the Eulerian (or pullback) map. We use a squared $L^2$-distance to measure the proximity between $m_1$ and $m_R$. The parameters $\beta_v>0$ and $\beta_w>0$ control the contribution of the regularization norms for $\vect{v}$ and $w$. The constraint on the divergence of $\vect{v}$ in~\eqref{e:varopt:constraint} allows us to control the compressibility of $\vect{y}$. If we set $w$ in~\eqref{e:varopt:constraint} to zero $\vect{y}$ is incompressible, i.e., for all $\vect{x}\in\Omega$, $\det\igrad\vect{y}(\vect{x})=1$, up to numerical accuracy~\cite{Gurtin:1981a}. By introducing a nonzero mass-source map $w$, we can relax this model to near-incompressible diffeomorphisms $\vect{y}$; the regularization on $w$ in~\eqref{e:varopt:objective} acts like a penalty on the divergence of $\vect{v}$; we use an $H^1$-norm.

Our solver supports different Sobolev \mbox{(semi-)}norms to regularize $\vect{v}$. The choice of the differential operator $\dop{B}$ not only depends on application requirements but is also critical from a theoretical point of view; an adequate choice guarantees existence and uniqueness of an optimal solution of the control problem~\cite{Barbu:2016a,Beg:2005a,Borzi:2002a,Chen:2011a,Lee:2010a} (subject to the smoothness properties of the images). We use an $H^1$-seminorm, i.e., $\dop{B} = \igrad$, if we consider the incompressibility constraint. If we neglect the incompressibility constraint, we use $\dop{B} = -\ilap$. We note that \claire{} also features $H^3$ regularization operators, and Helmholtz-type operators (i.e., regularization operators of the form $\dop{B} = -\ilap + \gamma\mat{I}$, $\mat{I} \defeq \diag{1,1,1}\in\ns{R}^{3,3}$, $\gamma > 0$, as used, e.g., in \cite{Beg:2005a}).

\subsection{Optimality Condition and Newton Step}
\label{s:optsys-and-newtonstep}

We use the method of Lagrange multipliers~\cite{Lions:1971a} to turn the constrained problem~\eqref{e:varopt} into an unconstrained one; neglecting boundary conditions, the Lagrangian functional is given by
\begin{align}
\label{e:lagrangian}
\fun{L}[m,\lambda,p,w,\vect{v}] \defeq &
  \half{1} \iom{(m_1 - m_R)^2}
+ \half{\beta_v} \big\langle\dop{B}[\vect{v}],\dop{B}[\vect{v}]\big\rangle_{L^2(\Omega)^s}
+ \half{\beta_w} \iom{\igrad w\cdot \igrad w + w^2}
\\\nonumber
+ &\int_0^1\!
\langle
\p_t m + \igrad m \cdot \vect{v},\lambda
\rangle_{L^2(\Omega)}\d{t}
+ \langle m(t=0) - m_T,\lambda(t=0)\rangle_{L^2(\Omega)}
- \langle \idiv \vect{v} - w,p\rangle_{L^2(\Omega)}
\end{align}

\noindent with Lagrange multipliers $\lambda:\bar{\Omega}\times[0,1]\rightarrow\ns{R}$ for the transport equation~\eqref{e:varopt:constraint}, and $p:\bar{\Omega}\rightarrow\ns{R}$ for the incompressibility constraint~\eqref{e:varopt:constraint}. Formally, we have to compute variations of $\fs{L}$ with respect to the state, adjoint, and control variables. We will only consider a reduced form (after eliminating the incompressibility constraint) of the optimality system---a system of nonlinear PDEs for $m$, $\lambda$, and $\vect{v}$. Details on how we formally arrive at this reduced from can be found in~\cite{Mang:2015a,Mang:2016a} (see also \secref{s:optimality-cond-dev} in the supplementary materials). We eliminate the state and adjoint variables, and iterate in the control space.

The evaluation of the reduced gradient $\vect{g}$ (the first variation of the Lagrangian $\fun{L}$ in~\eqref{e:lagrangian} with respect to $\vect{v}$) for a candidate $\vect{v}$ requires several steps. We first solve the transport equation~\eqref{e:varopt:constraint} with initial condition $m(\vect{x},t=0) = m_T(\vect{x})$ forward in time to obtain the state variable $m(\vect{x},t)$ for all $t$. Given $m$, we then compute the adjoint variable $\lambda(\vect{x},t)$ for all $t$  by solving the adjoint equation
\begin{subequations}
\label{e:adj}
\begin{align}
-\p_t\lambda - \idiv \lambda\vect{v} & = 0
&&\text{in}\;\;\Omega\times[0,1)
\label{e:adj-transport}
\\
\lambda &=  m_R - m
&&{\rm in} \;\; \Omega \times\{1\}
\label{e:adj-fc}
\end{align}
\end{subequations}

\noindent with periodic boundary conditions on $\p\Omega$ backward in time. Once we have the adjoint and state fields, we can evaluate the expression for the reduced gradient
\begin{equation}
\label{e:reduced-gradient}
\vect{g}(\vect{v}) \defeq
\beta_v\dop{A}[\vect{v}] + \dop{K}\big[\iut{\lambda\igrad m}\;\big].
\end{equation}

The differential operator $\dop{A} = \dop{B}^\ast\dop{B}$ in~\eqref{e:reduced-gradient} corresponds to the first variation of the regularization norm for $\vect{v}$ in~\eqref{e:lagrangian}, e.g., resulting in an elliptic ($\dop{A} = -\ilap$), biharmonic ($\dop{A} = \ilap^2$) , or triharmonic ($\dop{A} = \ilap^3$) control equation for $\vect{v}$, respectively. The operator $\dop{K}$ projects $\vect{v}$ onto the space of incompressible or near-incompressible velocity fields; we have $\dop{K}\defeq\mat{I}-\igrad(\beta_v(\beta_w(-\ilap+\id))^{-1}+\id)^{-1}\ilap^{-1}\idiv$ and $\dop{K}\defeq\id-\igrad\ilap^{-1}\idiv$ for the incompressible case (see \cite{Mang:2015a,Mang:2016a}). If we neglect the incompressibility constraint~\eqref{e:varopt:constraint}, $\dop{K}$ in~\eqref{e:reduced-gradient} is an identity operator. The dependence of $m$ and $\lambda$ on $\vect{v}$ is ``hidden'' in the transport and continuity equations~\eqref{e:varopt:constraint} and~\eqref{e:adj-transport}, respectively.

The first order optimality condition (control or decision equation) requires that $\vect{g}(\vect{v}^\star)=\vect{0}$ for an admissible solution $\vect{v}^\star$ to~\eqref{e:varopt}. Most available registration packages use gradient descent-type optimization schemes to find an optimal point~\cite{Beg:2005a,Hart:2009a,Vialard:2012a}. Newton-type methods are expected to yield better convergence rates~\cite{Nocedal:2006a,Boyd:2004a}. However, if they are implemented naively, they can become computationally prohibitive. The expressions associated with the Newton step of our control problem are formally obtained by computing second variations of the Lagrangian in~\eqref{e:lagrangian}. In \emph{full space methods} we find the Newton updates (i.e., the search direction) for the state, adjoint, and control variables of our control problem simultaneously. That is, we iterate on all variables at once. In \emph{reduced space methods} we only iterate on the control variable $\vect{v}$. Reduced space methods can be obtained from the full space KKT system by block elimination~\cite{Prudencio:2006a,Biros:1999a,Biros:2005a,Biros:2005b}. The associated reduced space Newton system for the incremental control variable $\vect{\tilde{v}}$ (search direction) is given by
\begin{equation}
\label{e:kkt-system}
\dop{H}\vect{\tilde{v}} = - \vect{g},
\end{equation}

\noindent where $\vect{g}$ is the reduced gradient in~\eqref{e:reduced-gradient}. The expression for the (reduced space) Hessian matvec (action of $\dop{H}$ on a vector $\vect{\tilde{v}}$) in~\eqref{e:kkt-system} is given by
\begin{equation}
\label{e:hessian-matvec}
\dop{H}[\vect{\tilde{v}}](\vect{v})
\defeq
\underbrace{
\vphantom{\int_0^1}
\beta_v\dop{A}[\vect{\tilde{v}}]
}_{\eqdef\dop{H}_{\text{reg}}[\tilde{v}]}
+
\underbrace{
\dop{K}\big[
\iut{\tilde{\lambda} \igrad m
+ \lambda \igrad \tilde{m}}
\;\big]}_{\eqdef\dop{H}_{\text{data}}[\tilde{v}](\vect{v})}
=\dop{H}_{\text{reg}}[\tilde{v}] +  \dop{H}_{\text{data}}[\tilde{v}](\vect{v}).
\end{equation}

\noindent We use the notation $\dop{H}[\vect{\tilde{v}}](\vect{v})$ to indicate that the Hessian matvec in~\eqref{e:hessian-matvec} is a function of $\vect{v}$ through a set of PDEs for $m(\vect{x},t)$, $\tilde{m}(\vect{x},t)$, $\lambda(\vect{x},t)$, and $\tilde{\lambda}(\vect{x},t)$. The space-time fields $m$ and $\lambda$ are found during the evaluation of~\eqref{e:varopt:objective} and~\eqref{e:reduced-gradient} for a candidate $\vect{v}$ as described above. What is missing to be able to evaluate~\eqref{e:hessian-matvec} are the fields $\tilde{m}(\vect{x},t)$ and $\tilde{\lambda}(\vect{x},t)$. Given $m(\vect{x},t)$, $\vect{v}(\vect{x})$, and $\vect{\tilde{v}}(\vect{x})$, we find $\tilde{m}(\vect{x},t)$ by solving
\begin{subequations}
\label{e:inc-state}
\begin{align}
\p_t \tilde{m} + \igrad \tilde{m} \cdot \vect{v}
+ \igrad m \cdot \vect{\tilde{v}} & = 0
&&\text{in}\;\; \Omega \times (0,1],
\label{e:inc-state-transport}\\
\tilde{m} &= 0
&&\text{in}\;\;\Omega\times\{0\},
\end{align}
\end{subequations}

\noindent forward in time. Now, given $\tilde{m}(\vect{x},t=1)$, $\lambda(\vect{x},t)$, $\vect{v}(\vect{x})$, and $\vect{\tilde{v}}(\vect{x})$ we solve
\begin{subequations}
\label{e:inc-adj}
\begin{align}
-\p_t \tilde{\lambda} - \idiv(\tilde{\lambda}\vect{v}
+ \lambda\vect{\tilde{v}}) & = 0
&&\text{in}\;\; \Omega \times [0,1),
\label{e:inc-adj-transport}
\\
\tilde{\lambda} &= -\tilde{m}
&&\text{in}\;\;\Omega\times\{1\},
\end{align}
\end{subequations}

\noindent for $\tilde{\lambda}(\vect{x},t)$ backward in time. 

\subsection{Numerics}
\label{s:numerics}

In the following, we describe our distributed-memory solver for 3D diffeomorphic image registration problems.

\subsubsection{Discretization}
\label{s:discretization}

We discretize in space on a regular grid $\Omega^h\in\ns{R}^{3,n_1,n_2,n_3}$ with grid points $\di{x}_{\vect{k}} \defeq 2\pi\vect{k}\oslash\vect{n}_x$, $\vect{k}=(k_1,k_2,k_3)^\T\in\ns{N}^3$, $-n_i/2 + 1 \leq k_i \leq n_i/2$, $i=1,2,3$, $\vect{n}_x\defeq (n_1,n_2,n_3)^\T \in\ns{N}^3$ and periodic boundary conditions; $\oslash$ denotes the Hadamard division. In the continuum, we model images as compactly supported (periodic), smooth functions. We apply Gaussian smoothing (in the spectral domain) with a bandwidth of $\vect{h}_x = (h_1,h_2,h_3)^\T\in\ns{R}^3$ and mollify the discrete data to meet these requirements. We rescale the images to an intensity range of $[0,1]\subset\ns{R}$ prior to registration. We use a trapezoidal rule for numerical quadrature and a spectral projection scheme for all spatial operations. The mapping between spectral and spatial domain is done using forward and inverse FFTs~\cite{Gholami:2016b}. All spatial derivatives are computed in the spectral domain; we first take the FFT, then apply the appropriate weights to the spectral coefficients, and then take the inverse FFT. This scheme allows us to efficiently and accurately apply differential operators and their inverses. Consequently, the main cost of our scheme is the solution of the transport equations~\eqref{e:varopt:constraint}, \eqref{e:adj-transport}, \eqref{e:inc-state-transport}, and~\eqref{e:inc-adj-transport}, and not the inversion of differential (e.g., elliptic or biharmonic) operators. We use a nodal discretization in time, which results in $n_t+1$ space-time fields for which we need to solve. We use a fully explicit, unconditionally stable semi-Lagrangian scheme~\cite{Staniforth:1991a,Falcone:1998a} to solve the transport equations that appear in our formulation (\eqref{e:varopt:constraint},~\eqref{e:adj-transport},~\eqref{e:inc-state-transport}, and~\eqref{e:inc-adj-transport}). This allows us to keep $n_t$ small (we found empirically that $n_t=4$ yields a good compromise between runtime and numerical accuracy). The time integration steps in our semi-Lagrangian scheme are implemented using a fully explicit 2nd order Runge--Kutta scheme. Interpolations are carried out using third-degree polynomials. Details for our semi-Lagrangian scheme can be found in~\cite{Mang:2016c,Mang:2017b,Gholami:2017a}.

\subsubsection{Newton--Krylov Solver}
\label{e:newton-krylov-solver}

A prototype implementation of our Newton--Krylov solver is described in~\cite{Mang:2015a,Mang:2017b}. We have already mentioned in~\secref{s:optsys-and-newtonstep} that we use a reduced space method. That is, we only iterate on the reduced space for the control variable $\di{v}\in\ns{R}^n$, $n=3n_1n_2n_3$. We globalize our method using an Armijo linesearch, resulting in the iterative scheme
\begin{equation}
\label{e:newton-step}
\di{v}_{k+1} = \di{v}_k + \alpha_k \di{\tilde{v}}_k, \qquad \di{H}_k \di{\tilde{v}}_k = - \di{g}_k,
\end{equation}

\noindent with iteration index $k$, step length $\alpha_k \geq 0$, iterate $\di{v}_k\in\ns{R}^n$, search direction $\di{\tilde{v}}_k\in\ns{R}^n$, reduced gradient $\di{g}_k\in\ns{R}^n$ (see \eqref{e:reduced-gradient} for the continuous equivalent), and reduced space Hessian $\di{H}_k \in\ns{R}^{n,n}$, where \[\di{H}_k = \di{H}_{\text{reg}} + \di{H}_{\text{data},k}.\] (See \eqref{e:hessian-matvec} for an expression for the Hessian matvec in the continuous setting.) We refer to the steps for updating $\di{v}_k$ as outer iterations and the steps for computing the search direction $\di{\tilde{v}}_k$ as inner iterations.

\begin{algorithm}
\caption{Inexact Newton--Krylov method (outer iterations). We use the relative norm of the reduced gradient (with tolerance $\epsilon_{\text{opt}} > 0$ as stopping criterion.}
\label{a:outer-iteration}
\algadjust
\begin{algorithmic}[1]
\STATE{initial guess $\di{v}_0 \leftarrow \dig{0}$, $k\leftarrow0$}
\STATE{$\di{m}_0 \leftarrow$ solve state equation in~\eqref{e:varopt:constraint} forward in time given $\di{v}_0$}
\STATE{$j_0 \leftarrow$ evaluate objective functional~\eqref{e:varopt:objective} given $\di{m}_0$ and $\di{v}_0$}
\STATE{$\dig{\lambda}_0 \leftarrow$ solve adjoint equation~\eqref{e:adj-transport} backward in time given $\di{v}_0$ and $\di{m}_0$}
\STATE{$\di{g}_0 \leftarrow$ evaluate reduced gradient~\eqref{e:reduced-gradient} given $\di{m}_0$, $\dig{\lambda}_0$ and $\di{v}_0$}
\WHILE{$\|\di{g}_k\|^2_2 > \|\di{g}_0\|_2^2\epsilon_{\text{opt}}$\label{alg:stop-newton}}
    \STATE{$\di{\tilde{v}}_k\leftarrow$ solve $\di{H}_k \di{\tilde{v}}_k = -\di{g}_k$ given $\di{m}_k$, $\dig{\lambda}_k$, $\di{v}_k$, and $\di{g}_k$ (see \algref{a:inner-iteration})}
    \STATE{$\alpha_k \leftarrow$ perform line search on $\di{\tilde{v}}_k$ subject to Armijo condition}
    \STATE{$\di{v}_{k+1} \leftarrow \di{v}_k + \alpha_k\di{\tilde{v}}_k$}
    \STATE{$\di{m}_{k+1} \leftarrow$ solve state equation~\eqref{e:varopt:constraint} forward in time given $\di{v}_{k+1}$}
    \STATE{$j_{k+1} \leftarrow$ evaluate~\eqref{e:varopt:objective} given $\di{m}_{k+1}$ and $\di{v}_{k+1}$}
    \STATE{$\dig{\lambda}_{k+1} \leftarrow$ solve adjoint equation~\eqref{e:adj-transport} backward in time given $\di{v}_{k+1}$ and $\di{m}_{k+1}$}
    \STATE{$\di{g}_{k+1} \leftarrow$ evaluate~\eqref{e:reduced-gradient} given $\di{m}_{k+1}$, $\dig{\lambda}_{k+1}$ and $\di{v}_{k+1}$}
    \STATE{$k \leftarrow k + 1$}
\ENDWHILE
\end{algorithmic}
\end{algorithm}

\begin{algorithm}
\caption{Newton step (inner iterations). We illustrate the solution of the reduced KKT system~\eqref{e:kkt-system} using a PCG method at a given outer iteration $k\in\ns{N}$. We use a superlinear forcing sequence to compute the tolerance $\eta_k$ for the PCG method (inexact solve).}
\label{a:inner-iteration}
\algadjust
\begin{algorithmic}[1]
\STATE{input: $\di{m}_k$, $\dig{\lambda}_k$, $\di{v}_k$, $\di{g}_k$, $\di{g}_0$}
\STATE{set $\epsilon_H \leftarrow\min(0.5,(\|\di{g}_k\|_2/\|\di{g}_0\|_2)^{1/2})$, $\di{\tilde{v}}_0 \leftarrow \dig{0}$, $\di{r}_0 \leftarrow - \di{g}_k$ \label{alg:init}}
\STATE{$\di{z}_0 \leftarrow$ apply preconditioner $\di{M}^{-1}$ to $\di{r}_0$}
\STATE{$\di{s}_0 \leftarrow \di{z}_0$, $\iota\leftarrow 0$}
\WHILE{$\iota < n$}
    \STATE{$\di{\tilde{m}}_\iota \leftarrow$ solve incremental state equation~\eqref{e:inc-state} forward in time given $\di{m}_k$, $\di{v}_k$ and $\di{\tilde{v}}_\iota$\label{alg:inc-state}}
    \STATE{$\dig{\tilde{\lambda}}_\iota \leftarrow$ solve incremental adjoint equation~\eqref{e:inc-adj} backward in time given $\dig{\lambda}_k$, $\di{v}_k$, $\di{\tilde{m}}_\iota$ and $\di{\tilde{v}}_\iota$\label{alg:inc-adjoint}}
    \STATE{$\di{\tilde{s}}_\iota \leftarrow$ apply $\di{H}_\iota$ to $\vect{s}_\iota$ given $\dig{\lambda}_k$, $\di{m}_k$, $\di{\tilde{m}}_\iota$ and $\dig{\tilde{\lambda}}_\iota$ (Hessian matvec; see \eqref{e:hessian-matvec})}
    \STATE{$\kappa_\iota \leftarrow\langle\di{r}_\iota,\di{z}_\iota\rangle/\langle\di{s}_\iota,\di{\tilde{s}}_\iota\rangle$,
            \quad $\di{\tilde{v}}_{\iota+1} \leftarrow \di{\tilde{v}}_\iota + \kappa_\iota\di{s}_\iota$,
            \quad $\di{r}_{\iota+1}\leftarrow\di{r}_\iota-\kappa_\iota\di{\tilde{s}}_\iota$}
    \STATE{\textbf{if} $\|\di{r}_{\iota+1}\|_2 < \epsilon_H$ \textbf{break}\label{alg:forcing}}
    \STATE{$\di{z}_{\iota+1} \leftarrow$ apply preconditioner $\di{M}^{-1}$ to $\di{r}_{\iota+1}$}
    \STATE{$\mu_\iota \leftarrow \langle\di{z}_{\iota+1},\di{r}_{\iota+1}\rangle/\langle\di{z}_\iota,\di{r}_\iota\rangle$, \quad $\di{s}_{\iota+1} \leftarrow \di{z}_{\iota+1} + \mu_\iota\di{s}_\iota$, \quad $\iota \leftarrow \iota + 1$}
\ENDWHILE
\STATE{output: $\di{\tilde{v}}_k \leftarrow\di{\tilde{v}}_{\iota+1}$}
\end{algorithmic}
\end{algorithm}

In what follows, we drop the dependence on the (outer) iteration index $k$ for notational convenience. The data term $\di{H}_{\text{data}}$ of the reduced space Hessian $\di{H}$ in~\eqref{e:newton-step} involves inverses of the state and adjoint operators (a consequence of the block elimination in reduced space methods; see~\secref{s:optsys-and-newtonstep}). This makes $\di{H}$ a nonlocal, dense operator that is too large to be computed and stored. (We have seen in~\secref{s:optsys-and-newtonstep} that each matvec given by~\eqref{e:hessian-matvec} requires the solution of~\eqref{e:inc-state} forward in time and~\eqref{e:inc-adj} backward in time; see also lines~\ref{alg:inc-state} and \ref{alg:inc-adjoint} in \algref{a:inner-iteration}. So, to form $\di{H}$ we require a total of $2n$ PDE solves per outer iteration $k$.) Consequently, direct methods to solve the linear system in~\eqref{e:newton-step} are not applicable. We use iterative, matrix-free Krylov subspace methods instead. They only require an expression for the action of $\di{H}$ on a vector, which is precisely what we are given in~\eqref{e:hessian-matvec}. We use a PCG method~\cite{Hestenes:1952a} under the assumption that $\di{H}$ is a symmetric, positive (semi-)definite operator. To reduce computational costs, we do not solve the linear system in~\eqref{e:newton-step} exactly. Instead, we use a tolerance $\epsilon_H>0$ that is proportional to the norm of $\di{g}$ (see lines~\ref{alg:init} and~\ref{alg:forcing} in \algref{a:inner-iteration}; details can be found in~\cite{Dembo:1982a,Eisenstat:1996a} and~\cite[p.~166ff]{Nocedal:2006a}). We summarize the steps for the outer and inner iterations of our Newton--Krylov solver in~\algref{a:outer-iteration} and~\algref{a:inner-iteration}, respectively.

Since we are solving a non-convex problem it is not guaranteed that the Hessian $\di{H}$ is positive definite. As a remedy, we use a Gauss--Newton approximation $\di{\tilde{H}}$ to $\di{H}$; doing so guarantees that $\di{\tilde{H}} \succeq 0$ far away from the (local) optimum. This corresponds to dropping all terms in~\eqref{e:hessian-matvec} and~\eqref{e:inc-adj} that involve the adjoint variable $\lambda$. We expect the rate of convergence of our solver to drop from quadratic to superlinear. As $\lambda$ tends to zero (i.e., the mismatch goes to zero), we recover quadratic convergence. We terminate the inversion if the $\ell^2$-norm of the gradient in~\eqref{e:reduced-gradient} is reduced by a factor of $\epsilon_g>0$, i.e., if $\|\di{g}_k\|_2^2 <  \epsilon_{\text{opt}}\|\di{g}_0\|_2^2$, where $\di{g}_k\in\ns{R}^n$ is the gradient at (outer) iteration $k\in\ns{N}_0$ and $\di{g}_0$ is the gradient for the initial guess $\di{v}_0 = \vect{0}$ (see line~\ref{alg:stop-newton} in \algref{a:outer-iteration}). In most of our experiments, we use $\epsilon_g=\num{5E-2}$. We also provide an option to set a lower bound for the absolute $\ell^2$-norm of the gradient (the default value is $\num{1E-6}$). \claire{} also features other stopping criteria discussed in~\cite{Modersitzki:2009a,Gill:1981a,Mang:2015a} (not considered in this work).

\subsubsection{Preconditioners for Reduced Space Hessian}

We have seen that we need to solve two PDEs every time $\di{H}$ is applied to a vector. These PDE solves are the most expensive part of our solver. Consequently, we have to keep the number of Hessian matvecs small for our solver to be efficient. This necessitates the design of an effective preconditioner $\di{M}$. The speed of convergence of the linear solver used to compute the search direction $\di{\tilde{v}}$ in~\eqref{e:newton-step} depends on the distance of $\di{M}^{-1}\di{H}$ from identity; ideally, the spectrum of $\di{M}^{-1}\di{H}$ is clustered around one. We cannot form and store $\di{H}$ (too expensive). Moreover, we know that large eigenvalues of $\di{H}$ are associated with smooth eigenvectors~\cite{Mang:2015a}. Consequently, standard preconditioners for linear systems are not applicable. In our former work, we have considered two matrix-free preconditioners. Our first preconditioner is based on the (exact) inverse of the regularization operator $\di{H}_{\text{reg}}$; the regularization preconditioned Hessian is given by $\mat{I} + \di{H}_{\text{reg}}^{-1} \di{H}_{\text{data}}$. This is a common choice in PDE-constrained optimization problems~\cite{Alexanderian:2016a,BuiThanh:2013a}. We used this preconditioner in~\cite{Mang:2015a,Mang:2016a,Mang:2016c,Mang:2017a,Mang:2018a}.

\begin{remark}
$\di{H}_{\text{reg}}$ is a discrete representation of the regularization operator. The computational costs for inverting and applying this operator are negligible (two FFTs and a diagonal scaling). Notice that the operator $\di{H}_{\text{reg}}$ is singular if we consider a seminorm as regularization model in~\eqref{e:varopt}. We bypass this problem by setting the zero singular values of the regularization operator to one before computing the inverse.
\end{remark}

The second preconditioner uses an inexact inverse of a coarse grid approximation to the Hessian $\di{H}$. This preconditioner was proposed and tested in~\cite{Mang:2017b} for the 2D case. A similar preconditioner has been developed in~\cite{Adavani:2008a,Biros:2008a}. It is based on the conceptual idea that we can decompose the reduced Hessian $\di{H}$ into two operators $\di{H}_L$ and $\di{H}_H$ that act on the high and low frequency parts of a given vector $\di{\tilde{v}}$, respectively~\cite{Adavani:2008a,Biros:2008a,Giraud:2006a,Kaltenbacher:2003a,Kaltenbacher:2001a,King:1990a}. We denote the operators that project on the low and high frequency subspaces by $\di{F}_L : \ns{R}^n \to \ns{R}^n$ and $\di{F}_H : \ns{R}^n \to \ns{R}^n$, respectively. With $\di{F}_H + \di{F}_L = \di{I}$, we have
\[
\di{H}\di{e}_k
= (\di{F}_H + \di{F}_L)\di{H}(\di{F}_H + \di{F}_L)\di{e}_k
= \di{F}_H\di{H}\di{F}_H\di{e}_k + \di{F}_L\di{H}\di{F}_L\di{e}_k,
\]
\noindent under the assumption that the unit vector $\di{e}_k\in \ns{R}^n$, $(\di{e}_k)_i = 1$ if $k=i$ and $(\di{e}_k)_i = 0$ for $i\not=k$, $i,k=1,\ldots,n$, is an eigenvector of $\di{H}$ so that $(\di{F}_L\di{H}\di{F}_H)\di{e}_k = (\di{F}_H\di{H}\di{F}_L)\di{e}_k = \di{0}$. We note that this assumption will not hold in general. However, since we are only interested in developing a preconditioner, an approximate decomposition of the solution of the reduced space system is acceptable. Using this approximation we can represent the solution of $\di{H} \di{\tilde{v}} = -\di{g}$ as $\di{\tilde{v}} = \di{\tilde{v}}_L + \di{\tilde{v}}_H$, where $\di{\tilde{v}}_L$ and $\di{\tilde{v}}_H$ are found by solving \[\di{H}_L\di{\tilde{v}}_L = (\di{F}_L\di{H}\di{F}_L)\di{\tilde{v}}_L = -\di{F}_L\di{g} \quad\text{and}\quad \di{H}_H\di{\tilde{v}}_H = (\di{F}_H\di{H}\di{F}_H)\di{\tilde{v}}_H = -\di{F}_H\di{g},\] respectively.

We discuss how we use this decomposition to design a preconditioner, next. Let $\di{s}\in\ns{R}^n$ denote the vector we apply our preconditioner to. Since we use an approximation of the inverse of $\di{H}$, we have to design a scheme for approximately solving $\di{H} \di{u} = \di{s}$. We find the smooth part of $\di{u}$ by (iteratively) solving
\begin{equation}
\label{e:coarse-grid-hessian}
\di{\bar{H}}\di{\bar{u}}_L = \di{Q}_R\di{F}_L\di{s},
\end{equation}

\noindent where $\bar{\di{H}}\in\ns{R}^{c,c}$ and $\di{\bar{u}}_L\in\ns{R}^c$ represent coarse grid approximations of $\di{H}_L$ and $\di{u}_L$, respectively, and $\di{Q}_R : \ns{R}^n \to \ns{R}^{c}$ is a restriction operator. We do not iterate on the oscillatory components of $\di{s}$ (i.e., we replace $\di{H}_H$ by $\di{I}$). The solution $\di{u}$ of $\di{H} \di{u} = \di{s}$ is given by $\di{u} = \di{u}_L + \di{u}_H \approx \di{Q}_P\di{F}_L\di{\bar{u}}_L + \di{F}_H\di{s}$, where $\di{\bar{u}}_L \approx \di{\bar{H}}^{-1}\di{Q}_R\di{F}_L\di{s}$ and $\di{Q}_P : \ns{R}^{c} \to \ns{R}^n$ is a prolongation operator. We use spectral restriction and prolongation operators. The projection operators $\di{F}_L$ and $\di{F}_H$ are implemented as cut-off filters in the frequency domain.

An important aspect of our approach is that we do not apply our two-level preconditioner to the original Hessian $\di{H}$. Since we can invert $\di{H}_{\text{reg}}\succeq 0$ explicitly, we consider the (symmetric) regularization split-preconditioned system
$
(\di{I} + \di{H}_{\text{reg}}^{\nicefrac{-1}{2}} \di{H}_{\text{data}}\di{H}_{\text{reg}}^{\nicefrac{-1}{2}}) \di{w} = - \di{H}_{\text{reg}}^{\nicefrac{-1}{2}} \di{g}
$
instead, where $\di{w} \defeq \di{H}_{\text{reg}}^{\nicefrac{1}{2}} \di{\tilde{v}}$. Notice that the inverse of $\di{H}_{\text{reg}}$ acts as a smoother. This allows us to get away with not treating high-frequency errors in our scheme. Our approach can be interpreted as an approximate two-level multigrid V-cycle with an explicit (algebraic) smoother given by $\di{H}_{\text{reg}}^{\nicefrac{-1}{2}}$.

The final questions are how to discretize and solve~\eqref{e:coarse-grid-hessian}. We can use a Galerkin or a direct (non-Galerkin) discretization to implement the coarse grid operator $\di{\bar{H}}$. Using the fact that $\di{Q}_R$ and $\di{Q}_P$ are adjoint operators, the Galerkin discretization is formally given by $\di{\bar{H}} = \di{Q}_R\di{H}\di{Q}_P$~\cite[p.~75]{Briggs:2000a}. The drawback of using a Galerkin operator is that every matvec requires the solution of the incremental forward and adjoint equations on the fine grid. This is different if we directly discretize the matvec on the coarse grid. To save computational costs, we opt for this approach. For the iterative solver to approximately invert $\di{\bar{H}}$ we have tested several variants, all of which are available in \claire. We can use a nested PCG. This requires a tolerance $\epsilon_M > 0$ for the nested solver for the preconditioner that is only a fraction of the tolerance used to solve for the Newton step $\di{\tilde{v}}$ on the fine grid, i.e., $\epsilon_M = \kappa\epsilon_H$ with $\kappa\in(0,1)$. This is due to the fact that Krylov subspace methods are nonlinear operators. We refer to this solver as PCG($\kappa$). Another possibility is to use a semi-iterative Chebyshev method~\cite{Gutknecht:2002a} with a predefined number of iterations $k > 0$; this results in a fixed linear operator for a particular choice of eigenvalue bounds~\cite{Golub:1961a}. The eigenvalue bounds can be estimated using a Lanczos method. We refer to this strategy as CHEB($k$). If we would like to use PCG with a fixed number of iterations as a nested solver, we can also replace the solver for the Newton step with a flexible Krylov subspace method~\cite{Axelsson:1991a,Notay:2000a}. We observed that the performance of this approach deteriorates significantly as we reduce the regularization parameter. We disregard this approach.

\section{Implementation and Software Aspects}
\label{s:software}

We make \claire{} available under GPL license. \claire{} is written in C/C\texttt{++} and implements data parallelism via MPI. The source code can be downloaded from the \sw{github} repository~\cite{claire-web} at
\begin{center}
\url{https://github.com/andreasmang/claire}.
\end{center}

The URL for the deployment page of \claire{} is \url{https://andreasmang.github.io/claire}. Here, one can find a detailed documentation as well as use cases for \claire{}. In what follows, we \bipa\item describe implementation aspects, \item list features implemented in \claire{}, and \item provide information relevant to potential users of \claire{}\eipa. It is important to note that we will not be able to cover all implementation aspects, and we are continuously making improvements to our software. We refer the reader to the deployment page for updates and detailed information on how to compile, execute, and run \claire{} on various systems.

As we have mentioned above, \claire{} is written in C\texttt{++}. The main functionalities of \claire{} are implemented in \texttt{CLAIRE.cpp}. Different formulations are implemented using derived classes. The distance measures and regularization operators supported by \claire{} are, like most of the building blocks of \claire{}, implemented through classes (again, using inheritance). We provide interfaces to the main \petsc{} functionalities through functions.

\subsection{Executables}

\claire{} has two main executables, \texttt{claire} and \texttt{clairetools}. The registration solver can be executed with the \texttt{claire} executable. The \texttt{clairetools} executable serves as a postprocessing tool that allows users to, e.g., compute deformation measures (examples include the deformation map $\vect{y}$, the determinant of the deformation gradient, or a RAVENS map), or transport images or label maps for the evaluation of registration performance. We will keep adding features to these executables in future releases. Both executables provide a help message that briefly explains to users how to control the behavior, how to set parameters, and what features are provided. To access this help message, the user can simply execute the binaries without any parameters or add a \texttt{-help} flag to the executable (i.e., for instance execute \texttt{claire -help} from the command line window). The main output of \texttt{claire} is the computed velocity field. These fields can subsequently be used within \texttt{clairetools} to compute additional outputs. We explain the most common options for both executables in greater detail on the deployment page / in the \texttt{README} files for the repository.

\subsection{External Dependencies and IO}

\claire{} depends on four main software packages. We use the \petsc{} library~\cite{Balay:2016a,petsc-web} for linear algebra, and \petsc{}'s \tao{} package~\cite{Munson:2017a,petsc-web} for numerical optimization (\tao{} is included in \petsc). We use the \accfft{} package~\cite{accfft-web,Gholami:2016b}---a parallel, open-source FFT library for CPU/GPU architectures developed in our group---to apply spectral operators. \accfft{} requires \fftw~\cite{Frigo:2005a,fftw-web}. We use \sw{niftilib}~\cite{niftilib-web} for IO. As such, \claire{} currently supports IO of (uncompressed and compressed in \texttt{gzip} format) files in \texttt{nifti-1} (\texttt{*.nii} or \texttt{*.nii.gz}) and Analyze 7.5 (\texttt{*.hdr} and \texttt{*.img}/\texttt{*.img.gz}) format. The default output format of \claire{} is in \texttt{*.nii.gz}. We optionally support the \sw{PnetCDF} format (\texttt{*.nc})~\cite{Li:2003a,pnetcdf-web} for IO in parallel. The revision and version numbers for these libraries used in our experiments can be found in the references.

\subsection{Compilation and Installation}

Our solver supports single and double precision. (The precision is handed down from the \petsc{} library.) Our current software uses \texttt{make} for compilation. We provide scripts in the repository to download and compile the external libraries mentioned above using default settings that have worked most consistently on the systems on which we have executed \claire{}. Switches for controlling the precision are provided in the \texttt{makefile}. The user needs to compile \petsc{} and \fftw{} in single precision to be able to run \claire{} in single precision. We have compiled, tested, and executed \claire{} on HPC systems at TACC~\cite{tacc-web} (Stampede, Stampede 2, Lonestar 5, and Maverick), at HLRS (Hazelhen/CRAY XC40)~\cite{hlrs-web} and at RCDC~\cite{cacds-web} (Opuntia and Sabine). Specifications of some of these systems can be found in \secref{s:compute-systems}. While we recommend the execution of \claire{} on multicore systems (to reduce the runtime), it is not a prerequisite to have access to HPC systems. \claire{} has been successfully executed on personal computers and local compute servers with no internode communication. Large-scale systems are only required to significantly reduce the runtime or when considering large-scale applications (image sizes of $512^3$ and beyond). We provide additional help for compilation and installation of \claire{} in the repository.

\subsection{Parallel Algorithms and Computational Kernels}

The main computational kernels of \claire{} are FFTs (spectral methods) and scattered data interpolation operations (semi-Lagrangian solver; see~\cite{Mang:2016c,Mang:2017b,Gholami:2017a} for details). We use the \accfft{} package~\cite{accfft-web,Gholami:2016b} to perform spectral operations (a software package developed by our group). This package dictates the data layout on multicore systems: We partition the data based on a pencil decomposition for 3D FFTs~\cite{Grama:2003a,Czechowski:2012a}. Let $n_p=p_1p_2$ denote the number of MPI tasks. Then each MPI task gets $(\nicefrac{n_1}{p_1})\times(\nicefrac{n_2}{p_2})\times n_3$ grid points. That is, we partition the domain $\dig{\Omega}$ of size $3\times n_1\times n_2 \times n_3$ along the $x_1$- and $x_2$-axes into subdomains $\dig{\Omega}_i$, $i=1,\ldots,n_p$, of size $3\times (\nicefrac{n_1}{p_1})\times (\nicefrac{n_2}{p_2})\times n_3$. There is no partitioning in time.

The scalability of the 3D FFT is well explored~\cite{Grama:2003a,Czechowski:2012a,Gholami:2016b}. We refer the reader to~\cite{Gholami:2016b,Mang:2016a} for scalability results for \accfft. If we assume that the number of grid points $n_i$, $i=1,2,3$, is equal along each spatial direction, i.e., $\tilde{n} = n_1 = n_2 = n_3$, each 3D FFT requires $\dop{O}(\nicefrac{3\tilde{n}\log(\tilde{n})}{2n_p})$ computations and $\dop{O}(2\sqrt{\vphantom{b}n_p}t_s+(\nicefrac{2\tilde{n}^3}{n_p})t_w)$ communications, where $t_s> 0$ is the startup time for the data transfer and $t_w>0$ represents the per-word transfer time~\cite{Grama:2003a}.

\begin{figure}
\centering
\includegraphics[totalheight=3cm]
{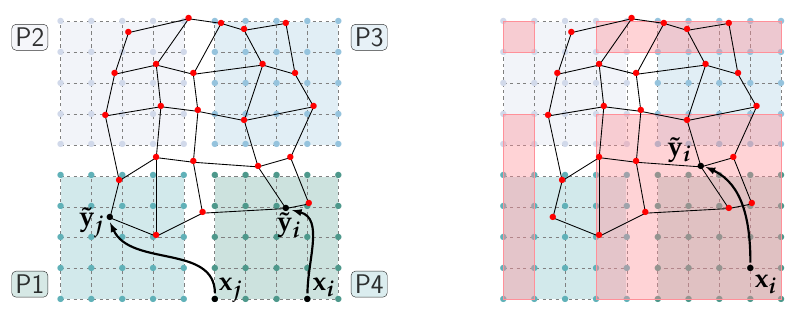}
\caption{2D illustration of the data layout and the communication steps for the evaluation of the interpolation kernel. The original grid at timepoint $t^{k+1}$ is distributed across $n_p=4$ processors P$i$, $i=1,2,3,4$. To solve the transport problem using a semi-Lagrangian scheme, we have to trace a characteristic for each grid point $\di{x}_{\vect{l}}$ backward in time (see~\cite{Mang:2016c,Mang:2017b,Gholami:2017a} for details). This requires a scattered data interpolation step. The deformed configuration of the grid (i.e., the departure points) originally owned by P4 (red points) are displayed in overlay. We illustrate three scenarios: The departure point is located (i) on P4 (left; $\di{x}_{\vect{i}}\to \di{\tilde{y}}_{\vect{i}}$), (ii) on a different processor P1 (left; $\di{x}_{\vect{j}}\to\di{\tilde{y}}_{\vect{j}}$), and (iii) between processors P3 and P4 (right). For the first case, no communication is required. The second case requires the communication of $\di{\tilde{y}}_{\vect{j}}$ to P1, and the communication of the interpolation result back to P4. For the third case, we add a ghost layer with a size equal to the support of the interpolation kernel (4 grid points in our case) to each processor; the evaluation of the interpolation happens on the same processor (like in the first case). Notice that the communication of the departure points (for the forward and backward characteristics) needs to be performed only once per Newton iteration, since the velocity field is stationary.\label{f:sl-scheme}}
\end{figure}

The parallel implementation of our interpolation kernel is introduced in~\cite{Mang:2016c} and improved in~\cite{Gholami:2017a}. We use a tricubic interpolation model to evaluate off-grid points in our semi-Lagrangian scheme (see \cite{Mang:2016c,Mang:2017b} for a detailed description of our solver). The polynomial is implemented in Lagrange form. The evaluation of the interpolation kernel requires the computation of 12 basis polynomials. The local support of the cubic basis is $4^3$ grid points. Overall, this results in a complexity of $\dop{O}(256\nicefrac{\tilde{n}^3}{n_p})$ computations. We have implemented an SIMD vectorization based on advanced vector extensions (AVX2) for Haswell architectures for the evaluation of the interpolation kernel (available for single precision only). Compared to our initial work in~\cite{Mang:2016c} our method is now bound by communication time instead of time spent in the interpolation. The communication costs are more difficult to estimate; they not only depend on the data layout but also on the characteristics obtained for a given velocity field. If a departure point is owned by the current processor, we require no communication. If the values for a departure point are owned by another processor/MPI task (the \emph{worker}), we communicate the coordinates from the \emph{owner} to the worker. We then evaluate the interpolation model on the worker and communicate the result back to the owner. This results in a communication cost of $4t_w$ per off-grid point not owned by a processor. To evaluate the interpolation model at off-grid points not owned by any MPI task (i.e., located in between the subdomains $\dig{\Omega}_i$), we add a layer of four ghost points (scalar values to be interpolated; see~\figref{f:sl-scheme}, right). This results in an additional communication cost of $n_s(2n_3(\nicefrac{n_1}{p_1} + \nicefrac{n_2}{p_2})t_w) + 4t_s$ for each MPI task for the four face neighbors, where $n_s$ is the size of layer for the ghost points (in our case four). The communication with the four corner neighbors can be combined with the messages of the edge neighbors, by appropriate ordering of the messages. Notice that the communication of the departure points (for the forward and backward characteristics) needs to be performed only once per Newton iteration, since our velocity field is stationary. We perform this communication when we evaluate the forward and the adjoint operators, i.e., during the evaluation of the objective functional and the reduced gradient.

\subsection{Memory Requirements}

In our most recent implementation, we have reduced the memory footprint for the Gauss--Newton approximation; we only store the time history of the state and incremental state variables. This is accomplished by evaluating the time integrals that appear in the reduced gradient in~\eqref{e:reduced-gradient} and the Hessian matvec in~\eqref{e:hessian-matvec} simultaneously with the time integration of the adjoint and incremental adjoint equations~\eqref{e:adj} and~\eqref{e:inc-adj}, respectively. With this we can reduce the memory pressure from $\dop{O}((2n_t +8)n_1n_2n_3)$ (full Newton) to $\dop{O}((n_t + 7)n_1n_2n_3)$ (Gauss--Newton) for the gradient (see \eqref{e:reduced-gradient}) and from $\dop{O}((4n_t + 13)n_1n_2n_3)$ (full Newton) to $\dop{O}((n_t + 10)n_1n_2n_3)$ (Gauss--Newton) for the Hessian matvec (see \eqref{e:hessian-matvec}), respectively. Notice that we require $0.5\times$ the memory of the Hessian matvec, if we consider the two-level preconditioner. The spectral preconditioner does not add to the memory pressure.

\subsection{Additional Software Features}

We provide schemes for automatically selecting an adequate regularization parameter. This a topic of research by itself~\cite{Hansen:1998a,Haber:2000b}. Disregarding theoretical requirements~\cite{Beg:2005a,Dupuis:1998a,Trouve:1998a}, one in practice typically selects an adequate regularization norm based on application requirements. From a practical point of view we are interested in computing velocities for which the determinant of the deformation map does not change sign/is strictly positive for every point inside the domain. This guarantees that the transformation is locally diffeomorphic (subject to numerical accuracy). Consequently, we determine the regularization parameter $\beta_v$ for the Sobolev norm for the velocity based on a binary search (this strategy was originally proposed in~\cite{Mang:2015a}; a similar strategy is described in~\cite{Haber:2000b}). We control the search based on a bound for the determinant of the deformation gradient. That is, we choose $\beta_v$ so that the determinant of the deformation gradient is bounded below by $\epsilon_J$ and bounded above by $1/\epsilon_J$, where $\epsilon_J \in (0,1)$ is a user defined parameter. This search is expensive, since it requires a repeated solution of the inverse problem. (For each trial $\beta_v$ we iterate until we meet the convergence criteria for our Newton solver and then use the obtained velocity as an initial guess for the next $\beta_v$.) We assume that, once we have found an adequate $\beta_v$, we can use this parameter for similar registration problems. Such cohort studies are quite typical in medical imaging.

\claire{} features several well established schemes to accelerate the rate of convergence and reduce the likelihood to get trapped in local minima. The user can choose between \bipa\item parameter continuation in $\beta_v$ (starting with a default value of $\beta_v=1$ we reduce $\beta_v$ until we reach a user defined parameter $\beta_v^\star$; we found this scheme to perform best), \item grid continuation, i.e., a coarse-to-fine multi-resolution scheme with a smoothing of $\sigma = 1$ voxels (consequently, the standard deviation increases for coarser grids), and \item scale continuation using a scale-space representation of the image data (again, coarse-to-fine).\eipa

We summarize the critical parameters of \claire{} in~\tabref{t:claire-parameters}.

\begin{table}
\centering\scriptsize
\caption{Parameters available in \claire{} (there are more, but these are the critical ones).\label{t:claire-parameters}}
\begin{tabular}{llll}\toprule
variable     & meaning                                  & suggested value                  & determined automatically \\
\midrule
$\beta_v$    & regularization parameter for $\vect{v}$  & ---                              & yes                      \\
$\beta_w$    & regularization parameter for $w$         & \num{1e-4}                       & no                       \\
$\epsilon_g$ & relative tolerance for gradient          & \num{5e-2}                       & no                       \\
$n_t$        & number of time steps                     & 4                                & no                       \\
$\epsilon_j$ & bound for $\det\igrad\vect{y}^{-1}$      & 0.25 ($H^1$-div) or 0.1 ($H^2$)  & no                       \\
\bottomrule
\end{tabular}
\end{table}

\section{Experiments}
\label{s:experiments}

\begin{figure*}
\centering
\includegraphics[width=0.8\textwidth]
{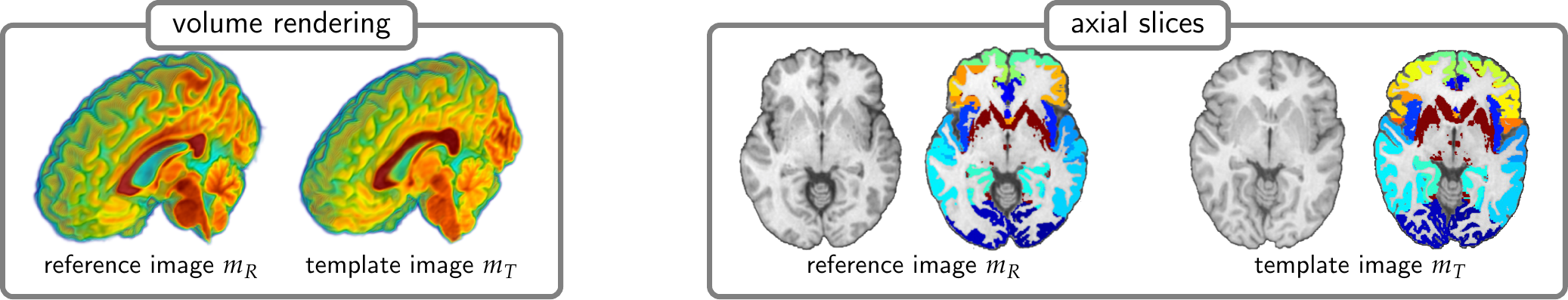}
\caption{Illustration of exemplary images from the \nirep{} data~\cite{Christensen:2006a}. Left: Volume rendering of an exemplary reference image $m_R(\vect{x})$ (dataset {\tt na01}) and an exemplary template image $m_T(\vect{x})$ (dataset {\tt na03}), respectively. Right: Axial slice for these datasets together with label maps associated with these data. Each color corresponds to one of the 32 individual anatomical gray matter regions that serve as a ground truth to evaluate the registration performance.}
\label{f:nirep-data}
\end{figure*}

We evaluate the registration accuracy for 16 segmented MRI brain volumes~\cite{Christensen:2006a}. Details on the considered datasets can be found in~\secref{s:data}. We showcase two exemplary datasets in~\figref{f:nirep-data}. Notice that these datasets have been rigidly preregistered. We directly apply our method to this data (without an additional affine preregistration step). The runs were executed on the RCDC's Opuntia server or on TACC's Lonestar 5 system. The specs of these systems can be found below. Notice that we accompany this document with \emph{supplementary materials} that provide more detailed results for some of the experiments conducted in this study.

For \claire{} we consider two models: \bipa \item \textbf{\textit{$H^1$-div regularization}}: $H^1$-seminorm for the regularization model for the velocity field (controlled by $\beta_v$; $\dop{B} = \igrad$) in combination with a penalty for the divergence of $\vect{v}$ (controlled by $\beta_w$, which is fixed to $\beta_w = \num{1e-4}$). \item \textbf{\textit{$H^2$ regularization}}: $H^2$-seminorm for the regularization model for the velocity field (controlled by $\beta_v$; $\dop{B} = -\ilap$). No penalty for the divergence of $\vect{v}$ is added\eipa.

\subsection{Setup, Implementation, and Hardware}
\label{s:compute-systems}

We execute the runs on RCDC's Opuntia system (Intel ten-core Xeon E5-2680v2 at \SI{2.8}{\giga\hertz} with \SI{64}{\giga\byte} memory; 2 sockets for a total of 20 cores~\cite{cacds-web}) and TACC's Lonestar 5 system (2-socket Xeon E5-2690 v3 (Haswell) with 12 cores/socket, 64 GB memory per node~\cite{tacc-web}). Our code is written in C++ and uses MPI for parallelism. It is compiled with the default Intel compilers available on these systems at the time (Lonestar 5: Intel 16.0.1 and Cray MPICH 7.3.0; Opuntia: Intel PSXE 2016, Intel ICS 2016, and Intel MPI 5.1.1). We use \claire{} commit \texttt{v0.07-131-gbb7619e} to perform the experiments. For the software packages/libraries used in combination with \claire{}, we refer the reader to \secref{s:software}. The versions of the libraries used for our runs can be found in the references.

\subsection{Real and Synthetic Data}
\label{s:data}

We report results for the \nirep{} (\iquote{Non-Rigid Image Registration Evaluation Project}) data~\cite{Christensen:2006a}. This repository contains 16 rigidly aligned T1-weighted MRI brain datasets ({\tt na01}--{\tt na16}; image size: $256\times300\times256$ voxels) of different individuals. Each dataset comes with 32 labels of anatomical gray matter regions. (Additional information on the datasets, the imaging protocol, and the preprocessing can be found in~\cite{Christensen:2006a}.) We illustrate an exemplary dataset in \figref{f:nirep-data}. The initial Dice score (before registration) for the combined label map (i.e., the union of the 32 individual labels) is on average \num{5.181288528e-01} (mean) with a maximum of \num{5.622982959E-01} (dataset \texttt{na08}) and a minimum of \num{4.379558331E-01} (dataset \texttt{na14}). We generate the data for grids not corresponding to the original resolution based on a cubic interpolation scheme.

The scalability runs reported in \secref{s:scalability} are based on synthetic test data. We use a template image $m_T(\vect{x}) = ((\sin x_1)(\sin x_1) + (\sin x_2)(\sin x_2) + (\sin x_3)(\sin x_3))/3$. The reference image $m_R(\vect{x})$ is computed by solving the forward problem for a predefined velocity field $\vect{v}^\star(\vect{x}) = (v_1^\star(\vect{x}), v_2^\star(\vect{x}),v_3^\star(\vect{x}))^\T$, where $v_1^\star(\vect{x}) = \sin x_3 \cos x_2 \sin x_2$, $v_2^\star(\vect{x}) = \sin x_1 \cos x_3 \sin x_3$, and $v_3^\star(\vect{x}) = \sin x_2 \cos x_1 \sin x_1$.

\subsection{Convergence: Preconditioner}
\label{s:precond-performance}

We study the performance of different variants of our preconditioner for the reduced space Hessian.

\ipoint{Setup} We solve the KKT system in~\eqref{e:hessian-matvec} at a true solution $\vect{v}^\star$. This velocity $\vect{v}^\star$ is found by registering two neuroimaging datasets from \nirep{} (\texttt{na01} and \texttt{na02}). The images are downsampled to a resolution of $128\times150\times128$ (half the original resolution). We consider an $H^1$-div regularization model with $\beta_v = \num{1E-2}$ and $\beta_w = \num{1E-4}$ and an $H^2$ regularization model with $\beta_v = \num{1E-4}$ with a tolerance $\epsilon_g = \num{1E-2}$ to compute $\vect{v}^\star$. Once we have found $\vect{v}^\star$, we generate a synthetic reference image $m_R$ by transporting the reference image using $\vect{v}^\star$. We use the velocity $\vect{v}^\star$ as an initial guess for our solver, and iteratively solve for the search direction $\vect{\tilde{v}}$ using different variants of our preconditioner. The number of time steps for the PDE solves is set to $n_t=4$. We fix the tolerance for the (outer) PCG method to $\epsilon_H = \num{1E-3}$. We consider an inexact Chebyshev semi-iterative method with a fixed number of $k\in\{5,10,20\}$ iterations (denoted by CHEB($k$)) and a nested PCG method with a tolerance of $\epsilon_P = 0.1\epsilon_H$ (denoted by PCG(\num{1E-1})) for the iterative inversion of the preconditioner. Details can be found in \secref{s:methods}. We compare these strategies to a spectral preconditioner (inverse of the regularization operator $\dop{A}$; used in~\cite{Mang:2016c,Gholami:2017a,Mang:2018a}). We study the rate of convergence of the PCG solver for a vanishing regularization parameter $\beta_v$. We consider mesh sizes of $128\times150\times128$ and $256\times300\times256$. We execute \claire{} on a single node of Opuntia with 20 MPI tasks.

\begin{figure}
\centering
\includegraphics[width=0.99\textwidth]
{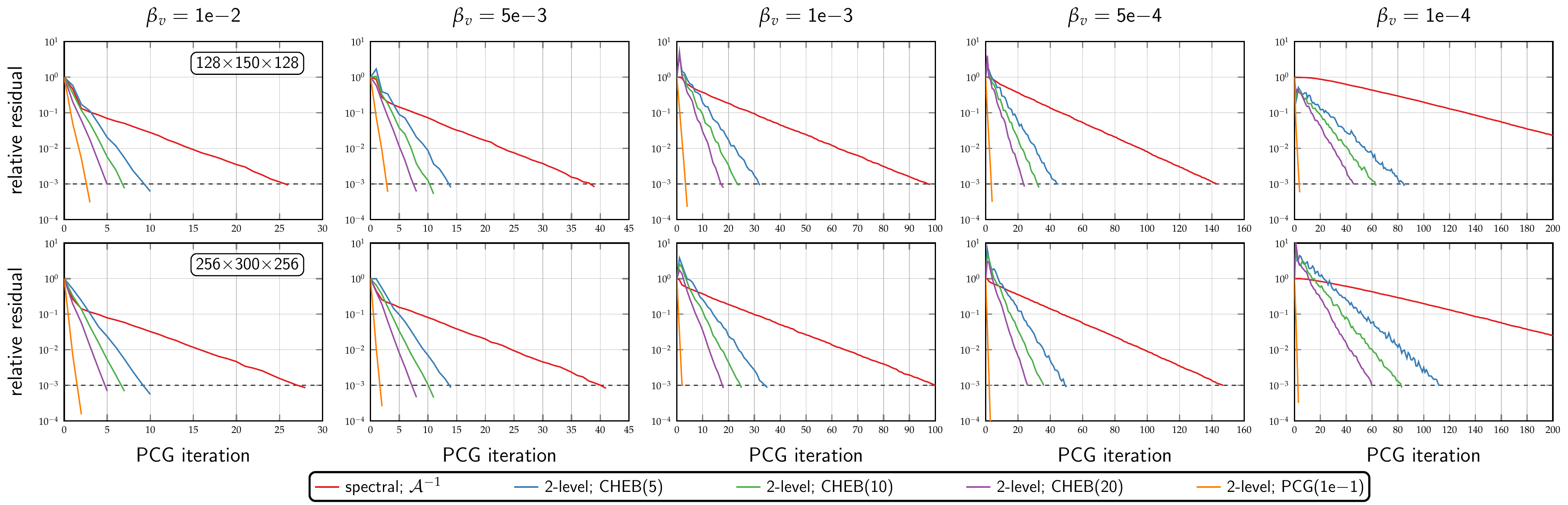}
\caption{Convergence of Krylov solver for different variants of the preconditioner for the reduced space Hessian. We consider an $H^2$-seminorm as regularization model for the velocity (neglecting the incompressibility constraint). We report results for different regularization parameters $\beta_v\in\{\num{1E-2},\num{5E-3},\num{1E-3},\num{5E-4},\num{1E-4}\}$. We report the trend of the relative residual for the outer Krylov method (PCG) versus the iteration count. We report results for the spectral preconditioner and the two-level preconditioner. We use different iterative algorithms to compute the action of the inverse of the preconditioner: CHEB($k$) with $k\in\{5,10,20\}$ refers to a CHEB method with a fixed number of $k$ iterations; PCG(\num{1E-1}) refers to a PCG method with a tolerance that is $0.1\times$ smaller than the tolerance used for the outer PCG method.\label{f:convergence-nirep-krylov-h2s}}
\end{figure}

\begin{figure}
\centering
\includegraphics[width=0.78\textwidth]
{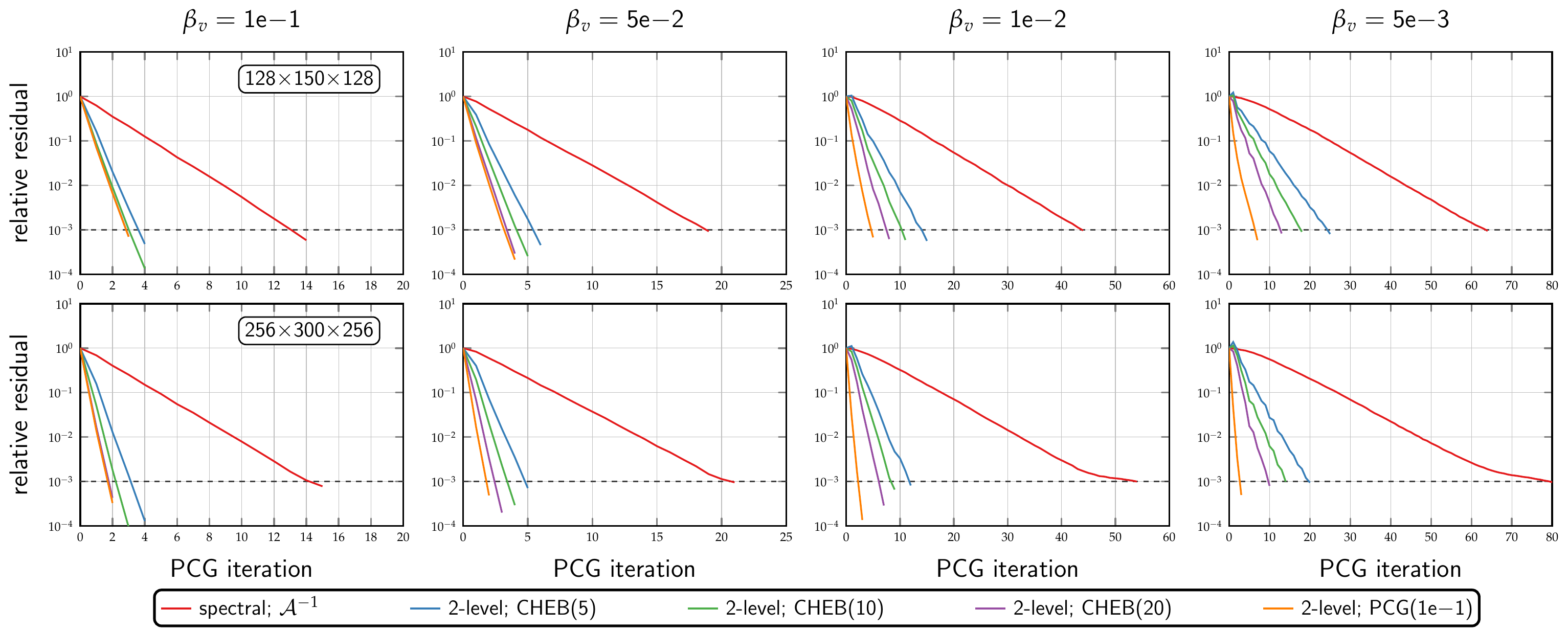}
\caption{Convergence of Krylov solver for different variants of the preconditioner for the reduced space Hessian. We consider an $H^1$-div regularization model with an $H^1$-seminorm for the velocity. We report results for different regularization parameters $\beta_v\in\{\num{1E-1},\num{5E-2},\num{1E-2},\num{5E-3}\}$. We set $\beta_w = \num{1e-4}$. We report the trend of the relative residual for the outer Krylov method (PCG) versus the iteration count. We report results for the spectral preconditioner and the two-level preconditioner. We use different algorithms to compute the action of the inverse of the preconditioner: CHEB($k$) with $k\in\{5,10,20\}$ refers to a CHEB method with a fixed number of $k$ iterations; PCG(\num{1E-1}) refers to a PCG method with a tolerance that is $0.1\times$ smaller than the tolerance used for the outer PCG method.\label{f:convergence-nirep-krylov-h1s-div}}
\end{figure}

\ipoint{Results} We display the trend of the residual with respect to the (outer) PCG iterations in \figref{f:convergence-nirep-krylov-h2s} ($H^2$-seminorm for $\vect{v}$, i.e., $\dop{B} = -\ilap$, with $\beta_v\in\{\num{1E-2},\num{5E-3},\num{1E-3},\num{5E-4},\num{1E-4}\}$) and in \figref{f:convergence-nirep-krylov-h1s-div} ($H^1$-div regularization model with an $H^1$-seminorm for $\vect{v}$, i.e., $\dop{B} = \igrad$ with penalty for $\idiv \vect{v}$, with $\beta_v\in\{\num{1E-1},\num{5E-2},\num{1E-2},\num{5E-3}\}$ and $\beta_w = \num{1e-4}$), respectively. Detailed results for these runs can be found in \tabref{t:convergence-nirep-krylov-h2s} and \tabref{t:convergence-nirep-krylov-h1s-div} in the supplementary materials.

\ipoint{Observations} The most important observations are:
\begin{itemize}[align=left,leftmargin=1.8em,itemindent=0pt,labelsep=0pt,labelwidth=1.2em]
\item The PCG method converges significantly faster for the two-level preconditioner.
\item The performance of all preconditioners considered in this study is \emph{not} independent of the regularization parameter $\beta_v$. The workload increases significantly for vanishing regularity of the velocity $\vect{v}$ for all preconditioners. The plots in \figref{f:convergence-nirep-krylov-h2s} and \figref{f:convergence-nirep-krylov-h1s-div} imply that the convergence of the PCG method for the two-level preconditioner is less sensitive to (or even independent of) $\beta_v$. The work goes to the inversion of the reduced space Hessian on the coarse grid (cf. \tabref{t:convergence-nirep-krylov-h2s} and \tabref{t:convergence-nirep-krylov-h1s-div} in the supplementary materials for details). If we further reduce the regularization parameter (below \num{1e-5} for the $H^2$-regularization model and below \num{1e-4} for the $H^1$-div regularization model) the performance of our preconditioners deteriorates further; the runtime becomes impractical for all preconditioners.
\item The rate of convergence of the PCG method is (almost) independent of the mesh size for all preconditioners. We note that we apply a smoothing of $\sigma = 2$ along each spatial dimension so that the input image data is resolved on the coarse grid of size $128\times150\times128$. The same frequency content is presented to the solver on the fine grid of size $256\times300\times256$.
\item The PCG method converges significantly faster if we consider an $H^1$-regularization model for $\vect{v}$. This is a direct consequence of fact that the condition number of the Hessian increases with the order of the regularization operator $\dop{A}$.
\item The differences of the performance of the preconditioners are less pronounced for an $H^1$-div regularization model for $\vect{v}$ than for an $H^2$-regularization model. For an $H^2$ regularization model with $\beta_v=\num{1e-4}$ we require more than 200 iterations for the spectral preconditioner.
\item Considering runtime (not reported here), we obtain a speedup of up to 2.9 for the $H^2$-regularization model (see \runref{20} in \tabref{t:convergence-nirep-krylov-h2s} in the supplementary materials) and a speedup of up to 2.6 for the $H^1$-div regularization model (see \runref{40} in \tabref{t:convergence-nirep-krylov-h1s-div} in the supplementary materials). The coarser the grid, the less effective is the two-level preconditioner, especially for vanishing regularization parameters $\beta_v$. This is expected, since we cannot resolve high-frequency components of the fine level on the coarse level. Secondly, we do not use a proper (algorithmic) smoother in our scheme to reduce the high-frequency errors.
\item The performance of the CHEB and the nested PCG method for iteratively inverting the reduced space Hessian are similar. There are differences in terms of the mesh size. For a coarser grid ($128\times150\times128$) the CHEB seems to perform slightly better. For a grid size of $256\times300\times256$ the nested PCG method is slightly better.
\end{itemize}

\ipoint{Conclusions} \bipa\item The two-level preconditioner is more effective than the spectral preconditioner. \item The nested PCG method is more effective than the CHEB method on a finer grid (and does not require a repeated estimation of the spectrum of the Hessian operator). \item The PCG method converges faster if we consider an $H^1$-div regularization model for $\vect{v}$. \item Designing a preconditioner that delivers a good performance for vanishing regularization parameters requires more work\eipa.

\subsection{Convergence: Newton--Krylov Solver}
\label{s:convergence-nks}

We study the rate of convergence of our Newton--Krylov solver for the entire inversion. We consider the neuroimaging data described in \secref{s:data}. We report additional results for a synthetic test problem (ideal case) in the supplementary materials.

\ipoint{Setup} We register the datasets \texttt{na02} through \texttt{na16} (template images) with \texttt{na01} (reference image). We execute the registration in full resolution ($256\times300\times256$; \num{58982400} unknowns). We consider an $H^1$-div regularization model ($H^1$-seminorm for $\vect{v}$ with $\beta_v=\num{1E-2}$ and $\beta_w = \num{1E-4}$; the parameters are chosen empirically). The number of Newton iterations is limited to 50 (not reached). The number of Krylov iterations is limited to 100 (not reached). We use a tolerance of $\num{5e-2}$ and $\num{1e-6}$ (the latter is not reached) for the relative reduction and the absolute $\ell^2$-norm of the reduced gradient as a stopping criterion, respectively. We use $n_t=4$ time steps for numerical time integration. We compare results obtained for the two-level preconditioner to results obtained using a spectral preconditioner (inverse of the regularization operator). We use a nested PCG method with a tolerance of $\epsilon_P = 0.1\epsilon_H$ for computing the action of the inverse of the two-level preconditioner. We do not perform any parameter, scale, or grid continuation. (We note that we observed that these continuation schemes are critical when performing runs for smaller regularization parameters.) We compare results obtained for single (\SI{32}{\bit}) and double (\SI{64}{\bit}) precision. We execute these runs on TACC's Lonestar 5 system (see \secref{s:compute-systems} for specs).

\begin{figure}
\centering
\includegraphics[width=0.99\textwidth]
{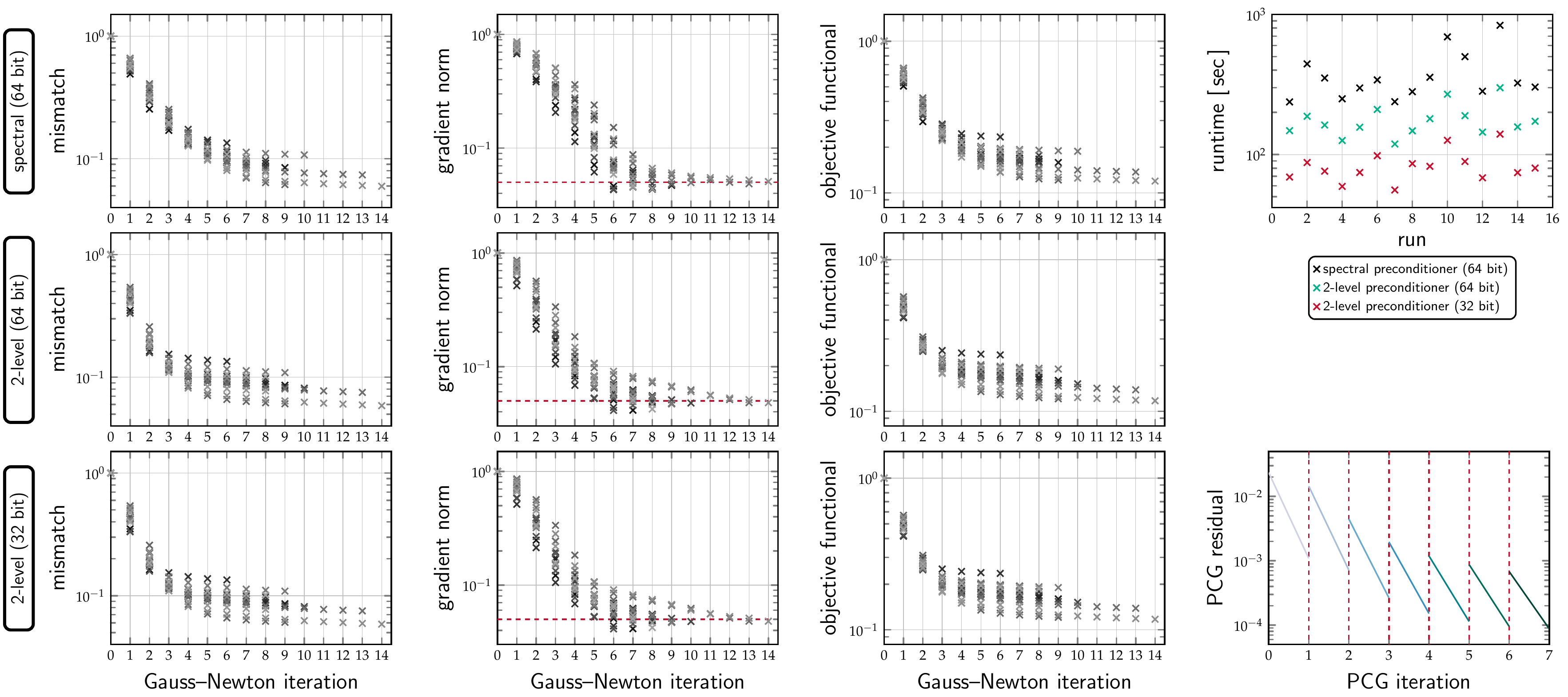}
\caption{Convergence of \claire{}\,'s Newton--Krylov solver for neuroimaging data for different realizations of the preconditioner. Top row: inverse regularization operator. Middle and bottom row: two-level preconditioner using PCG(\num{1E-1}) for double (\SI{64}{\bit}; middle row) and single (\SI{32}{\bit}; bottom row) precision, respectively. We report results for 15 multi-subject brain registration problems ({\tt na02} through {\tt na16} of the \nirep{} repository registered to {\tt na01}). Each of these 15 registrations is plotted in a different shade of gray. We plot (from left to right) the relative reduction of ($i$) the mismatch (squared $L^2$-distance between the images to be registered), ($ii$) the reduced gradient, and ($iii$) the objective functional, with respect to the Gauss--Newton iterations. We use a relative change of the gradient of \num{5e-2} as a stopping criterion (dashed red line in second column). We also report the runtime in seconds for each registration problem (right plot at top) and an exemplary plot of the reduction of the residual of the (outer) PCG solver per Newton iteration (right plot at bottom; the Newton iterations are separated by vertical dashed lines). The runs are performed on one node of TACC's Lonestar 5 system. The results reported here correspond to those in \tabref{t:nirep-convergence-h1sdiv-spectral-ls5}, \tabref{t:nirep-convergence-h1sdiv-2level-pcg0d1-dbl-ls5}, and \tabref{t:nirep-convergence-h1sdiv-2level-pcg0d1-sgl-ls5} in the supplementary materials.\label{f:convergence-nirep-fullsolve}}
\end{figure}

\ipoint{Results} We show convergence plots for all datasets in \figref{f:convergence-nirep-fullsolve}. We plot the relative reduction of the mismatch (left column), the relative reduction of the gradient (middle column), and the relative reduction of the objective functional (right column) with respect to the Gauss--Newton iterations. The top row shows results for the spectral preconditioner; the other two rows show results for the two-level preconditioner (middle row: double precision; bottom row: single precision). The runtime for the inversion is reported in the plot at the top right of \figref{f:convergence-nirep-fullsolve}. An exemplary trend for the residual of the PCG method per Gauss--Newton iteration is displayed at the bottom right of \figref{f:convergence-nirep-fullsolve}. These plots summarize results reported in the supplementary materials; results for the spectral preconditioner are reported in \tabref{t:nirep-convergence-h1sdiv-spectral-ls5}; results for the two-level preconditioner are reported in \tabref{t:nirep-convergence-h1sdiv-2level-pcg0d1-dbl-ls5} (double precision) and \tabref{t:nirep-convergence-h1sdiv-2level-pcg0d1-sgl-ls5} (double precision). We also report a comparison of the performance of our solver for single (\SI{32}{\bit}) and double (\SI{64}{\bit}) precision in \tabref{t:dbl-vs-sgl} for two exemplary images of the \nirep{} dataset.

\ipoint{Observations} The most important observations are the following:
\begin{itemize}[align=left,leftmargin=1.8em,itemindent=0pt,labelsep=0pt,labelwidth=1.2em]
\item Switching from double to single precision does not affect the convergence of our solver (see \figref{f:convergence-nirep-fullsolve}; detailed results are reported in \tabref{t:dbl-vs-sgl} in the supplementary materials).
\item The two-level preconditioner executed with single precision yields a speedup of up to 6$\times$ (with an average speedup of $4.4\pm 0.8$) compared to our baseline method (spectral preconditioner executed in double precision)~\cite{Mang:2016c,Gholami:2017a} (see \figref{f:convergence-nirep-fullsolve} top right). Switching from single to double precision yields a speedup of more than 2$\times$ (detailed results are reported in \tabref{t:dbl-vs-sgl} in the supplementary materials).
\item The average runtime of our improved solver is $\SI{85}{\second} \pm \SI{22}{\second}$ with a maximum of \SI{140}{\second} (see \runref{13} in \tabref{t:nirep-convergence-h1sdiv-2level-pcg0d1-sgl-ls5} in the supplementary materials for details) and a minimum of \SI{56}{\second} (see \runref{7} in \tabref{t:nirep-convergence-h1sdiv-2level-pcg0d1-sgl-ls5} in the supplementary materials for details).
\item We obtain a very similar convergence behavior for the outer Gauss--Newton iterations for different variants of our solver (see \figref{f:convergence-nirep-fullsolve}). We can reduce the $\ell^2$-norm of the gradient by $\num{5e-2}$ in 6 to 14 Gauss--Newton iterations (depending on the considered pair of images).
\item The mismatch between the deformed template image and the reference image stagnates once we have reduced the gradient by more than one order of magnitude (for the considered regularization parameter).
\item We oversolve the reduced space KKT system if we consider a superlinear forcing sequence in combination with a nested PCG method (see \figref{f:convergence-nirep-fullsolve} bottom right). This is different for synthetic data (we report exemplary results in the supplementary materials).
\end{itemize}

\ipoint{Conclusions} \bipa \item Our improved implementation of \claire{} yields an overall speedup of $4\times$ for real data if executed on a single resolution level. \item Executing \claire{} in single precision does not deteriorate the performance of our solver (if we consider an $H^1$-regularization model for the velocity)\eipa.

\begin{figure}
\begin{minipage}[]{0.60\textwidth}
\null\vspace{0.3cm}
\centering
\includegraphics[width=0.99\textwidth]
{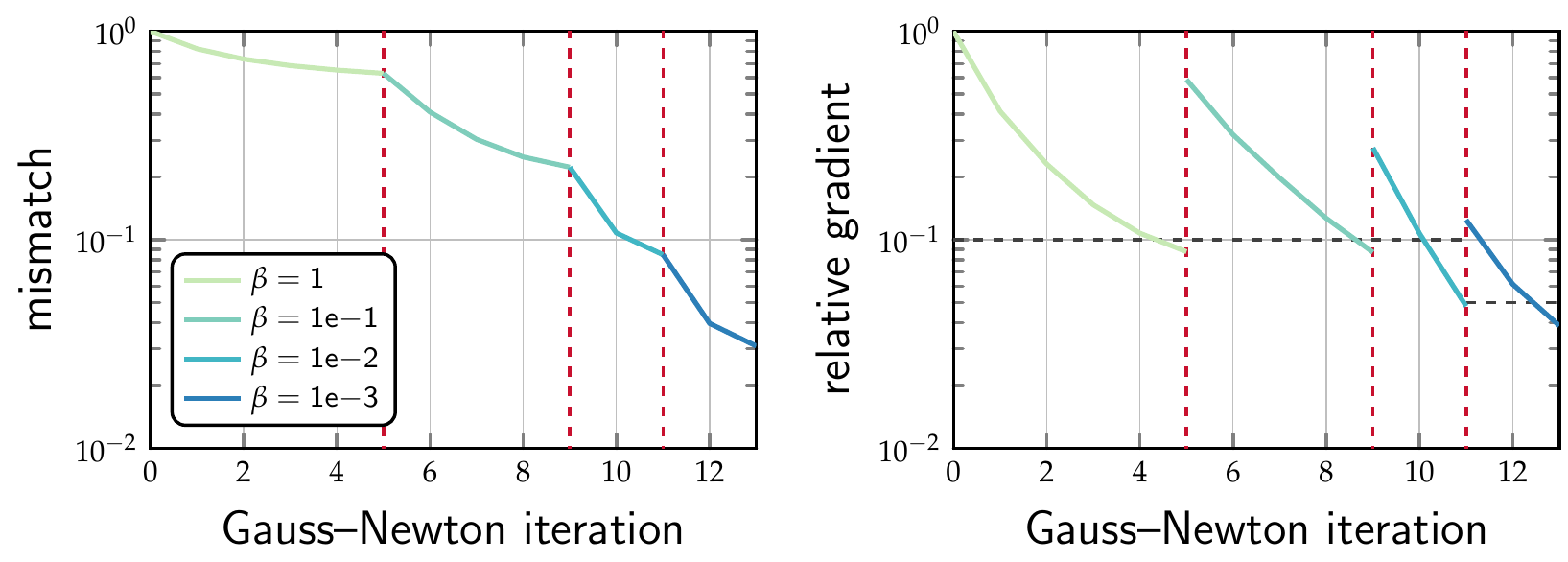}
\end{minipage}
\begin{minipage}[]{0.38\textwidth}
\caption{Convergence results for the parameter continuation scheme implemented in \claire{}. We report results for the registration of {\tt na11} to {\tt na01}. We report the reduction of the mismatch (left) the reduced gradient (right) per level (different regularization parameters) versus the cumulative number of Gauss--Newton iterations. (We require 5, 4, 2, and 2 Gauss--Newton iterations per level, respectively.) The individual levels are separated by vertical, dashed lines. The horizontal dashed lines in the right plot show the tolerance for the relative reduction of the gradient for the inversion.\label{f:convergence-nirep-parametercont}}
\end{minipage}
\end{figure}

\subsection{Time-to-Solution}
\label{s:time-to-solution}

We study the performance of \claire{}. We note that the \demons{} algorithm requires between approximately \SI{30}{\second} (3 levels with 15, 10, and 5 iterations) and \SI{3600}{\second} (3 levels with 1500, 1000 and 500 iterations) until `'convergence`' on the same system (depending on the parameter choices; see \secref{s:reg-quality} for details).

\begin{remark}
Since we perform a fixed number of iterations for the \demons{} algorithm, the runtime only depends on the execution time of the operators. The regularization parameters control the support of the Gaussian smoothing operator; the larger the parameters, the longer the execution time. This is different for \claire{}; large regularization parameters result in fast convergence and, hence, yield a short execution time. A simple strategy to obtain competitive results in terms of runtime would be to also execute \claire{} for a fixed number of iterations. We prefer to use a tolerance for the relative reduction of the gradient, instead, since it yields consistent results across different datasets.
\end{remark}

\ipoint{Setup} We use the dataset \texttt{na03}, \texttt{na10}  and \texttt{na11} as template images, and register them to \texttt{na01} (reference image). We consider and $H^1$-div regularization model ($H^1$-seminorm for $\vect{v}$ with $\beta_v\in\{\num{1e-2},\num{1e-3}\}$ and $\beta_w = \num{1E-4}$; these parameters are chosen empirically). The number of Newton iterations is limited to 50 (not reached). The number of Krylov iterations is limited to 100 (not reached). We use a tolerance of $\num{5e-2}$ for the relative reduction of the $\ell^2$-norm of the gradient and a tolerance of $\num{1e-6}$ (not reached) for its $\ell^2$-norm as a stopping criterion. We use $n_t=4$ time steps for numerical time integration. We compare results obtained for the two-level preconditioner (runs executed in single precision) to results obtained using a spectral preconditioner (inverse of the regularization operator; runs executed in double precision; the baseline method is described in~\cite{Mang:2016c}). We use a nested PCG method with a tolerance of $\epsilon_P = 0.1\epsilon_H$ for computing the action of the inverse of the two-level preconditioner. We execute \claire{} using a parameter continuation scheme. That is, we run the inversion until convergence for a sequence of decreasing regularization parameters (one order of magnitude, starting with $\beta_v=\num{1e0}$) until we reach the target regularization parameter. We execute these runs on one node of the Opuntia system using 20 MPI tasks (see \secref{s:compute-systems} for specs).

\begin{figure}
\centering
\includegraphics[width=0.99\textwidth]
{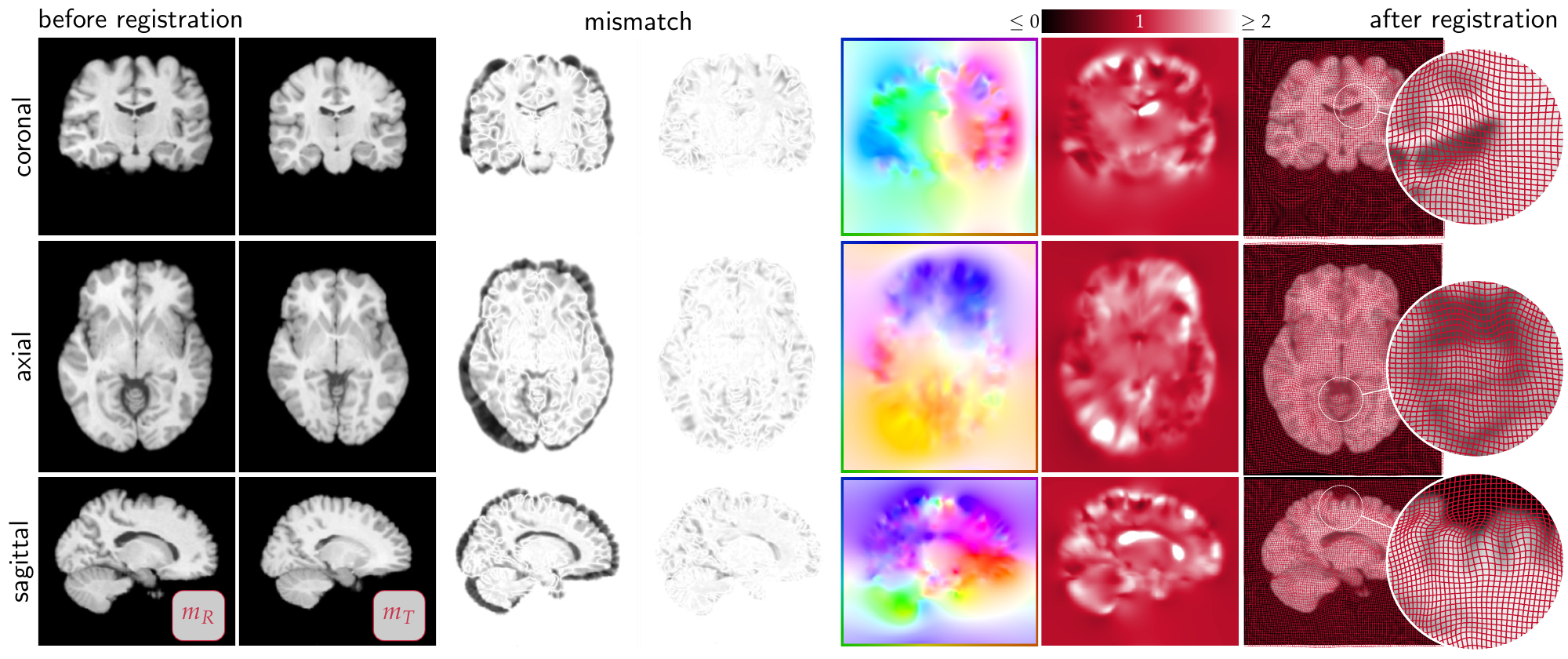}
\caption{Exemplary results for the parameter continuation scheme implemented in \claire{}. We consider the datasets {\tt na10} (template image) to {\tt na01} (reference image). We show (from top to bottom) coronal, axial and sagittal slices. The three columns on the left show the original data (left: reference image $m_R$; middle: template image $m_T$; right: mismatch between $m_R$ and $m_T$ before registration). The four columns on the right show results for the parameter continuation scheme (\runref{9} in \tabref{t:results-nirep-runtime}; from left to right: mismatch between $m_R$ and $m_1$ (after registration); a map of the orientation of $\vect{v}$; a map of the determinant of the deformation gradient (the color bar is shown at the top); and a deformed grid illustrating the in plane components of $\vect{y}$).\label{f:nirep-parameter-cont-na01-na10-h1sdiv}}
\end{figure}

\ipoint{Results} We report the results in \tabref{t:results-nirep-runtime}. We report the number of Gauss--Newton iterations, the number of Hessian matrix vector products (per level), the number of PDE solves (per level), the relative reduction of the mismatch, the $\ell^2$-norm of the reduced gradient, the relative reduction of the $\ell^2$-norm of the gradient, the runtime, and the associated speedup compared to a full solve disregarding any acceleration schemes. We showcase the trend of the mismatch and the $\ell^2$-norm of the gradient for different levels of the parameter continuation scheme in \figref{f:convergence-nirep-parametercont}. We report exemplary convergence results for the parameter continuation scheme in \figref{f:convergence-nirep-parametercont}. We show exemplary registration results for the parameter continuation in \figref{f:nirep-parameter-cont-na01-na10-h1sdiv} (for the registration of \texttt{na10} to \texttt{na01}).

\begin{table}
\caption{We compare different schemes implemented in \claire{} for stabilizing and accelerating the computations. We consider two datasets as a template image ({\tt na03} and {\tt na10}). We use an $H^1$-div regularization model with $\beta_w=\num{1E-4}$. We consider regularization parameters $\beta_v=\num{1e-2}$ and $\beta_v = \num{1e-3}$. We execute the inversion with a spectral preconditioner (double precision) to establish a baseline (\runref{1}, \runref{4}, \runref{7}, and \runref{10}; corresponds to the method presented in~\cite{Mang:2016c}). The remaining results are obtained with a two-level preconditioner using a nested PCG method with a tolerance of $0.1\epsilon_H$ to compute the action of the inverse of the preconditioner. For each dataset and each choice of $\beta_v$ we report results for a two-level preconditioner without any accelerations and a parameter continuation (PC) scheme. We report (from left to right) the number of Gauss--Newton iterations per level (\#iter; the total number for the entire inversion is the sum), the number of Hessian matvecs per level (\#matvecs; the total number for the entire inversion is the sum), the number of PDE solves (on the fine grid; \#PDE), the relative reduction of the mismatch, the absolute $\ell^2$-norm of the reduced gradient ($\|\di{g}^\star\|_2$), and the relative $\ell^2$-norm of the reduced gradient after convergence ($\|\di{g}^\star\|_{\text{rel}}$). We also report the runtime (in seconds) as well as the speedup compared to our baseline method presented in \cite{Mang:2016c}.\label{t:results-nirep-runtime}}
\resetrunid\centering\scriptsize
\begin{tabular}{llllrrrllllr}\toprule
       &               & $\beta_v$  &     & \#iter      & \#matvecs   & \#PDE    & mismatch                 & $\|\di{g}^\star\|_2$     & $\|\di{g}^\star\|_{\text{rel}}$ & runtime            & speedup  \\\midrule
\runid & \texttt{na03} & \num{1E-2} & --- & 9           & 83          & 187      & \num{8.465094487461E-02} & \num{4.630653548871E-04} & \num{4.710900411598E-02}        & \num{6.049768e+02} &          \\
\runid &               &            & --- & 9           & 9           & 39       & \num{8.597876018175E-02} & \num{4.653556640269E-04} & \num{4.734200401881E-02}        & \num{1.220632e+02} &  5.0     \\
\runid &               &            & PC  & 4,3,2       & 4,3,2       & 46       & \num{9.839029297325E-02} & \num{8.662219911489E-04} & \num{4.772381156960E-02}        & \num{9.325386e+01} &  6.5     \\
\midrule
\runid &               & \num{1E-3} & --- & 7           & 128         & 273      & \num{2.880830095328E-02} & \num{3.968630078429E-04} & \num{4.937184876106E-02}        & \num{8.974052e+02} &          \\
\runid &               &            & --- & 12          & 12          & 73       & \num{2.561926085770E-02} & \num{3.718886300290E-04} & \num{4.626490460159E-02}        & \num{7.169832e+02} &  1.3     \\
\runid &               &            & PC  & 4,3,2,2     & 4,3,2,2     & 56       & \num{3.368572976549E-02} & \num{8.252661203242E-04} & \num{4.546738044469E-02}        & \num{1.607081e+02} &  5.6     \\
\midrule\midrule
\runid & \texttt{na10} & \num{1E-2} & --- & 7           & 52          & 121      & \num{9.670308556120E-02} & \num{4.977777085935E-04} & \num{4.907589446181E-02}        & \num{3.843565e+02} &          \\
\runid &               &            & --- & 7           & 7           & 31       & \num{9.619403630495E-02} & \num{4.987742868252E-04} & \num{4.918244481087E-02}        & \num{9.347386e+01} &  4.1     \\
\runid &               &            & PC  & 3,3,2       & 3,3,2       & 42       & \num{1.095430627465E-01} & \num{9.546505170874E-04} & \num{4.975514858961E-02}        & \num{9.044438e+01} &  4.2     \\
\midrule
\runid &               & \num{1E-3} & --- & 7           & 134         & 285      & \num{3.168779332810E-02} & \num{3.458229436107E-04} & \num{4.240992672708E-02}        & \num{1.041283e+03} &          \\
\runid &               &            & --- & 8           & 16          & 51       & \num{3.110102564096E-02} & \num{3.853500820696E-04} & \num{4.727265238762E-02}        & \num{4.776165e+02} &  2.2     \\
\runid &               &            & PC  & 3,3,2,2     & 3,3,2,3     & 54       & \num{3.775701671839E-02} & \num{7.409917889163E-04} & \num{3.861953318119E-02}        & \num{1.873683e+02} &  5.6     \\
\bottomrule
\end{tabular}
\end{table}

\ipoint{Observations} The most important observations are the following:
\begin{itemize}[align=left,leftmargin=1.8em,itemindent=0pt,labelsep=0pt,labelwidth=1.2em]
\item The parameter continuation scheme in $\beta_v$ yields a speedup between $4\times$ and $6\times$ (\runref{3}, \runref{6}, \runref{9}, and \runref{12}  in~\tabref{t:results-nirep-runtime}) even if we reduce the target regularization parameter from $\num{1e-2}$ to $\num{1e-3}$. The runtime for this accelerated scheme ranges from \SI{9.044438e+01}{\second} (\runref{9}) and \SI{1.873683e+02}{\second} (\runref{12}) depending on problem and parameter selection.
\item The results obtained for the different schemes are qualitatively and quantitatively very similar. We obtain similar values for the relative mismatch, e.g., between \num{1.095430627465E-01} and \num{9.619403630495E-02} for $\beta_v=\num{1e-2}$ and between \num{3.775701671839E-02} and \num{3.110102564096E-02} for $\beta_v = \num{1e-3}$ for the registration of \texttt{na10} to \texttt{na01}.
\end{itemize}

\ipoint{Conclusions} \bipa \item Introducing the parameter continuation stabilizes the computations (similar results can be observed for grid and scale continuation schemes; not reported here). While the speedup for the preconditioner deteriorates as we reduce $\beta_v$ (see, e.g., \runref{2} and \runref{5} in~\tabref{t:results-nirep-runtime}), we can observe a speedup of about 5$\times$ for the parameter continuation scheme irrespective of $\beta_v$. We note that for small regularization parameters it is critical to execute \claire{} using a parameter continuation scheme. That is, for certain problems we observed a stagnation in the reduction of the gradient if \claire{} is executed without a parameter continuation scheme for small regularization parameters. We attribute this behavior to the accumulation of numerical errors in our scheme. This observation requires further exploration. \item Depending on the desired mismatch and regularity requirements, we achieve a runtime that is almost competitive with the \demons{} algorithm using the same system (i.e., the same number of cores). The peak performance in terms of speedup for \claire{} was achieved when using a grid continuation scheme (results not reported here), with a speedup of up to 17$\times$. However, as the regularity of the solution reduces, this speedup drops significantly; the parameter continuation is more stable. We expect to obtain a similar speedup with improved stability if we combine grid and parameter continuation. Designing an effective algorithm that combines these two approaches requires more work\eipa.

\subsection{Registration Quality}
\label{s:reg-quality}

We study registration accuracy for multi-subject image registration problems based on the \nirep{} dataset (see \secref{s:data}). We compare results for our method to different variants of the diffeomorphic \demons{} algorithm.

\ipoint{Setup} We consider the entire \nirep{} data repository. We register the dataset \texttt{na02} through \texttt{na16} (template images) to \texttt{na01} (reference image). The data has been rigidly preregistered~\cite{Christensen:2006a}. We do not perform an additional affine preregistration step. Each dataset comes with a label map that contains 32 labels (ground truth segmentations) identifying distinct gray matter regions (see \figref{f:nirep-data} for an example). We quantify registration accuracy based on the Dice coefficient (the optimal value is one) for these labels after registration. For ease of presentation we limit the evaluation to the union of the 32 labels (we report results for the individual 32 labels for \claire{} in \figref{f:claire-quality-individual-labels} of the supplementary materials). We assess the regularity of the computed deformation map based on the extremal values for the determinant of the deformation gradient. The analysis is limited to the foreground of the reference image (i.e., the area occupied by brain, identified by thresholding using a threshold of 0.05). We compare the performance of our method against different variants of the diffeomorphic \demons{} algorithm. We execute all runs on one node of the Opuntia system using 20 MPI tasks (see \secref{s:compute-systems} for specs).

\begin{itemize}[align=left,leftmargin=1.8em,itemindent=0pt,labelsep=0pt,labelwidth=1.2em]
\item \textit{\textbf{\demons{}}:} We consider (\emph{non}-)\emph{symmetric diffeomorphic} ((\textbf{S})\textbf{DDEM}; \emph{diffeomorphic update rule})~\cite{Vercauteren:2007a,Vercauteren:2009a}, and the (\emph{non}-)\emph{symmetric log-domain diffeomorphic} \demons{} algorithm ((\textbf{S})\textbf{LDDDEM}; (\emph{symmetric}) \emph{log-domain update rule})~\cite{Vercauteren:2008a}. We have tested different settings for these methods (see below). We limit our study to the default parameters suggested in the literature, online resources, and the manual of the software. We use the code available at \cite{demons-web}. We compile in release mode, with the {\tt -O3} option. The code has been linked against \itk\ version 4.9.1~\cite{Johnson:2015a,itk-web}. Notice that the implementation uses multithreading based on {\tt pthreads} to speed up the computations. We use the default setting, which corresponds to the number of threads being equal to the number of cores of the system. We use the \emph{symmetrized force} for the symmetric strategies. We consider the \emph{gradient of the deformed template} as a force for the non-symmetric strategies. We use a nearest-neighbor interpolation model to transform the label maps. We perform various runs to identify adequate parameters. For the first set of runs we use a three-level grid continuation scheme with 15, 10, and 5 iterations per level (the default), respectively. We estimate an optimal combination of regularization parameters $\sigma_u \geq 0$, $\sigma_d \geq 0$, and $\sigma_v \geq 0$ based on an exhaustive search. This search is limited to the datasets \texttt{na01} (reference image) and \texttt{na02} (template image). We define the optimal regularization parameter to be the one that yields the highest Dice score subject to the map $\vect{y}$ being diffeomorphic. We note that accurately computing $\det\igrad\vect{y}$ is challenging. The values reported in this study have to be considered with the numerical accuracy in mind. For \demons{} we report the values generated by the software. We refine this parameter search by increasing the number of iterations per level by a factor of 2, 5, 10, and 100 to make sure that we have `'converged`' to an `'optimal`' solution. We apply the best variants identified by this exhaustive search to the entire \nirep{} data.
\item \textit{\textbf{\claire{}}:} We consider an $H^1$-div regularization model ($H^1$-seminorm for $\vect{v}$, i.e., $\dop{A} = -\ilap$, with an additional penalty for $\idiv \vect{v}$). We set the regularization parameter for the penalty for the divergence of $\vect{v}$ to $\beta_w = \num{1E-4}$. To select an adequate regularization parameter $\beta_v$, we use a binary search. We set the bounds for the determinant of the deformation gradient to 0.25 and 0.30, respectively. We set the number of time steps of the SL scheme to $n_t=4$. The number of maximal iterations is set to 50 (not reached). The number of Krylov iterations is limited to 100 (not reached). We use a tolerance of $\num{5e-2}$ and $\num{1e-6}$ for the relative and absolute reduction of the reduced gradient as a stopping criterion. We use $n_t=4$ time steps for numerical integration. We run the registration on full resolution and (based on the experiments in \secref{s:time-to-solution}) use a parameter continuation scheme in $\beta_v$ to solve the registration problem. Probing for an optimal regularization parameter is expensive. We limit this estimation to the datasets {\tt na01} (reference image) and {\tt na02} (template image), assuming that we can estimate an adequate parameter for a particular application based on a subset of images. We execute \claire{} on the remaining images using the identified parameters. We compute $\det\igrad\vect{y}$ directly from $\vect{v}$ by solving a transport equation (see~\cite{Mang:2015a,Mang:2018a} for details). We transport the label maps to generate results that are consistent with the values reported for the determinant of the deformation map. This requires an additional smoothing (standard deviation: one voxel) and thresholding (threshold: 0.5) step.
\end{itemize}

\begin{figure}
\begin{minipage}[c]{0.6\textwidth}
\centering\vspace{0.3cm}
\includegraphics[width=0.98\textwidth]
{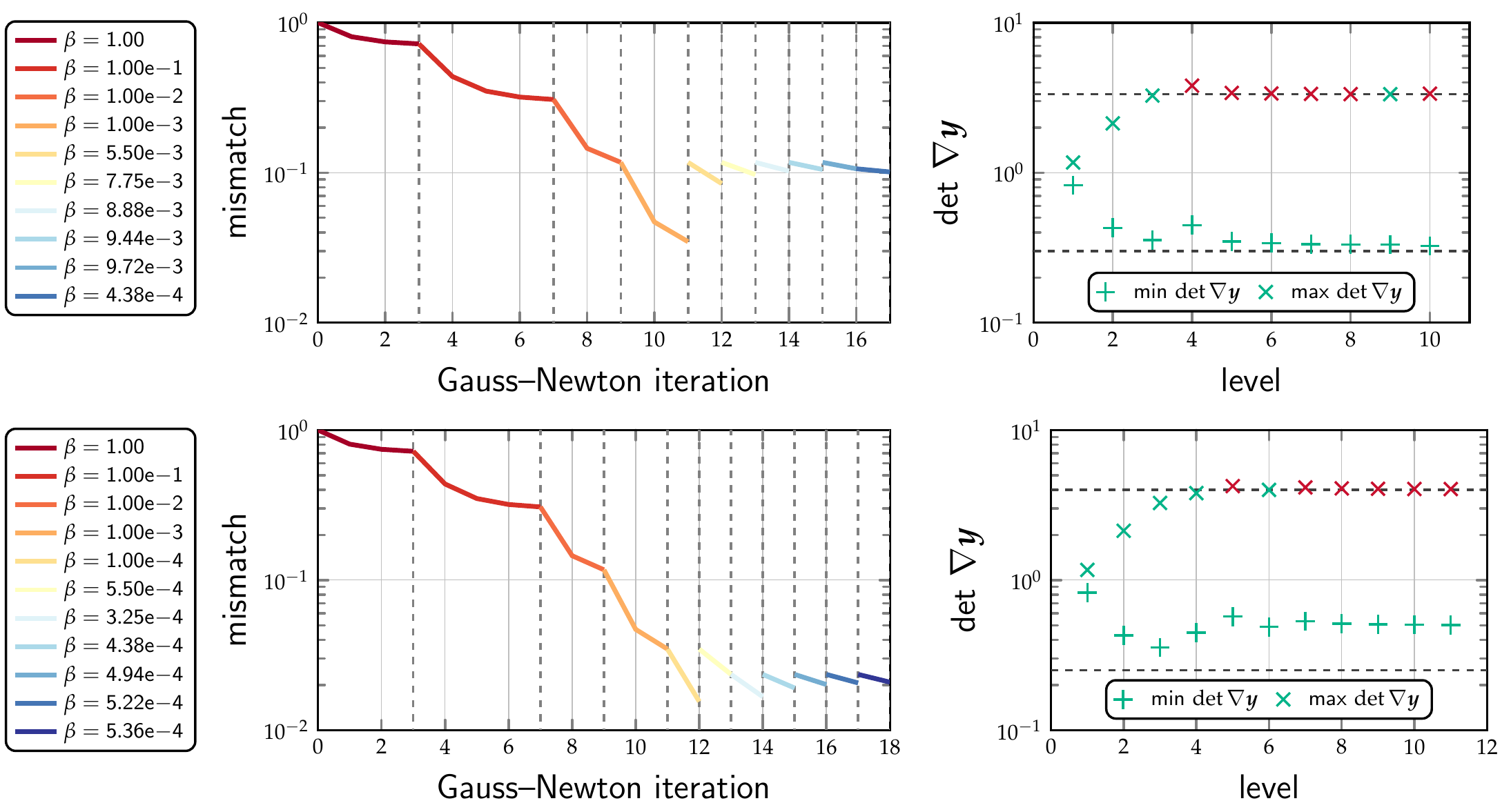}
\end{minipage}
\begin{minipage}[c]{0.37\textwidth}
\caption{Estimation of the regularization parameter $\beta_v$. We use an $H^1$-div regularization model with $\beta_w=\num{1E-4}$. We show the trend of the mismatch with respect to the Gauss--Newton iterations (left column) and the trend of the extremal values of the determinant of the deformation gradient with respect to the continuation level (right column). The top block shows results for a bound of 0.3 for $\min\det\igrad\vect{y}$. The bound for the bottom row is 0.25. These bounds are illustrated as dashed gray lines in the plots on the right. Here, we show (per continuation level) the trend of $\max\det\igrad\vect{y}$ (marker: $\times$) and $\min\det\igrad\vect{y}$ (marker: $+$). If the bounds are violated, we display the marker in red. We separate the continuation levels with a vertical gray line in the plots for the mismatch; the color of the line corresponds to a particular regularization parameter (see legend).\label{f:claire-estimate-regularization-parameter-na01-na02-h1sdiv}}
\end{minipage}
\end{figure}

\begin{figure}
\begin{minipage}[c]{0.37\textwidth}
\caption{Convergence results for \claire{} and SDDEM. We report the trend of the mismatch (left) and the Dice coefficient (right) versus the outer iterations. For \claire{}, we solve this problem more accurately than in the other runs on the real data to show the asymptotic behaviour of our solver. We do not perform any grid, scale, or parameter continuation for both methods. We consider the datasets {\tt na01} (reference image) and {\tt na02} (template image).\label{f:convergence-nirep-claire-vs-demons}}
\end{minipage}
\begin{minipage}[c]{0.6\textwidth}
\centering
\includegraphics[width=0.98\textwidth]
{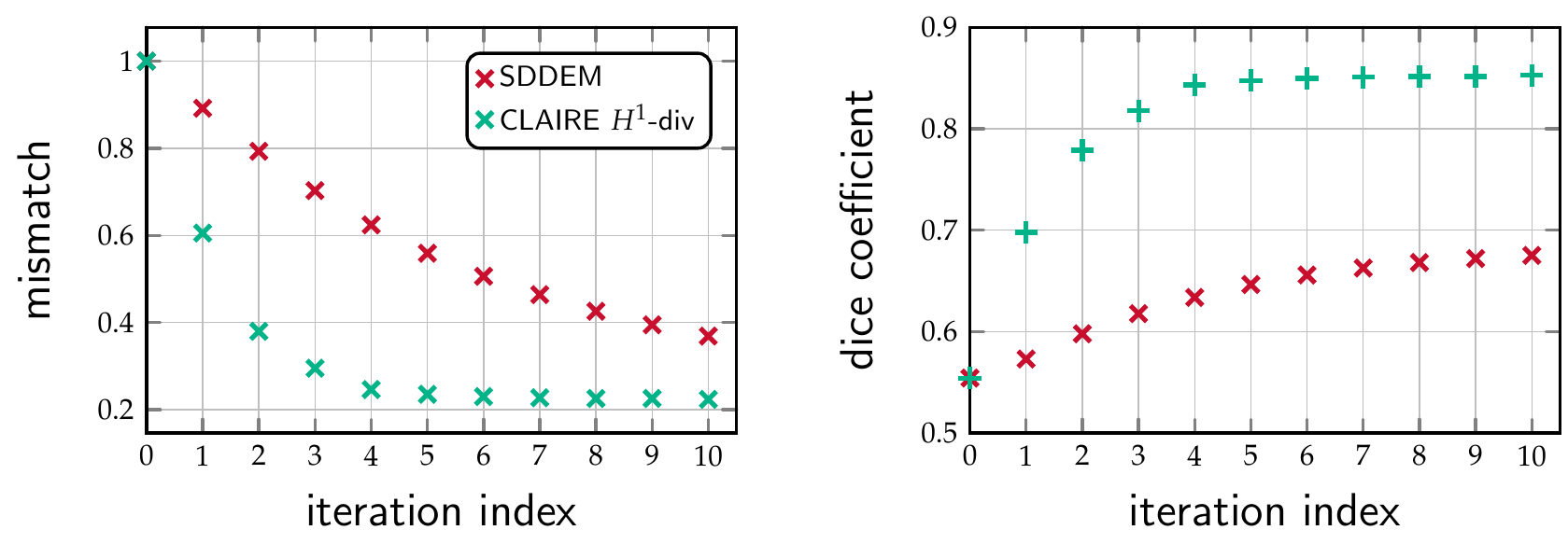}
\end{minipage}
\end{figure}

\begin{figure}
\begin{minipage}[c]{0.55\textwidth}
\centering\vspace{0.3cm}
\includegraphics[width=0.98\textwidth]
{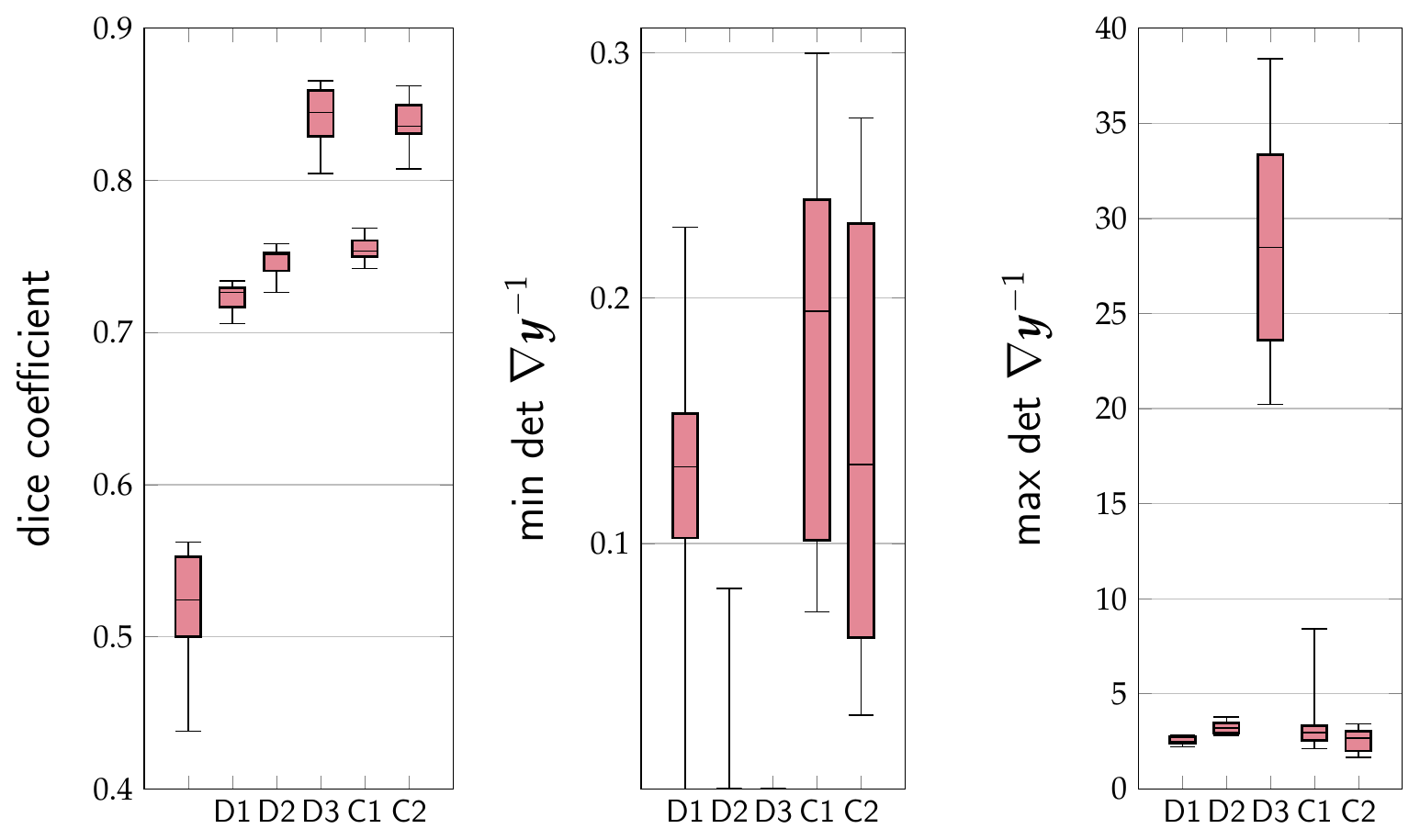}
\end{minipage}
\begin{minipage}[c]{0.43\textwidth}
\caption{Registration results for the \nirep{} data. We consider three variants of the diffeomorphic \demons{} algorithm: D1 corresponds to SDDEM with $(\sigma_u,\sigma_d) = (0,3.5)$, D2 to SDDEM with $(\sigma_u,\sigma_d) = (0,3.0)$, and D3 to SDDEM with $(\sigma_u,\sigma_d) = (0,1.0)$. These choices are based on an exhaustive search (we refer the interested reader to the supplementary materials for details). For \claire{} we use two different choices of the regularization parameter for the $H^1$-div regularization model (C1 corresponds to \claire{} with $\beta_v = \num{9.718750e-03}$ and C2 to \claire{} with $\beta_v = \num{5.500000e-04}$; these parameters are determined via a binary search (see \figref{f:claire-estimate-regularization-parameter-na01-na02-h1sdiv})). We report results for the entire \nirep{} dataset. The plot on the left shows the Dice coefficient (on the very left, we also provide a box plot for the Dice coefficient before registration). This coefficient is computed for the union of all gray matter labels (to simplify the analysis). The middle and right box plot show the extremal values for the determinant of the deformation gradient.\label{f:claire-vs-demons-quality-nirep}}
\end{minipage}
\end{figure}

\ipoint{Results} We illustrate the search for an optimal regularization parameter for \claire{} in \figref{f:claire-estimate-regularization-parameter-na01-na02-h1sdiv}. We showcase an exemplary result for the rate of convergence of SDDEM and \claire{} in \figref{f:convergence-nirep-claire-vs-demons} (the software is executed at full image resolution). We summarize exemplary registration results for all datasets in \figref{f:claire-vs-demons-quality-nirep}. Here, D1, D2, D3, C1, and C2 correspond to different variants of the \demons{} algorithm and \claire{}. C1 corresponds to \claire{} with regularization parameter of \num{9.718750e-03} ($\epsilon_J = 0.3$) and C2 to \claire{} with a regularization parameter of \num{5.500000e-04} ($\epsilon_J = 0.25$). The first \demons{} variant D1 is SDDEM with $(\sigma_u,\sigma_d) = (0,3.5)$ (smooth setting). It yields results that are competitive with \claire{} in terms of the determinant of the deformation gradient. The second variant D2 is SDDEM with $(\sigma_u,\sigma_d) = (0,3.0)$, which gave us the best result (highest attainable Dice score with the determinant of the deformation gradient not changing sign for the training data \texttt{na01} and \texttt{na02}). The third variant D3 is SDDEM with $(\sigma_u,\sigma_d) = (0,1.0)$ (aggressive setting). We achieve results that are competitive with \claire{} in terms of the Dice score. We execute the \demons{} algorithm with a three-level grid continuation scheme with 150, 100, and 50 iterations per level, respectively.

We refer the interested reader to the supplementary materials for more detailed results for these runs and an additional insight into the parameter search we have conducted to identify the best variant of the \demons{} algorithm. Detailed results for the \claire{} variant C1 are reported in \tabref{t:nirep-regquality-h1sdiv-3d3e-3}. Detailed results for the \claire{} variant C2 are reported in \tabref{t:nirep-regquality-h1sdiv-5d5e-4}. For \claire{}, we report Dice coefficients for the individual 32 gray matter labels in \figref{f:claire-quality-individual-labels}. Results for probing for adequate regularization parameters $\sigma_u$, $\sigma_d$, and $\sigma_v$ for different variants of the \demons{} algorithm are reported in \tabref{t:ddem:sweep} and \tabref{t:lddem:sweep} (exhaustive search). Building up on these results we extend this search by additionally increasing the iteration count. These results are reported in \tabref{t:dem-iterations}. We determined that SDDEM gives us the best results in terms of the Dice coefficient. Detailed results for the variants D1, D2, and D3 can be found in \tabref{t:demons-all-nirep-data}.

\ipoint{Observations} The most important observations are the following:
\begin{itemize}[align=left,leftmargin=1.8em,itemindent=0pt,labelsep=0pt,labelwidth=1.2em]
\item \claire{} yields a smaller mismatch/higher Dice coefficient with a better control of the determinant of the deformation gradient (see \figref{f:claire-vs-demons-quality-nirep}). We obtain an average Dice coefficient of \num{8.376830e-01} with $(\min,\max) = (\num{4.137990e-01}, \num{1.106347e+01})$ as extremal values for the determinant of the deformation gradient (on average). The Dice score for the best variant of the \demons{} algorithm, SDDEM, is $\num{8.422854506E-01}$. To attain this score we have to commit to nondiffeomorphic deformation maps (as judged by the values for the determinant of the deformation gradient reported by the \demons{} software). An extension of \claire{}, which we did not consider in this work, is to enable a monitor for the determinant of the deformation gradient that increases the regularization parameter if we hit the bound we used to estimate $\beta_v$. This would prevent the outliers we observe in this study, without having to probe for a new regularization parameter for each individual dataset.
\item For \claire{}, the average runtime (across all registrations) is \SI{1.084713e+02}{\second} and \SI{2.425679e+02}{\second} for $\beta_v = \num{9.718750e-03}$ and $\beta_v = \num{5.500000e-04}$, respectively. This is between $1.5\times$ and $5\times$ slower than the \demons{} algorithm if we execute \demons{} using 15, 10, and 5 iterations per level. Notice that \demons{} is executed for a fixed number of iterations. The runs reported here use 10$\times$ more iterations per level (which slightly improves the performance of \demons{}; we refer the interested reader to \tabref{t:dem-iterations} in the supplementary materials for details). This increases the runtime of the \demons{} algorithm by roughly a factor of 10. \claire{} uses a relative tolerance for the gradient as a stopping criterion. Moreover, \demons{} uses a grid continuation scheme. We execute these runs on the fine resolution, and perform a parameter continuation instead (since we observed it is more stable for vanishing $\beta_v$; see \secref{s:time-to-solution}).
\item On the fine grid (single-level registration), \claire{} converges significantly faster than the \demons{} algorithm. We reach a Dice score of more than 0.8 for \claire{} after only three Gauss--Newton iterations (see \figref{f:convergence-nirep-claire-vs-demons}).
\end{itemize}

\ipoint{Conclusions} With \claire{} we achieve \bipa\item a computational performance that is close to that of the \demons{} algorithm ($1.5\times$ to $5\times$ slower for the fastest setting we used for \demons{}) with \item a registration quality that is superior (higher Dice coefficient with a better behaved determinant of the deformation gradient)\eipa.

\subsection{Scalability}
\label{s:scalability}

We study strong scaling of our new implementation of \claire{} for up to \num{3221225472} unknowns for a synthetic test problem consisting of smooth trigonometric functions (see \secref{s:data}).

\ipoint{Setup} We consider grid sizes $128^3$, $256^3$, $512^3$, and $1024^3$. We use an $H^1$-div regularization model with $\beta_w=\num{1E-3}$ and $\beta_w=\num{1E-4}$. We use the two-level preconditioner with a nested PCG method with a tolerance of $0.1\epsilon_H$ to compute the action of the inverse of the preconditioner. We set the tolerance for the stopping condition for the relative reduction of the reduced gradient to $\num{1E-2}$ (with an absolute tolerance of $\num{1E-6}$ (not reached)). We execute the runs on TACC's Lonestar 5 system (see \secref{s:compute-systems} for specs).

\begin{figure}
\centering
\includegraphics[width=0.8\textwidth]
{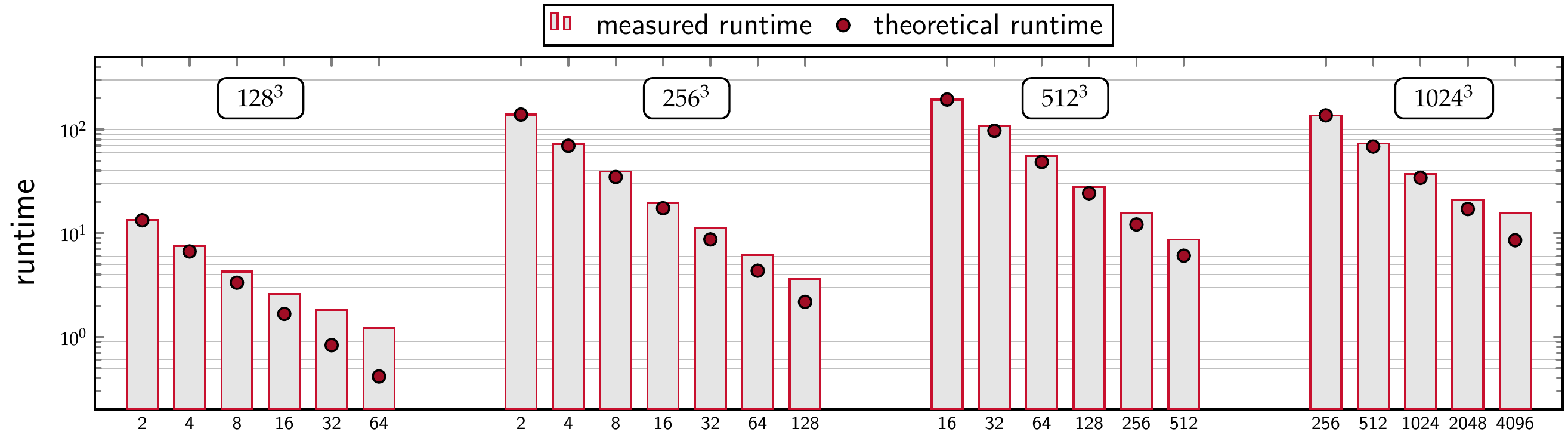}
\caption{Strong scaling results for a synthetic test problem on TACC's Lonestar 5 system (see \secref{s:compute-systems} for specs). We use 12 MPI tasks per node. We report the runtime (time-to-solution) for the entire inversion (in seconds). Our Newton--Krylov solver converges in three iterations (with three Hessian matvecs and a total of 15 PDE solves on the fine level). We consider grid sizes $128^3$, $256^3$, $512^3$, and $1024^3$ (from left to right). The largest run uses 4096 MPI tasks on 342 compute nodes (we solve for \num{3221225472} unknowns).\label{f:claire-strong-scaling}}
\end{figure}

\ipoint{Results} We report strong scaling results for \claire{} in \figref{f:claire-strong-scaling}. We report the time-to-solution and compare it to the runtime we expect theoretically. We report detailed results, which form the basis of the runtime reported in \figref{f:claire-strong-scaling}, in~\tabref{t:scalability-synthetic-ls5}. Here, we report the execution time of the FFT and the interpolation kernels on the coarse (two-level preconditioner) and fine grid, the runtime of our solver (time-to-solution), and the strong scaling efficiency of our improved implementation of \claire{}. We refer the reader to \cite{Mang:2016c,Gholami:2017a} more detailed results on the scalability of our original implementation of \claire{}

\begin{table}
\caption{Scalability results for \claire{} for a synthetic test problem. We report strong scaling results for up to \num{3221225472} unknowns (grid sizes: $128^3$, $256^3$, $152^3$, and $1024^3$). We execute these runs on TACC's Lonestar 5 system (see \secref{s:compute-systems} of the main manuscript for the specs). We consider an $H^1$-div regularization model with $\beta_v = \num{1e-3}$ and $\beta_w = \num{1e-4}$. We use a two-level preconditioner with a nested PCG method. We terminate the inversion if the gradient is reduced by $\num{1e-2}$. We execute these runs in single precision. We use 12 MPI tasks per node. We report the execution time for the FFT and the interpolation (on the coarse and the fine grid; in seconds), the runtime of the solver (time-to-solution; in seconds), and the strong scaling efficiency.\label{t:scalability-synthetic-ls5}}
\resetrunid\centering\scriptsize
\begin{tabular}{rrrRlrlrlrlrlL}\toprule
grid      & run    & nodes & tasks  & \multicolumn{4}{l}{fine grid}                                   & \multicolumn{4}{l}{coarse grid}                                     & runtime            & efficiency \\
\midrule
\multicolumn{4}{c}{\null}  & \multicolumn{2}{l}{FFT}     & \multicolumn{2}{l}{interpolation} & \multicolumn{2}{l}{FFT}     & \multicolumn{2}{l}{interpolation}              &                    &            \\
\midrule
$128^3$   & \runid &     1 &      2 & \num{4.247586e+00} & (32\%) & \num{2.817606e+00} & (21\%)       & \num{2.213636e+00} & (17\%) & \num{1.726949e+00} & (13\%)           & \num{1.332395e+01} &            \\
          & \runid &     1 &      4 & \num{2.445537e+00} & (32\%) & \num{1.457371e+00} & (19\%)       & \num{1.294543e+00} & (17\%) & \num{9.215655e-01} & (12\%)           & \num{7.488414e+00} & 89\%       \\
          & \runid &     1 &      8 & \num{1.346370e+00} & (32\%) & \num{8.137712e-01} & (19\%)       & \num{7.315106e-01} & (17\%) & \num{5.128467e-01} & (12\%)           & \num{4.259472e+00} & 78\%       \\
          & \runid &     2 &     16 & \num{7.388508e-01} & (28\%) & \num{5.686014e-01} & (22\%)       & \num{4.366422e-01} & (17\%) & \num{3.114562e-01} & (12\%)           & \num{2.594392e+00} & 64\%       \\
          & \runid &     3 &     32 & \num{4.160264e-01} & (23\%) & \num{3.906395e-01} & (21\%)       & \num{3.778143e-01} & (21\%) & \num{2.545202e-01} & (14\%)           & \num{1.819628e+00} & 46\%       \\
          & \runid &     6 &     64 & \num{3.119893e-01} & (26\%) & \num{3.451383e-01} & (28\%)       & \num{1.522639e-01} & (13\%) & \num{1.215253e-01} & (10\%)           & \num{1.211427e+00} & 34\%       \\
\midrule
$256^3$   & \runid &     1 &      2 & \num{5.550650e+01} & (40\%) & \num{2.770376e+01} & (20\%)       & \num{2.079800e+01} & (15\%) & \num{1.471913e+01} & (11\%)           & \num{1.392794e+02} &            \\
          & \runid &     1 &      4 & \num{2.698466e+01} & (37\%) & \num{1.405047e+01} & (19\%)       & \num{1.177014e+01} & (16\%) & \num{7.591873e+00} & (10\%)           & \num{7.231571e+01} & 96\%       \\
          & \runid &     1 &      8 & \num{1.449565e+01} & (37\%) & \num{7.703082e+00} & (20\%)       & \num{6.300269e+00} & (16\%) & \num{4.136640e+00} & (11\%)           & \num{3.919360e+01} & 89\%       \\
          & \runid &     2 &     16 & \num{6.870668e+00} & (35\%) & \num{3.503119e+00} & (18\%)       & \num{3.412268e+00} & (18\%) & \num{2.126711e+00} & (11\%)           & \num{1.946021e+01} & 89\%       \\
          & \runid &     3 &     32 & \num{4.058555e+00} & (36\%) & \num{1.938044e+00} & (17\%)       & \num{2.014238e+00} & (18\%) & \num{1.153192e+00} & (10\%)           & \num{1.131656e+01} & 77\%       \\
          & \runid &     6 &     64 & \num{2.151984e+00} & (35\%) & \num{1.041131e+00} & (17\%)       & \num{1.052281e+00} & (17\%) & \num{6.378655e-01} & (10\%)           & \num{6.141388e+00} & 71\%       \\
          & \runid &    11 &    128 & \num{1.204325e+00} & (33\%) & \num{6.262128e-01} & (17\%)       & \num{5.915005e-01} & (16\%) & \num{3.903782e-01} & (11\%)           & \num{3.628621e+00} & 60\%       \\
          & \runid &    22 &    256 & \num{7.077870e-01} & (30\%) & \num{4.376023e-01} & (18\%)       & \num{3.341010e-01} & (14\%) & \num{2.577584e-01} & (11\%)           & \num{2.337411e+00} & 47\%       \\
\midrule
$512^3$   & \runid &     2 &     16 & \num{8.013097e+01} & (41\%) & \num{3.261554e+01} & (17\%)       & \num{3.391836e+01} & (17\%) & \num{1.852011e+01} & (10\%)           & \num{1.942798e+02} &            \\
          & \runid &     3 &     32 & \num{4.523542e+01} & (41\%) & \num{1.790180e+01} & (16\%)       & \num{1.944877e+01} & (18\%) & \num{9.856170e+00} & (\phantom{0}9\%) & \num{1.089049e+02} & 89\%       \\
          & \runid &     6 &     64 & \num{2.212043e+01} & (40\%) & \num{8.871014e+00} & (16\%)       & \num{1.027469e+01} & (19\%) & \num{5.084882e+00} & (\phantom{0}9\%) & \num{5.539895e+01} & 88\%       \\
          & \runid &    11 &    128 & \num{1.069926e+01} & (38\%) & \num{4.299843e+00} & (15\%)       & \num{5.591518e+00} & (20\%) & \num{2.681119e+00} & (10\%)           & \num{2.809870e+01} & 86\%       \\
          & \runid &    22 &    256 & \num{5.698979e+00} & (37\%) & \num{2.262838e+00} & (15\%)       & \num{3.158469e+00} & (20\%) & \num{1.576031e+00} & (10\%)           & \num{1.559486e+01} & 78\%       \\
          & \runid &    43 &    512 & \num{2.997216e+00} & (35\%) & \num{1.450647e+00} & (17\%)       & \num{1.403620e+00} & (16\%) & \num{9.389393e-01} & (11\%)           & \num{8.659814e+00} & 70\%       \\
\midrule
$1024^3$  & \runid &    22 &    256 & \num{5.688231e+01} & (42\%) & \num{2.156576e+01} & (16\%)       & \num{2.679725e+01} & (20\%) & \num{1.135775e+01} & (\phantom{0}8\%) & \num{1.367671e+02} &            \\
          & \runid &    43 &    512 & \num{2.853562e+01} & (39\%) & \num{1.060915e+01} & (14\%)       & \num{1.703841e+01} & (23\%) & \num{6.417702e+00} & (\phantom{0}9\%) & \num{7.339710e+01} & 93\%       \\
          & \runid &    86 &   1024 & \num{1.451158e+01} & (39\%) & \num{5.225239e+00} & (14\%)       & \num{7.749017e+00} & (21\%) & \num{3.290969e+00} & (\phantom{0}9\%) & \num{3.735317e+01} & 92\%       \\
          & \runid &   171 &   2048 & \num{7.222220e+00} & (35\%) & \num{3.263639e+00} & (16\%)       & \num{4.129433e+00} & (20\%) & \num{2.150985e+00} & (10\%)           & \num{2.077145e+01} & 82\%       \\
          & \runid &   342 &   4096 & \num{4.488453e+00} & (28\%) & \num{2.300169e+00} & (15\%)       & \num{3.201794e+00} & (21\%) & \num{1.756468e+00} & (11\%)           & \num{1.553637e+01} & 55\%       \\
\bottomrule
\end{tabular}
\end{table}

\ipoint{Observations} The most important observations are the following:
\begin{itemize}[align=left,leftmargin=1.8em,itemindent=0pt,labelsep=0pt,labelwidth=1.2em]
\item We obtain a good strong scaling efficiency that is at the order of 60\%.
\item The strong scaling results are in accordance with the performance reported in~\cite{Mang:2016a,Gholami:2017a}. The key difference is that the scalability of our new solver is dominated by the coarse grid discretization within the preconditioner. That is, we do not observe the scalability reported in~\cite{Mang:2016a,Gholami:2017a} if we execute \claire{} with the same amount of resources for a given resolution of the data. However, if we compare the scalability results reported in~\cite{Gholami:2017a} with a resolution that matches the coarse grid in the preconditioner, we can observe a similar strong scaling efficiency.
\item We can solve clinically relevant problems in about \SI{2}{\second} if we execute \claire{} with 256 MPI tasks (see \runref{14} in \tabref{t:scalability-synthetic-ls5}).
\item We can solve problems with up to \num{3221225472} unknowns in less then \SI{5}{\second} with 4096 MPI tasks on 342 compute nodes on TACC's Lonestar 5 system (see \runref{25} in \tabref{t:scalability-synthetic-ls5}). The solver converges in \SI{1.367671e+02}{\second} if we execute the run on 22 nodes with 256 MPI tasks.
\end{itemize}

\ipoint{Conclusions} With \claire{} we deploy a solver that scales on HPC platforms. \claire{} approaches run-times that represent a significant step towards providing \iquote{real-time} capabilities for clinically relevant problem sizes (inversion for $\sim$50 million unknowns in \SI{2.34}{\second} using 256 MPI tasks; see also~\cite{Mang:2016c,Gholami:2017a}). \claire{} provides fast solutions on moderately sized clusters (which could potentially be deployed to hospitals). We note that \claire{} does not require a cluster. It can be executed on individual compute nodes. Further accelerations on reduced hardware resources form the basis of our current work. \claire{} can also be used to solve diffeomorphic image registration problems of unprecedented scale, something that is of interest for whole body imaging~\cite{Tarnoki:2015a,Lecouvet:2016a} or experimental, high-resolution microscopic imaging~\cite{Kutten:2017a,Chung:2013a,Tomer:2014a}. The largest problem we have solved with our original implementation of \claire{} is \num{25769803776} unknowns (see~\cite{Gholami:2017a}). To the best of our knowledge, \claire{} is the only software for large deformation diffeomorphic registration with these capabilities.

\section{Conclusions}
\label{s:conclusions}

With this publication we release \claire{}, a memory-distributed algorithm for stationary velocity field large deformation diffeomorphic image registration in 3D. This work builds up on our former contributions on constrained large deformation diffeomorphic image registration~\cite{Mang:2015a,Mang:2016a,Mang:2016c,Mang:2017b,Gholami:2017a,Mang:2018a}. We have performed a detailed benchmark study of the performance of \claire{} on synthetic and real data. We have studied the convergence for different schemes for preconditioning the reduced space Hessian in~\secref{s:precond-performance}. We have examined the rate of convergence of our Gauss--Newton--Krylov solver in \secref{s:convergence-nks}. We have reported results for different schemes available in \claire{} in \secref{s:time-to-solution} to study the time-to-solution. We have compared the registration quality obtained with \claire{} to different variants of the diffeomorphic \demons{} algorithm in \secref{s:reg-quality}. We have also reported strong scaling results for our improved memory-distributed solver on supercomputing platforms (see \secref{s:scalability}). We note that we accompany this work with supplementary materials that provides a more detailed picture about the performance of our method. The most important conclusions are the following:
\begin{itemize}[align=left,leftmargin=1.8em,itemindent=0pt,labelsep=0pt,labelwidth=1.2em]
\item \claire{} delivers high-fidelity results with well-behaved deformations. Our results are in accordance with observations we have made for the two-dimensional case~\cite{Mang:2016a}. Our $H^1$-div formulation outperforms the diffeomorphic \demons{} algorithm in terms of data fidelity and deformation regularity (as judged by the higher dice score and more well-behaved extremal values for the determinant of the deformation gradient; see \figref{f:claire-vs-demons-quality-nirep} in \secref{s:reg-quality}).
\item Our Gauss--Newton--Krylov solver converges after only a few iterations to high-fidelity results. The rate of convergence of \claire{} is significantly better than that of the \demons{} algorithm (if we run the code on a single resolution level; see \figref{f:convergence-nirep-claire-vs-demons} in \secref{s:reg-quality}).
\item \claire{} introduces different acceleration schemes. These schemes not only stabilize the computations but also lead to a reduction in runtime (see \tabref{t:results-nirep-runtime} in \secref{s:time-to-solution}). \claire{} delivers a speedup of $5\times$ for the parameter continuation. We observed a speedup of up to $17\times$ when considering a grid continuation scheme (results not reported here). We disregarded this scheme, because we observed a significant dependence of the performance on the regularity of the velocity. Combining the grid and parameter continuation scheme may yield an even better performance. Designing an effective schedule for a combined scheme remains subject to future work.
\item Our two-level preconditioner is effective. We achieve the best performance if we compute the action of its inverse with a nested PCG method. This allows us to avoid a repeated estimation of spectral bounds of the reduced space Hessian operator, which is necessary if we consider a semi-iterative Chebyshev method. For real data, we achieve a moderate speedup of about $4\times$ for the entire inversion compared to our prior work~\cite{Mang:2016a}. Moreover, we saw that the performance of our schemes for preconditioning the reduced space Hessian is \emph{not} independent of the regularization parameter for the velocity. Designing a preconditioner that yields a good performance for vanishing regularity of the objects requires more work.
\item \claire{} delivers good scalability results. In this work, we showcase results for up to $\num{3221225472}$ unknowns on 342 compute nodes of TACC's Lonestar 5 system executed with 4096 MPI tasks. This demonstrates that we can tackle applications that require the registration of high-resolution imaging data such as, e.g., CLARITY imaging (a new optical imaging technique that delivers sub-micron resolution~\cite{Kutten:2017a,Chung:2013a,Tomer:2014a}). Further, we demonstrated that \claire{} can deliver runtimes that represent a significant step towards providing \iquote{real-time} capabilities for clinically relevant problem sizes (inversion for $\sim$50 million unknowns in about \SI{2}{\second} using 256 MPI tasks). To the best of our knowledge, \claire{} is the only software with these capabilities. We emphasize that \claire{} does not need to be executed on an HPC system; it can be executed on a standard compute node with a single core. Further runtime accelerations on limited hardware resources form the basis of our current work.
\end{itemize}

With this work we have identified several aspects of \claire{} that need to be improved. The time-to-solution on a single workstation is not yet fully competitive with the diffeomorphic \demons{} algorithm. We are currently working on improvements to our computational kernels to further reduce the execution time of \claire{}. In addition to algorithmic improvements, we are also actively working on a GPU implementation of \claire{}. In our scheme, we fix the parameter that controls the penalty on the divergence of the velocity; we only search for an adequate the regularization parameter for the velocity automatically (using a binary search). We found that this scheme works well in practice. Introducing this penalty not only yields better behaved deformation map (determinant of deformation gradient remains close to one) but also stabilizes the computations~\cite{Mang:2017b}. Designing an efficient method to automatically identify both parameters requires more work. As we have mentioned in the limitations, \claire{} does not support time dependent (nonstationary) velocities. We note that certain applications may benefit from nonstationary $\vect{v}$. In this work, we have demonstrated experimentally that if we are only interested in registering two images, stationary $\vect{v}$ produce good results. This is in accordance with observations  made in our past work~\cite{Mang:2015a,Mang:2016a,Mang:2017c,Scheufele:2019a,Gholami:2019a} as well as observations made by other groups~\cite{Arsigny:2006a,Ashburner:2007a,Hernandez:2009a,Lorenzi:2013a,Lorenzi:2013b,Vercauteren:2009a}. The design of efficient numerical schemes for nonstationary (time dependent) velocities is something we will address in our future work. Moreover, we are currently adding support for new distance measures to enable multi-modal registration.

\paragraph{Acknowledgements.} We would like to thank Anna-Lena Belgardt for suggesting the name \claire{}.

\FloatBarrier
\newpage

\setcounter{section}{0}
\renewcommand*{\thesection}{S\arabic{section}}
\renewcommand*{\thefigure}{S\arabic{figure}}
\renewcommand*{\thetable}{S\arabic{table}}

\begin{center}
\textbf{SUPPLEMENTARY MATERIAL\\CLAIRE: A DISTRIBUTED-MEMORY SOLVER FOR CONSTRAINED LARGE DEFORMATION DIFFEOMORPHIC IMAGE REGISTRATION}
\end{center}

In the following sections we provide a more detailed picture of the results reported in~\secref{s:experiments} of our manuscript. We also provide additional insight into the formal derivation and definition of the optimality conditions.

\section{Preconditioning}

We provide detailed results for the study of the performance of the preconditioner reported in \secref{s:precond-performance} of the main manuscript. We execute the runs CACDS's Opuntia system (see \secref{s:compute-systems} for the specs). We report results for an $H^2$ regularization model for $\beta_v\in\{\num{1e-2},\num{1e-3},\num{1e-4}\}$ in \tabref{t:convergence-nirep-krylov-h2s}. We report results for an $H^1$-div regularization model for $\beta_v\in\{\num{1e-1},\num{5e-2},\num{1e-2},\num{5e-3}\}$ in \tabref{t:convergence-nirep-krylov-h1s-div}. As a baseline, we consider the spectral preconditioner used in our prior work~\cite{Mang:2016c,Gholami:2017a}. We report the number of Hessian matvecs on the fine and the coarse grid (in brackets), the number of PDE solves (on the fine grid), the runtime in seconds and the speedup compared to the baseline method (spectral preconditioner).

\FloatBarrier

\begin{table}
\caption{Rate of convergence for the iterative inversion of the Hessian operator for different realizations of the preconditioner. We consider an $H^2$ regularization model (seminorm). We report results for the spectral preconditioner and different variants of the two-level preconditioner. We consider an inexact Chebyshev semi-iterative methods, CHEB($k$), with a fixed number of $k\in\{5,10,20\}$ iterations and a PCG method with a tolerance that is \num{1e-1} times smaller than the tolerance of the (outer) PCG method. The relative tolerance for the (outer) PCG method is set to $\num{1E-3}$. We report the number of Hessian matvecs on the fine (and the coarse) grid, the number of PDE solves (on the fine grid), and the runtime (in seconds). We report results for a grid size of $128\times150\times128$ (left columns) and $256\times300\times256$ (right columns). The results correspond to those reported in~\figref{f:convergence-nirep-krylov-h2s} in the main manuscript. We execute our solver on CACDS's Opuntia system (see \secref{s:compute-systems} of the main manuscript for the specs).\label{t:convergence-nirep-krylov-h2s}}
\resetrunid\centering\scriptsize
\begin{tabular}{lllrrlLlrrlL}\toprule
\multicolumn{7}{r}{$128\times150\times128$}                                                                &        & \multicolumn{4}{r}{$256\times300\times256$}                    \\
\midrule
$\beta_v$  & solver             &        & \#matvecs               & \#PDE & runtime            & speedup  &        & \#matvecs                    & \#PDE & runtime            & speedup \\
\midrule
\num{1e-2} & ---                & \runid &  26\;\;\phantom{(0000)} &  56   & \num{1.745817e+01} &          & \runid &  28\;\;\phantom{(0000)} &  60   & \num{1.847025e+02} &         \\
           & CHEB(\phantom{0}5) & \runid &  10\;\;(\phantom{00}70) &  24   & \num{1.351365e+01} & 1.3      & \runid &  10\;\;(\phantom{00}70) &  24   & \num{1.306160e+02} & 1.4     \\
           & CHEB(10)           & \runid &   7\;\;(\phantom{00}87) &  18   & \num{1.278577e+01} & 1.4      & \runid &   7\;\;(\phantom{00}87) &  28   & \num{1.187935e+02} & 1.6     \\
           & CHEB(20)           & \runid &   5\;\;(\phantom{0}115) &  14   & \num{1.304398e+01} & 1.3      & \runid &   5\;\;(\phantom{0}115) &  14   & \num{1.213128e+02} & 1.5     \\
           & PCG(\num{1e-1})    & \runid &   3\;\;(\phantom{0}131) &  10   & \num{1.287856e+01} & 1.4      & \runid &   2\;\;(\phantom{00}97) &   8   & \num{9.914814e+01} & 1.9     \\
\midrule
\num{1e-3} & ---                & \runid &  98\;\;\phantom{(0000)} & 200   & \num{6.424825e+01} &          & \runid & 100\;\;\phantom{(0000)} & 204   & \num{6.541756e+02} &         \\
           & CHEB(\phantom{0}5) & \runid &  32\;\;(\phantom{0}202) &  68   & \num{3.893008e+01} & 1.7      & \runid &  35\;\;(\phantom{0}220) &  74   & \num{4.178054e+02} & 1.6     \\
           & CHEB(10)           & \runid &  24\;\;(\phantom{0}274) &  52   & \num{3.970342e+01} & 1.6      & \runid &  25\;\;(\phantom{0}285) &  54   & \num{3.826560e+02} & 1.7     \\
           & CHEB(20)           & \runid &  18\;\;(\phantom{0}388) &  40   & \num{4.159431e+01} & 1.5      & \runid &  18\;\;(\phantom{0}388) &  40   & \num{3.931752e+02} & 1.7     \\
           & PCG(\num{1e-1})    & \runid &   4\;\;(\phantom{0}559) &  12   & \num{4.505822e+01} & 1.4      & \runid &   2\;\;(\phantom{0}311) &   8   & \num{2.232896e+02} & 2.9     \\
\midrule
\num{1e-4} & ---                & \runid & 347\;\;\phantom{(0000)} & 698   & \num{2.232703e+02} &          & \runid & 356\;\;\phantom{(0000)} & 716   & \num{2.272152e+03} &         \\
           & CHEB(\phantom{0}5) & \runid &  85\;\;(\phantom{0}520) & 174   & \num{1.022998e+02} & 2.2      & \runid & 112\;\;(\phantom{0}682) & 228   & \num{1.299544e+03} & 1.7     \\
           & CHEB(10)           & \runid &  63\;\;(\phantom{0}703) & 130   & \num{9.917461e+01} & 2.3      & \runid &  83\;\;(\phantom{0}923) & 170   & \num{1.235113e+03} & 1.8     \\
           & CHEB(20)           & \runid &  46\;\;(\phantom{0}976) &  96   & \num{1.026121e+02} & 2.2      & \runid &  60\;\;(1270)           & 124   & \num{1.289689e+03} & 1.8     \\
           & PCG(\num{1e-1})    & \runid &   4\;\;(1717)           &  12   & \num{1.304624e+02} & 1.7      & \runid &   3\;\;(1439)           &  10   & \num{9.476400e+02} & 2.4     \\
\bottomrule
\end{tabular}
\end{table}

\begin{table}
\caption{Rate of convergence for the iterative inversion of the Hessian operator for different realizations of the preconditioner. We consider an $H^1$-div regularization model ($H^1$-seminorm for $\vect{v}$; $\beta_w=\num{1E-4}$). We report results for the spectral preconditioner and different variants of the two-level preconditioner. We consider an inexact Chebyshev semi-iterative methods, CHEB($k$), with a fixed number of $k\in\{5,10,20\}$ iterations and a PCG method with a tolerance that is \num{1e-1} times smaller than the tolerance of the (outer) PCG method. The relative tolerance for the (outer) PCG method is set to $\num{1E-3}$. We report the number of Hessian matvecs on the fine (and the coarse) grid, the number of  PDE solves, and the runtime (in seconds). We report results for a grid size of $128\times150\times128$ (left columns) and $256\times300\times256$ (right columns). The results correspond to those reported in \figref{f:convergence-nirep-krylov-h1s-div} of the main manuscript. We execute our solver on CACDS's Opuntia system.\label{t:convergence-nirep-krylov-h1s-div}}
\resetrunid\centering\scriptsize
\begin{tabular}{lllrrlLlrrlL}\toprule
\multicolumn{7}{r}{$128\times150\times128$}                                                             &         & \multicolumn{4}{r}{$256\times300\times256$}                    \\\midrule
$\beta_v$  & solver             &        & \#matvecs             & \#PDE & runtime            & speedup &         & \#matvecs               & \#PDE & runtime            & speedup \\\midrule
\num{1e-1} &  ---               & \runid & 14\;\;\phantom{(000)} &  32   & \num{1.031646e+01} &         & \runid  &  15\;\;\phantom{(000)}  &  34   & \num{1.118666e+02} &         \\
           & CHEB(\phantom{0}5) & \runid &  4\;\;(\phantom{0}34) &  12   & \num{6.986409e+00} & 1.5     & \runid  &   4\;\;(\phantom{0}34)  &  12   & \num{6.775912e+01} & 1.7     \\
           & CHEB(10)           & \runid &  4\;\;(\phantom{0}54) &  12   & \num{8.444916e+00} & 1.2     & \runid  &   3\;\;(\phantom{0}43)  &  10   & \num{6.483268e+01} & 1.7     \\
           & CHEB(20)           & \runid &  3\;\;(\phantom{0}73) &  10   & \num{9.136916e+00} & 1.1     & \runid  &   2\;\;(\phantom{0}52)  &   8   & \num{6.172103e+01} & 1.8     \\
           & PCG(\num{1e-1})    & \runid &  3\;\;(\phantom{0}54) &  10   & \num{7.731989e+00} & 1.3     & \runid  &   2\;\;(\phantom{0}39)  &   8   & \num{5.441639e+01} & 2.1     \\
\midrule
\num{5e-2} & ---                & \runid & 19\;\;\phantom{(000)} &  42   & \num{1.359233e+01} &         & \runid  &  21\;\;\phantom{(000)}  &  46   & \num{1.511006e+02} &         \\
           & CHEB(\phantom{0}5) & \runid &  6\;\;(\phantom{0}46) &  16   & \num{9.404672e+00} & 1.4     & \runid  &   5\;\;(\phantom{0}40)  &  14   & \num{7.742016e+01} & 2.0     \\
           & CHEB(10)           & \runid &  5\;\;(\phantom{0}65) &  14   & \num{1.007616e+01} & 1.3     & \runid  &   4\;\;(\phantom{0}54)  &  12   & \num{7.999322e+01} & 1.9     \\
           & CHEB(20)           & \runid &  4\;\;(\phantom{0}94) &  12   & \num{1.145486e+01} & 1.2     & \runid  &   3\;\;(\phantom{0}73)  &  10   & \num{8.409937e+01} & 1.8     \\
           & PCG(\num{1e-1})    & \runid &  4\;\;(100)           &  12   & \num{1.201304e+01} & 1.1     & \runid  &   2\;\;(\phantom{0}55)  &   8   & \num{6.452471e+01} & 2.3     \\
\midrule
\num{1e-2} & ---                & \runid & 44\;\;\phantom{(000)} &  92   & \num{3.001269e+01} &         & \runid  &  54\;\;\phantom{(000)}  & 112   & \num{3.758063e+02} &         \\
           & CHEB(\phantom{0}5) & \runid & 15\;\;(100)           &  34   & \num{2.022719e+01} & 1.5     & \runid  &  12\;\;(\phantom{0}82)  &  28   & \num{1.632247e+02} & 2.3     \\
           & CHEB(10)           & \runid & 11\;\;(131)           &  26   & \num{1.951458e+01} & 1.5     & \runid  &   9\;\;(109)            &  22   & \num{1.570416e+02} & 2.4     \\
           & CHEB(20)           & \runid &  8\;\;(178)           &  20   & \num{2.067585e+01} & 1.5     & \runid  &   7\;\;(157)            &  18   & \num{1.717738e+02} & 2.2     \\
           & PCG(\num{1e-1})    & \runid &  5\;\;(279)           &  10   & \num{2.636161e+01} & 1.1     & \runid  &   3\;\;(183)            &  10   & \num{1.568739e+02} & 2.4     \\
\midrule
\num{5e-3} & ---                & \runid & 64\;\;\phantom{(000)} & 132   & \num{4.355236e+01} &         & \runid  &  80\;\;\phantom{(000)}  & 164   & \num{5.523839e+02} &         \\
           & CHEB(\phantom{0}5) & \runid & 25\;\;(160)           &  54   & \num{3.199040e+01} & 1.4     & \runid  &  20\;\;(130)            &  44   & \num{2.584596e+02} & 2.1     \\
           & CHEB(10)           & \runid & 18\;\;(208)           &  40   & \num{3.050017e+01} & 1.4     & \runid  &  14\;\;(164)            &  32   & \num{2.317423e+02} & 2.4     \\
           & CHEB(20)           & \runid & 13\;\;(283)           &  30   & \num{3.247983e+01} & 1.3     & \runid  &  10\;\;(220)            &  24   & \num{2.361081e+02} & 2.3     \\
           & PCG(\num{1e-1})    & \runid &  7\;\;(556)           &  18   & \num{4.830415e+01} & 0.9     & \runid  &   3\;\;(262)            &  10   & \num{2.097867e+02} & 2.6     \\
\bottomrule
\end{tabular}
\end{table}

We refer to \secref{s:precond-performance} of the main manuscript for a discussion of these results.

\section{Convergence: Newton--Krylov Solver (Smooth Data)}

We augment the convergence results reported in \secref{s:convergence-nks} for real data with results for a synthetic test problem (see above). The problem is discretized on a grid of size $256^3$. We consider an $H^2$-regularization model (seminorm; $\beta_v = \num{1e-4}$) and an $H^1$-div regularization model ($H^1$-seminorm for $\vect{v}$; $\beta_v=\num{1E-2}$, $\beta_w = \num{1E-4}$). We run the registration at full resolution (\num{50331648} unknowns). The number of Newton iterations is limited to 50 (not reached). We use a superlinear forcing sequence and limit the number of Krylov iterations to 100 (not reached). The tolerance for the relative change of the gradient is $\num{1e-03}$; the absolute tolerance for the norm of the gradient is $\num{1e-06}$. The number of time steps for the PDE solves is set to $n_t=4$. We use 20 cores (64GB compute nodes) resulting in a processor layout of $5\times4$ ($\sim\num{2555904}$ unknowns per core). We do not perform any parameter, scale, or grid continuation. The same synthetic test problem is also considered for the scalability study of our solver in \secref{s:scalability} of the main manuscript.

\begin{figure}
\centering
\includegraphics[width=0.99\textwidth]
{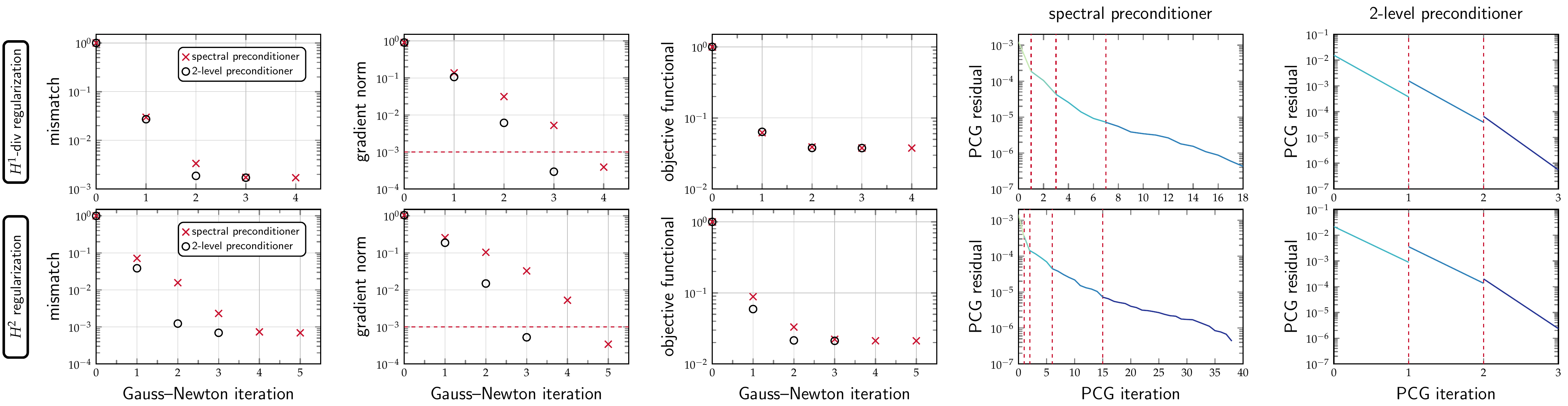}
\caption{Convergence of Newton--Krylov solver for a synthetic test problem for different regularization norms and variants for the preconditioner. The top row shows results for an $H^1$-div regularization model ($H^1$-seminorm for $\vect{v}$ with $\beta_v=\num{1E-2}$ and $\beta_w = \num{1E-4}$). The bottom row shows results for an $H^2$ regularization model (seminorm; $\beta_v = \num{1e-4}$). We plot (from left to right) the relative reduction of ($i$) the mismatch ($L^2$-distance between the images to be registered), ($ii$) the reduced gradient, and ($iii$) the objective functional, with respect to the Gauss--Newton iterations. We use a relative change of the gradient by $\num{1e-3}$ as a stopping criterion (dashed red line in second column). The two plots on the right show the convergence of the PCG solver per Gauss--Newton iteration for different realization of the preconditioner, respectively. The vertical, dashed red lines separate the individual Gauss--Newton iterations; the PCG iteration index is cumulative.\label{f:convergence-syn-fullsolve}}
\end{figure}

We report results in \figref{f:convergence-syn-fullsolve}. The top row shows results for the $H^1$-div regularization model and the bottom row for the $H^2$ regularization model. We plot the relative reduction of the mismatch, the reduced gradient, and the objective functional with respect to the Gauss--Newton iteration index. We also report results for the convergence of the PCG solver for different realizations of the preconditioner.

We can see that Newton--Krylov solver converges after only 3 to 4 Gauss--Newton iterations. We can reduce the gradient by three orders of magnitude in less than 5 Gauss--Newton iterations (about one order of magnitude per Gauss--Newton iteration if we consider a two-level preconditioner). We require one Hessian matvec per Gauss--Newton iteration for the two-level preconditioner in combination with a nested PCG method. The residual in the PCG method drops rapidly for both preconditioners. We observe a better search direction per Gauss--Newton iteration if we consider the two-level preconditioner in combination with a nested PCG method (we slightly oversolve the KKT system for the nested preconditioner). The reduced gradient drops more rapidly. The trend of the objective and the mismatch is more or less similar.

\section{Convergence: Newton--Krylov Solver (Real Data)}

We report detailed results for the convergence study of \claire{}\,'s Newton--Krylov solver in \secref{s:convergence-nks} of the main manuscript. We compare the performance of our new and improved solver to the performance of the solver used in our past work~\cite{Mang:2016c,Gholami:2017a}. We execute the runs TACC's Lonestar 5 system (see \secref{s:compute-systems} of the main manuscript for the specs). We consider all datasets of the \nirep{} repository (see \secref{s:data} of the main manuscript). As a baseline, we report results for the spectral preconditioner in \tabref{t:nirep-convergence-h1sdiv-spectral-ls5} (double precision). We report results for the two-level preconditioner in \tabref{t:nirep-convergence-h1sdiv-2level-pcg0d1-dbl-ls5} (double precision) and \tabref{t:nirep-convergence-h1sdiv-2level-pcg0d1-sgl-ls5} (single precision). We report the number of Newton iterations, the number of Hessian matvecs (on the fine and the coarse grid), and the number of PDE solves (on the fine grid) until convergence. We also report the relative mismatch, the absolute $\ell^2$-norm of the gradient (after registration), the relative change of the $\ell^2$-norm of the gradient, and the runtime (in seconds). For the runs for the two-level preconditioner we also report the achieved speedup with respect to each individual dataset. A direct comparison of our single and double implementation can be found in \tabref{t:dbl-vs-sgl}.

\begin{table}
\caption{Computational performance for the \nirep{} data for \claire{}. We consider the spectral preconditioner (inverse of the regularization operator). We consider an $H^1$-div regularization model ($H^1$-seminorm for $v$) with $\beta_v=\num{1e-2}$ and $\beta_w=\num{1e-4}$. We execute the runs on TACC's Lonestar 5 system in double precision (see \secref{s:compute-systems} of the main manuscript for the specs). We terminate the inversion if the gradient is reduced by \num{5e-2}. For each image registration pair {\tt na02} through {\tt na16}  to {\tt na01}, we report (i) the number of Gauss--Newton iterations until convergence, (ii) the number of Hessian matvecs, (iii) the number of PDE solves, (iv) the relative reduction of the mismatch, (v) the $\ell^2$-norm of the gradient after registration, (vi) the relative change of the $\ell^2$-norm of the gradient $\|\di{g}^\star\|_{\text{rel}}\defeq\|\di{g}^\star\|_2/\|\di{g}_0\|_2$, and (vii) the runtime (in seconds).\label{t:nirep-convergence-h1sdiv-spectral-ls5}}
\resetrunid\centering\scriptsize
\begin{tabular}{llrrrlllL}\toprule
       &               & \#iter & \#matvecs & \#PDE & mismatch                 & $\|\di{g}^\star\|_2$     & $\|\di{g}^\star\|_{\text{rel}}$ & runtime            \\\midrule
\runid & \texttt{na02} &  6     &  44       & 103   & \num{1.049965142958E-01} & \num{4.357365100927E-04} & \num{4.631785816778E-02}        & \num{2.376812e+02} \\
\runid & \texttt{na03} &  9     &  83       & 187   & \num{8.465094487460E-02} & \num{4.630653548871E-04} & \num{4.710900411598E-02}        & \num{4.440620e+02} \\
\runid & \texttt{na04} &  8     &  63       & 145   & \num{9.183222046396E-02} & \num{4.507574446373E-04} & \num{4.492712886368E-02}        & \num{3.521112e+02} \\
\runid & \texttt{na05} &  6     &  46       & 107   & \num{1.343885658861E-01} & \num{4.077175630690E-04} & \num{4.313448048189E-02}        & \num{2.502169e+02} \\
\runid & \texttt{na06} &  8     &  54       & 127   & \num{9.688378794386E-02} & \num{3.963745677822E-04} & \num{4.350073752212E-02}        & \num{2.988379e+02} \\
\runid & \texttt{na07} &  8     &  62       & 143   & \num{8.467089519133E-02} & \num{4.741323056129E-04} & \num{4.869817751748E-02}        & \num{3.410715e+02} \\
\runid & \texttt{na08} &  6     &  39       &  93   & \num{9.607811691334E-02} & \num{4.352931161878E-04} & \num{4.376957034245E-02}        & \num{2.385397e+02} \\
\runid & \texttt{na09} &  7     &  50       & 117   & \num{9.000616126230E-02} & \num{4.383875056436E-04} & \num{4.544159952194E-02}        & \num{2.811073e+02} \\
\runid & \texttt{na10} &  9     &  63       & 147   & \num{6.153417125619E-02} & \num{4.949717372295E-04} & \num{4.839976709155E-02}        & \num{3.563325e+02} \\
\runid & \texttt{na11} & 13     & 126       & 281   & \num{7.383833811430E-02} & \num{4.880549933022E-04} & \num{4.853651215705E-02}        & \num{6.918718e+02} \\
\runid & \texttt{na12} & 10     &  93       & 209   & \num{1.075030116929E-01} & \num{4.837212809706E-04} & \num{4.968675767516E-02}        & \num{4.997932e+02} \\
\runid & \texttt{na13} &  7     &  52       & 121   & \num{9.670308556120E-02} & \num{4.977777085935E-04} & \num{4.907589446181E-02}        & \num{2.832756e+02} \\
\runid & \texttt{na14} & 15     & 154       & 347   & \num{5.886181254449E-02} & \num{5.125269344282E-04} & \num{4.908259033548E-02}        & \num{8.367538e+02} \\
\runid & \texttt{na15} &  8     &  59       & 137   & \num{7.740261422476E-02} & \num{4.269922279308E-04} & \num{4.403600200771E-02}        & \num{3.242809e+02} \\
\runid & \texttt{na16} &  7     &  56       & 129   & \num{1.036062493205E-01} & \num{4.517379945876E-04} & \num{4.604843203529E-02}        & \num{3.034767e+02} \\\bottomrule
\end{tabular}
\end{table}

\begin{table}
\caption{Computational performance for the \nirep{} data for \claire{}. We consider the 2-level preconditioner with PCG(\num{1E-1}) as a solver. We execute the runs on TACC's Lonestar 5 system in double precision (see \secref{s:compute-systems} of the main manuscript for the specs). We consider an $H^1$-div regularization model ($H^1$-seminorm for $v$) with $\beta_v=\num{1e-2}$ and $\beta_w=\num{1e-4}$.  We terminate the inversion if the gradient is reduced by \num{5e-2}. For each image registration pair na$i$ to na01, we report (i) the number of Gauss--Newton iterations until convergence, (ii) the number of Hessian matvecs (the number of matvecs for the iterative inversion of the preconditioner is reported in brackets), (iii) the number of PDE solves, (iv) the relative reduction of the mismatch, (v) the $\ell^2$-norm of the gradient after registration, (vi) the relative change of the $\ell^2$-norm of the gradient $\|\di{g}^\star\|_{\text{rel}}\defeq\|\di{g}^\star\|_2/\|\di{g}_0\|_2$, (vii) the min, mean, and max values for $\det\igrad\vect{y}^{-1}$, and (viii) the runtime in seconds (the speedup we obtain for the entire inversion compared to the results reported in \tabref{t:nirep-convergence-h1sdiv-spectral-ls5} is given in brackets).\label{t:nirep-convergence-h1sdiv-2level-pcg0d1-dbl-ls5}}
\resetrunid\centering\scriptsize
\begin{tabular}{llrrrlllLR}\toprule
       &               & \#iter & \#matvecs    & \#PDE & mismatch                 & $\|\di{g}^\star\|_2$     & $\|\di{g}^\star\|_{\text{rel}}$ & runtime            & speedup \\\midrule
\runid & \texttt{na02} &  7     &  7\;\;(153)  & 31    & \num{1.031853094846E-01} & \num{3.887923940403E-04} & \num{4.132779913265E-02}        & \num{1.481538e+02} & 1.6     \\
\runid & \texttt{na03} &  9     &  9\;\;(190)  & 39    & \num{8.597876018175E-02} & \num{4.653556640269E-04} & \num{4.734200401881E-02}        & \num{1.876604e+02} & 2.4     \\
\runid & \texttt{na04} &  8     &  8\;\;(164)  & 35    & \num{9.186817751865E-02} & \num{4.867450449800E-04} & \num{9.186817751865E-02}        & \num{1.625234e+02} & 2.2     \\
\runid & \texttt{na05} &  6     &  6\;\;(132)  & 27    & \num{1.344639401091E-01} & \num{4.133006310169E-04} & \num{4.372514116772E-02}        & \num{1.260529e+02} & 2.0     \\
\runid & \texttt{na06} &  8     &  8\;\;(150)  & 35    & \num{9.758281342349E-02} & \num{4.554048537695E-04} & \num{4.997910718887E-02}        & \num{1.567854e+02} & 1.9     \\
\runid & \texttt{na07} & 10     & 10\;\;(210)  & 43    & \num{8.142204007764E-02} & \num{4.656481631524E-04} & \num{4.782677037915E-02}        & \num{2.100570e+02} & 1.6     \\
\runid & \texttt{na08} &  6     &  6\;\;(117)  & 27    & \num{9.467525923614E-02} & \num{4.093837635766E-04} & \num{4.116433449223E-02}        & \num{1.189912e+02} & 2.0     \\
\runid & \texttt{na09} &  7     &  7\;\;(135)  & 31    & \num{8.914189823343E-02} & \num{4.395306810622E-04} & \num{4.556009678495E-02}        & \num{1.477237e+02} & 1.9     \\
\runid & \texttt{na10} &  9     &  9\;\;(164)  & 39    & \num{6.073033658971E-02} & \num{4.858411318000E-04} & \num{4.750695010230E-02}        & \num{1.803303e+02} & 2.0     \\
\runid & \texttt{na11} & 13     & 13\;\;(269)  & 55    & \num{7.518610185000E-02} & \num{4.840489165193E-04} & \num{4.813811239239E-02}        & \num{2.697554e+02} & 2.6     \\
\runid & \texttt{na12} &  9     &  9\;\;(191)  & 39    & \num{1.089612545918E-01} & \num{4.562418838900E-04} & \num{4.686413605912E-02}        & \num{1.896031e+02} & 2.6     \\
\runid & \texttt{na13} &  7     &  7\;\;(144)  & 31    & \num{9.619024641798E-02} & \num{4.989570803847E-04} & \num{4.919216870342E-02}        & \num{1.446507e+02} & 2.0     \\
\runid & \texttt{na14} & 14     & 14\;\;(303)  & 59    & \num{5.862523626273E-02} & \num{5.022851847551E-04} & \num{4.810178021652E-02}        & \num{3.000502e+02} & 2.8     \\
\runid & \texttt{na15} &  8     &  8\;\;(150)  & 35    & \num{7.725422907455E-02} & \num{4.516550497233E-04} & \num{4.657949577395E-02}        & \num{1.576384e+02} & 2.1     \\
\runid & \texttt{na16} &  8     &  8\;\;(174)  & 35    & \num{1.010861334171E-01} & \num{4.152152845790E-04} & \num{4.232544758473E-02}        & \num{1.728787e+02} & 1.8     \\\bottomrule
\end{tabular}
\end{table}

\begin{table}
\caption{Computational performance for the \nirep{} data for \claire{}. We consider the 2-level preconditioner with PCG(\num{1E-1}) as a solver. We execute the runs on TACC's Lonestar 5 system in single precision (see \secref{s:compute-systems} of the main manuscript for the specs). We consider an $H^1$-div regularization model ($H^1$-seminorm for $v$) with $\beta_v=\num{1e-2}$ and $\beta_w=\num{1e-4}$. We terminate the inversion if the gradient is reduced by \num{5e-2}. For each image registration pair na$i$ to na01, we report (i) the number of Gauss--Newton iterations until convergence, (ii) the number of Hessian matvecs (the number of matvecs for the iterative inversion of the preconditioner is reported in brackets), (iii) the number of PDE solves, (iv) the relative reduction of the mismatch, (v) the $\ell^2$-norm of the gradient after registration, (vi) the relative change of the $\ell^2$-norm of the gradient $\|\di{g}^\star\|_{\text{rel}}\defeq\|\di{g}^\star\|_2/\|\di{g}_0\|_2$, (vii) the min, mean, and max values for $\det\igrad\vect{y}^{-1}$, and (viii) the runtime in seconds (the speedup we obtain for the entire inversion compared to the results reported in \tabref{t:nirep-convergence-h1sdiv-spectral-ls5} is given in brackets).\label{t:nirep-convergence-h1sdiv-2level-pcg0d1-sgl-ls5}}
\resetrunid\centering\scriptsize
\begin{tabular}{llrrrlllLR}\toprule
       &               & \#iter & \#matvecs   & \#PDE & mismatch                 & $\|\di{g}^\star\|_2$     & $\|\di{g}^\star\|_{\text{rel}}$ & runtime            & speedup \\\midrule
\runid & \texttt{na02} &  7     &  7\;\;(153) & 31    & \num{1.031870767474E-01} & \num{3.886843041983E-04} & \num{4.132203385234E-02}        & \num{6.915890e+01} & 3.4     \\
\runid & \texttt{na03} &  9     &  9\;\;(190) & 39    & \num{8.597768843174E-02} & \num{4.652219358832E-04} & \num{4.733553156257E-02}        & \num{8.779980e+01} & 5.1     \\
\runid & \texttt{na04} &  8     &  8\;\;(164) & 35    & \num{9.186743944883E-02} & \num{4.866061499342E-04} & \num{4.850723221898E-02}        & \num{7.615361e+01} & 4.6     \\
\runid & \texttt{na05} &  6     &  6\;\;(133) & 27    & \num{1.344650089741E-01} & \num{4.133087350056E-04} & \num{4.373261705041E-02}        & \num{5.929297e+01} & 4.2     \\
\runid & \texttt{na06} &  8     &  8\;\;(150) & 35    & \num{9.757937490940E-02} & \num{4.552870814223E-04} & \num{4.997305572033E-02}        & \num{7.451897e+01} & 4.0     \\
\runid & \texttt{na07} & 10     & 10\;\;(210) & 43    & \num{8.142353594303E-02} & \num{4.655551747419E-04} & \num{4.782346263528E-02}        & \num{9.816511e+01} & 3.5     \\
\runid & \texttt{na08} &  6     &  6\;\;(117) & 27    & \num{9.467566758394E-02} & \num{4.092739545740E-04} & \num{4.115952551365E-02}        & \num{5.589632e+01} & 4.3     \\
\runid & \texttt{na09} &  7     &  7\;\;(135) & 31    & \num{8.914212882519E-02} & \num{4.394194111228E-04} & \num{4.555451497436E-02}        & \num{8.618371e+01} & 3.3     \\
\runid & \texttt{na10} &  9     &  9\;\;(164) & 39    & \num{6.072994694114E-02} & \num{4.857277090196E-04} & \num{4.750218242407E-02}        & \num{8.248141e+01} & 4.3     \\
\runid & \texttt{na11} & 13     & 13\;\;(268) & 55    & \num{7.518398016691E-02} & \num{4.839375033043E-04} & \num{4.813356697559E-02}        & \num{1.264818e+02} & 5.5     \\
\runid & \texttt{na12} &  9     &  9\;\;(191) & 39    & \num{1.089613735676E-01} & \num{4.560852830764E-04} & \num{4.685466736555E-02}        & \num{8.927655e+01} & 5.6     \\
\runid & \texttt{na13} &  7     &  7\;\;(144) & 31    & \num{9.619294852018E-02} & \num{4.988027503714E-04} & \num{4.918412491679E-02}        & \num{6.829672e+01} & 4.1     \\
\runid & \texttt{na14} & 14     & 14\;\;(303) & 59    & \num{5.862844735384E-02} & \num{5.021372344345E-04} & \num{4.809351265430E-02}        & \num{1.399514e+02} & 6.0     \\
\runid & \texttt{na15} &  8     &  8\;\;(150) & 35    & \num{7.725407183170E-02} & \num{4.515327163972E-04} & \num{4.657289758325E-02}        & \num{7.434582e+01} & 4.4     \\
\runid & \texttt{na16} &  8     &  8\;\;(174) & 35    & \num{1.010846123099E-01} & \num{4.150946624577E-04} & \num{4.231951385736E-02}        & \num{8.013536e+01} & 3.8     \\\bottomrule
\end{tabular}
\end{table}

\begin{table}
\caption{Exemplary results for the performance of our solver using double (\SI{64}{\bit}) and single (\SI{32}{\bit}) precision. We consider an two dataset of the \nirep{} repository---{\tt na02} and {\tt na01}. We use an $H^1$-div regularization model ($H^1$-seminorm for $\vect{v}$ with $\beta_v=\num{1e-2}$ and $\beta_w=\num{1e-4}$). We perform these runs on one and eight nodes of TACC's Lonestar 5 system (see \secref{s:compute-systems} of the main manuscript for the specs). We report (from left to right) the number of Gauss--Newton iterations, the number of Hessian matvecs, the relative mismatch after registration, the norm of the reduced gradient after registration, the relative reduction of the norm of the reduced gradient, the runtime, and the speedup (when switching from double to single precision).\label{t:dbl-vs-sgl}}
\resetrunid\centering\scriptsize
\begin{tabular}{lllrrrrrllllr}\toprule
       & solver          &               & nodes & tasks & \#iter & \#matvecs    & \#PDE & mismatch                 & $\|\di{g}^\star\|_2$     & $\|\di{g}^\star\|_{\text{rel}}$ & runtime            & speedup  \\
\midrule
\runid & ---             & \SI{64}{\bit} & 1     &  24   &  6     & 44           & 103   & \num{1.049965142958E-01} & \num{4.357365100927E-04} & \num{4.631785816778E-02}        & \num{2.424476e+02} &          \\
\runid &                 & \SI{32}{\bit} &       &       &  6     & 44           & 103   & \num{1.049939095974E-01} & \num{4.355806158856E-04} & \num{4.630769789219E-02}        & \num{1.123044e+02} & 2.2      \\
\runid &                 & \SI{64}{\bit} & 8     & 192   &  6     & 44           & 103   & \num{1.049965142958E-01} & \num{4.357365100927E-04} & \num{4.631785816778E-02}        & \num{2.993062e+01} &          \\
\runid &                 & \SI{32}{\bit} &       &       &  6     & 44           & 103   & \num{1.049958541989E-01} & \num{4.357126017567E-04} & \num{4.631592333317E-02}        & \num{1.766337e+01} & 1.7      \\
\midrule
\runid & PCG(\num{1e-1}) & \SI{64}{\bit} & 1     &  24   &  7     & 7~(153)      &  31   & \num{1.031853094846E-01} & \num{3.887923940403E-04} & \num{4.132779913265E-02}        & \num{1.478262e+02} &          \\
\runid &                 & \SI{32}{\bit} &       &       &  7     & 7~(153)      &  31   & \num{1.031870990992E-01} & \num{3.886843624059E-04} & \num{4.132203757763E-02}        & \num{6.951523e+01} & 2.1      \\
\runid &                 & \SI{64}{\bit} & 8     & 192   &  7     & 7~(153)      &  31   & \num{1.031853094846E-01} & \num{3.887923940403E-04} & \num{4.132779913265E-02}        & \num{2.678688e+01} &          \\
\runid &                 & \SI{32}{\bit} &       &       &  7     & 7~(153)      &  31   & \num{1.031854003668E-01} & \num{3.887765633408E-04} & \num{4.132666066289E-02}        & \num{1.105036e+01} & 2.4      \\
\bottomrule
\end{tabular}
\end{table}

We refer to \secref{s:convergence-nks} of the main manuscript for a discussion of these results.

\section{Registration Quality}

\begin{figure}
\centering
\includegraphics[width=0.9\textwidth]
{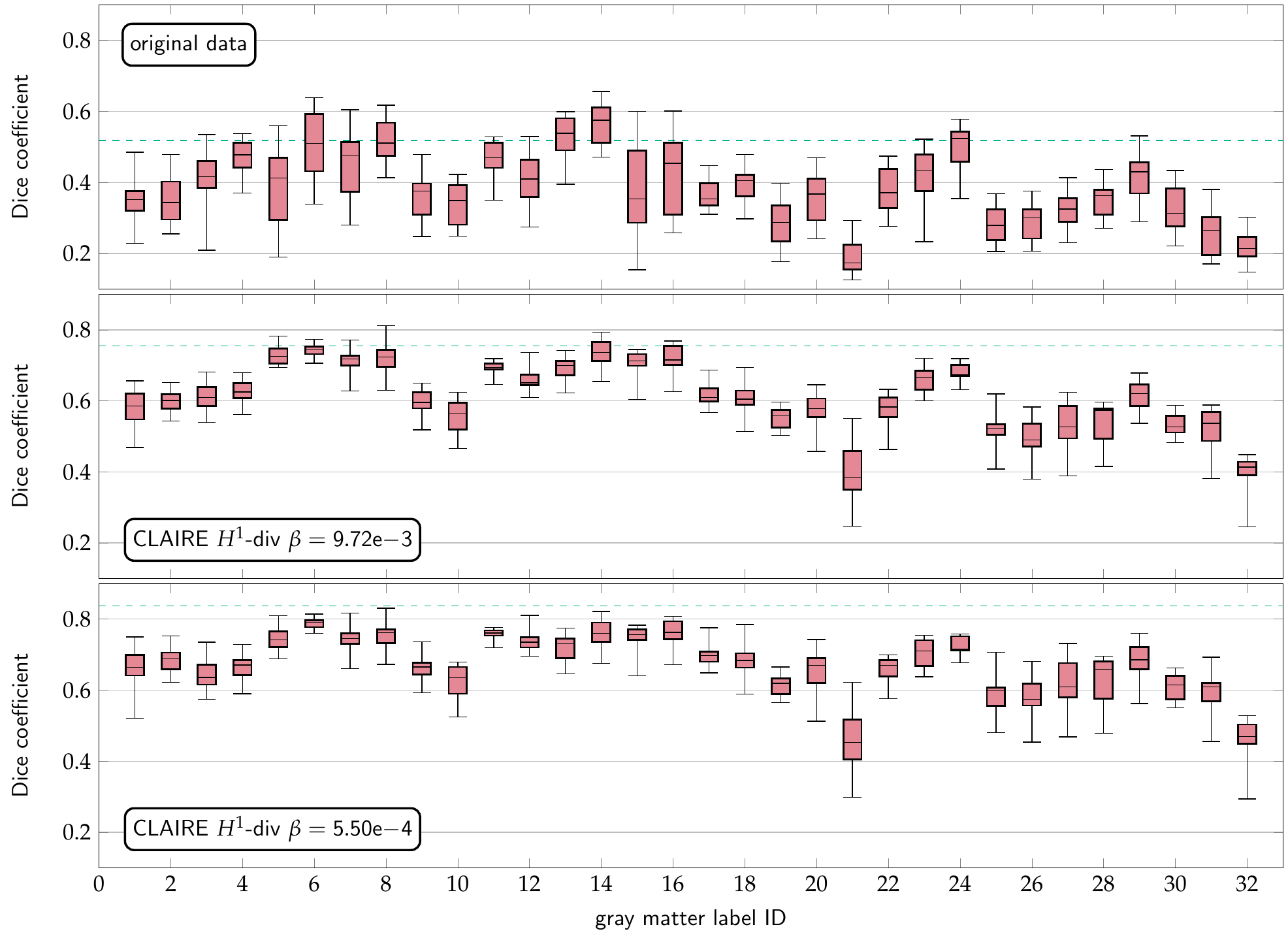}
\caption{Registration quality for the \nirep{} data for \claire{}. We report box plots for the Dice coefficient for the 32 individual gray matter labels averaged across 15 image pairs ({\tt na02} through {\tt na16} versus {\tt na01}; see \figref{f:nirep-data} of the main manuscript for an illustration). The box plots at the top represent the Dice coefficients before registration. The box plots in the middle and at the bottom represent the Dice coefficient after registration with \claire{}. We consider an $H^1$-div regularization model ($H^1$-seminorm for $\vect{v}$) with $\beta_w=\num{1e-4}$. The box plots in the middle correspond to results obtained for a regularization parameter $\beta_v =\num{9.718750e-03}$ (see also \tabref{t:nirep-regquality-h1sdiv-3d3e-3}); the box plots at the bottom correspond to results obtained for $\beta_v =\num{5.500000e-04}$ (see also \tabref{t:nirep-regquality-h1sdiv-5d5e-4}). The regularization parameters $\beta_v$ are determined via a binary search.\label{f:claire-quality-individual-labels}}
\end{figure}

We report additional results for the study of registration quality for the \nirep{} dataset for \claire{} and different variants of the \demons{} algorithm reported in \secref{s:reg-quality} of the main manuscript. We determine the regularization parameters for \claire{} using a binary search (see \secref{s:reg-quality} in the main manuscript for details on the implementation). We report results for \claire{} for $\beta_v =\num{9.718750e-03}$ in \tabref{t:nirep-regquality-h1sdiv-3d3e-3} and for $\beta_v =\num{5.500000e-04}$ in \tabref{t:nirep-regquality-h1sdiv-5d5e-4}. We report the relative change of the mismatch after registration. We also report the Dice score (before and after registration) and the false positive and false negative rate (after registration). These overlap scores are evaluated for the union of the 32 gray matter labels (for simplicity). In addition to that, we also provide the relative change of the $\ell^2$-norm of the reduced gradient, the extremal values for the determinant of the deformation gradient, and the runtime (in seconds). We report Dice scores for the 32 individual gray matter labels in \figref{f:claire-quality-individual-labels} (before and after registration).

The results for different variants of the diffeomorphic \demons{} algorithm are reported in \tabref{t:ddem:sweep}, \tabref{t:lddem:sweep}, and \tabref{t:demons-all-nirep-data}. We use a multi-resolution approach with 15, 10, and 5 iterations per level (default setting). We report registration quality as a function of the regularization weights $(\sigma_u, \sigma_d)$ for the DDEM algorithm in \tabref{t:ddem:sweep} ($\sigma_u$: smoothing for the updated field; $\sigma_d$: smoothing for the deformation field; units: voxel size along each spatial direction). We consider the \emph{diffeomorphic update rule} with forces based on the \emph{gradient of the deformed template image} (default method; left block in \tabref{t:ddem:sweep}) and \emph{symmetrized forces} (right block in \tabref{t:ddem:sweep}). We report registration quality as a function of the regularization parameters ($\sigma_u$, $\sigma_v$) for the LDDDEM algorithm in \tabref{t:lddem:sweep} ($\sigma_u$: smoothing for the updated field; $\sigma_v$: smoothing for the velocity field; units: voxel size along each spatial direction). We consider the \emph{log-domain update rule} with forces based on the \emph{gradient of the deformed template image} (LDDDEM; left block in \tabref{t:lddem:sweep}) and the \emph{symmetric log-domain update rule} with \emph{symmetrized forces} (SLDDDEM; default method; right block in \tabref{t:lddem:sweep}). For each variant of the \demons{} algorithm, we choose the regularization parameters that yield the highest Dice score (subject to the map being diffeomorphic as judged by the reported values for the determinant of the deformation gradient; we use the values reported by the \demons{} implementation). We refine the parameter search for the best \demons{} variants and identified parameters in \tabref{t:dem-iterations}. We additionally increase the number of iterations by a factor of 2, 5, 10, and 100. We can see that increasing the iteration count yields slightly better results. Based on these runs we found that SDDEM seems to give us the best results. We apply this method to the entire \nirep{} dataset.  We report these results in \tabref{t:demons-all-nirep-data}.

The results for the best variants of \demons{} and \claire{} across all \nirep{} datasets are summarized in \figref{f:claire-vs-demons-quality-nirep} in the main manuscript.

\begin{table}
\caption{Registration quality for the \nirep{} data for \claire{}. We consider an $H^1$-div regularization model ($H^1$-seminorm for $\vect{v}$) with $\beta_w=\num{1e-4}$. The regularization parameter $\beta_v$ is determined using a binary search with a bound of 0.3 for the determinant of the deformation gradient ($\beta_v =\num{9.718750e-03}$). We terminate the registration if the gradient is reduced by a factor of $\num{5e-2}$. We report (from left to right) the relative mismatch, the Dice coefficient (before and after registration), the false positive rate (after registration), the false negative rate (after registration), the relative reduction of the gradient, the extremal values of the determinant of the deformation gradient, and the overall runtime (in seconds). We execute the registration on CACDS's Opuntia server in single precision (see \secref{s:compute-systems} of the main manuscript for the specs).\label{t:nirep-regquality-h1sdiv-3d3e-3}}
\resetrunid\centering\scriptsize
\begin{tabular}{llllLlllllll}\toprule
\mcol{3}{\null}                                   & \multicolumn{2}{l}{dice}                & \mcol{3}{\null}                                                           & \multicolumn{3}{l}{$\det\igrad\vect{y}$}\\
\midrule
\mcol{2}{\null}        & mismatch                 & before             & after              & FP                 & FN                 & $\|\di{g}^\star\|_{\text{rel}}$ & $\min$             & $\operatorname{mean}$ & $\max$             & runtime     \\
\midrule
\runid & \texttt{na02} & \num{1.073665544391E-01} & \num{5.542543e-01} & \num{7.684721e-01} & \num{2.588482e-01} & \num{2.144801e-01} & \num{3.210427984595E-02}        & \num{3.271372e-01} & \num{1.012741e+00}    & \num{3.335465e+00} & \num{1.044427e+02} \\
\runid & \texttt{na03} & \num{8.903279155493E-02} & \num{5.039992e-01} & \num{7.517919e-01} & \num{3.124516e-01} & \num{2.095144e-01} & \num{3.058781102300E-02}        & \num{3.953025e-01} & \num{1.018710e+00}    & \num{4.995612e+00} & \num{1.116478e+02} \\
\runid & \texttt{na04} & \num{1.033979952335E-01} & \num{5.243747e-01} & \num{7.492278e-01} & \num{3.749245e-01} & \num{1.764035e-01} & \num{4.582517594099E-02}        & \num{2.939801e-01} & \num{1.022643e+00}    & \num{1.202991e+01} & \num{1.053583e+02} \\
\runid & \texttt{na05} & \num{1.387822479010E-01} & \num{5.550893e-01} & \num{7.533487e-01} & \num{2.966581e-01} & \num{2.164323e-01} & \num{2.975129708648E-02}        & \num{3.370395e-01} & \num{1.013372e+00}    & \num{3.806668e+00} & \num{9.688109e+01} \\
\runid & \texttt{na06} & \num{9.993377327919E-02} & \num{5.605455e-01} & \num{7.576347e-01} & \num{3.477614e-01} & \num{1.780914e-01} & \num{3.378171846271E-02}        & \num{3.366310e-01} & \num{1.020194e+00}    & \num{5.135905e+00} & \num{1.078640e+02} \\
\runid & \texttt{na07} & \num{8.585526794195E-02} & \num{5.311216e-01} & \num{7.643989e-01} & \num{3.255284e-01} & \num{1.799680e-01} & \num{2.912520989776E-02}        & \num{2.895173e-01} & \num{1.021152e+00}    & \num{3.603949e+00} & \num{1.089435e+02} \\
\runid & \texttt{na08} & \num{1.040807738900E-01} & \num{5.622983e-01} & \num{7.562517e-01} & \num{3.339812e-01} & \num{1.888828e-01} & \num{3.919803351164E-02}        & \num{2.635516e-01} & \num{1.014701e+00}    & \num{4.251909e+00} & \num{8.680117e+01} \\
\runid & \texttt{na09} & \num{9.102424234152E-02} & \num{5.077283e-01} & \num{7.540689e-01} & \num{3.123978e-01} & \num{2.057038e-01} & \num{2.864133380353E-02}        & \num{4.100031e-01} & \num{1.026398e+00}    & \num{6.399058e+00} & \num{1.122661e+02} \\
\runid & \texttt{na10} & \num{6.131324172020E-02} & \num{4.789702e-01} & \num{7.535548e-01} & \num{4.063941e-01} & \num{1.497460e-01} & \num{3.062749095261E-02}        & \num{4.692884e-01} & \num{1.032215e+00}    & \num{5.097562e+00} & \num{1.252207e+02} \\
\runid & \texttt{na11} & \num{7.730022817850E-02} & \num{4.579651e-01} & \num{7.468932e-01} & \num{3.774264e-01} & \num{1.790082e-01} & \num{3.116637282073E-02}        & \num{4.206548e-01} & \num{1.031427e+00}    & \num{1.182255e+01} & \num{1.272832e+02} \\
\runid & \texttt{na12} & \num{1.147219911218E-01} & \num{5.210443e-01} & \num{7.442749e-01} & \num{3.864193e-01} & \num{1.782619e-01} & \num{3.488916531205E-02}        & \num{3.475736e-01} & \num{1.021312e+00}    & \num{1.387510e+01} & \num{1.019495e+02} \\
\runid & \texttt{na13} & \num{9.955763071775E-02} & \num{5.307101e-01} & \num{7.518865e-01} & \num{3.646187e-01} & \num{1.779287e-01} & \num{3.214709088206E-02}        & \num{1.187021e-01} & \num{1.021972e+00}    & \num{7.354360e+00} & \num{1.033011e+02} \\
\runid & \texttt{na14} & \num{6.191835179925E-02} & \num{4.379558e-01} & \num{7.617315e-01} & \num{3.214263e-01} & \num{1.871133e-01} & \num{2.894376218319E-02}        & \num{3.298028e-01} & \num{1.034970e+00}    & \num{4.132410e+00} & \num{1.285522e+02} \\
\runid & \texttt{na15} & \num{8.581355959177E-02} & \num{4.987260e-01} & \num{7.420197e-01} & \num{3.583539e-01} & \num{1.987750e-01} & \num{3.940520435572E-02}        & \num{3.761488e-01} & \num{1.026163e+00}    & \num{5.796523e+00} & \num{1.035422e+02} \\
\runid & \texttt{na16} & \num{1.070288270712E-01} & \num{5.471500e-01} & \num{7.614410e-01} & \num{3.324996e-01} & \num{1.808062e-01} & \num{3.434212878346E-02}        & \num{3.567846e-01} & \num{1.022393e+00}    & \num{1.114069e+01} & \num{1.030157e+02} \\
\midrule
       & mean          & \num{9.514183e-02}       & \num{5.181288e-01} & \num{7.544664e-01} & \num{3.406460e-01} & \num{1.880744e-01} & \num{3.336907e-02}              & \num{3.381412e-01} & \num{1.022691e+00}    & \num{6.851845e+00} & \num{1.084713e+02} \\
\bottomrule
\end{tabular}
\end{table}

\begin{table}
\caption{Registration quality for the \nirep{} data for \claire{}. We consider an $H^1$-div regularization model ($H^1$-seminorm for $\vect{v}$) with $\beta_w=\num{1e-4}$. The regularization parameter $\beta_v$ is determined using a binary search with a bound of 0.25 for the determinant of the deformation gradient ($\beta_v =\num{5.500000e-04}$). We terminate the registration if the gradient is reduced by a factor of $\num{5e-2}$. We report (from left to right) the relative mismatch, the Dice coefficient (before and after registration), the false positive rate (after registration), the false negative rate (after registration), the relative reduction of the gradient, the extremal values of the determinant of the deformation gradient, and the overall runtime (in seconds). We execute the registration on CACDS's Opuntia server in single precision (see \secref{s:compute-systems} of the main manuscript for the specs).\label{t:nirep-regquality-h1sdiv-5d5e-4}}
\resetrunid\centering\scriptsize
\begin{tabular}{llllLlllllll}\toprule
\mcol{3}{\null}                                   & \multicolumn{2}{l}{dice}                & \mcol{3}{\null}                                                           & \multicolumn{3}{l}{$\det\igrad\vect{y}$}\\
\midrule
\mcol{2}{\null}        & mismatch                 & before             & after              & FP                 & FN                 & $\|\di{g}^\star\|_{\text{rel}}$ & $\min$             & $\operatorname{mean}$ & $\max$             & runtime     \\
\midrule
\runid & \texttt{na02} & \num{2.808164060116E-02} & \num{5.542543e-01} & \num{8.623370e-01} & \num{1.634204e-01} & \num{1.181392e-01} & \num{4.096447676420E-02}        & \num{4.746194e-01} & \num{1.010938e+00}    & \num{3.921268e+00} & \num{2.094917e+02} \\
\runid & \texttt{na03} & \num{2.684931084514E-02} & \num{5.039992e-01} & \num{8.327253e-01} & \num{2.320891e-01} & \num{1.210366e-01} & \num{4.089813679457E-02}        & \num{4.815015e-01} & \num{1.019933e+00}    & \num{7.199511e+00} & \num{2.187025e+02} \\
\runid & \texttt{na04} & \num{3.406877070665E-02} & \num{5.243747e-01} & \num{8.335772e-01} & \num{2.894796e-01} & \num{7.848105e-02} & \num{4.548207297921E-02}        & \num{3.410495e-01} & \num{1.030558e+00}    & \num{2.443181e+01} & \num{2.059377e+02} \\
\runid & \texttt{na05} & \num{4.021398350596E-02} & \num{5.550893e-01} & \num{8.541855e-01} & \num{1.964896e-01} & \num{1.080372e-01} & \num{4.759633541107E-02}        & \num{4.185913e-01} & \num{1.012542e+00}    & \num{5.214825e+00} & \num{2.031211e+02} \\
\runid & \texttt{na06} & \num{2.653059177101E-02} & \num{5.605455e-01} & \num{8.435593e-01} & \num{2.508350e-01} & \num{8.758519e-02} & \num{3.097482770681E-02}        & \num{5.229256e-01} & \num{1.020347e+00}    & \num{7.559669e+00} & \num{2.965448e+02} \\
\runid & \texttt{na07} & \num{2.504829131067E-02} & \num{5.311216e-01} & \num{8.512566e-01} & \num{2.374117e-01} & \num{8.303735e-02} & \num{3.945634141564E-02}        & \num{2.927341e-01} & \num{1.022311e+00}    & \num{3.658169e+00} & \num{2.156184e+02} \\
\runid & \texttt{na08} & \num{2.667831443250E-02} & \num{5.622983e-01} & \num{8.538373e-01} & \num{2.251264e-01} & \num{8.733844e-02} & \num{3.002730011940E-02}        & \num{3.273349e-01} & \num{1.015046e+00}    & \num{3.916790e+00} & \num{3.202449e+02} \\
\runid & \texttt{na09} & \num{2.871445752680E-02} & \num{5.077283e-01} & \num{8.389239e-01} & \num{2.305911e-01} & \num{1.108487e-01} & \num{4.478953778744E-02}        & \num{5.286912e-01} & \num{1.031077e+00}    & \num{1.013840e+01} & \num{2.196310e+02} \\
\runid & \texttt{na10} & \num{2.064519934356E-02} & \num{4.789702e-01} & \num{8.185274e-01} & \num{3.337895e-01} & \num{7.594708e-02} & \num{3.699371218681E-02}        & \num{6.026544e-01} & \num{1.037862e+00}    & \num{7.747598e+00} & \num{2.252374e+02} \\
\runid & \texttt{na11} & \num{2.617189474404E-02} & \num{4.579651e-01} & \num{8.293756e-01} & \num{2.913948e-01} & \num{8.505969e-02} & \num{4.601575061679E-02}        & \num{3.423367e-01} & \num{1.039678e+00}    & \num{2.184145e+01} & \num{2.294581e+02} \\
\runid & \texttt{na12} & \num{3.236174210906E-02} & \num{5.210443e-01} & \num{8.416418e-01} & \num{2.746737e-01} & \num{7.384552e-02} & \num{3.268095478415E-02}        & \num{5.069317e-01} & \num{1.030813e+00}    & \num{3.320764e+01} & \num{4.285205e+02} \\
\runid & \texttt{na13} & \num{3.190912306309E-02} & \num{5.307101e-01} & \num{8.075655e-01} & \num{3.152662e-01} & \num{1.092479e-01} & \num{4.736632853746E-02}        & \num{3.319518e-01} & \num{1.028464e+00}    & \num{8.115709e+00} & \num{2.127395e+02} \\
\runid & \texttt{na14} & \num{2.085218951106E-02} & \num{4.379558e-01} & \num{8.302875e-01} & \num{2.620444e-01} & \num{1.041733e-01} & \num{3.850091621280E-02}        & \num{3.312725e-01} & \num{1.039822e+00}    & \num{4.341471e+00} & \num{2.420825e+02} \\
\runid & \texttt{na15} & \num{2.927804738283E-02} & \num{4.987260e-01} & \num{8.318148e-01} & \num{2.475439e-01} & \num{1.116772e-01} & \num{4.277917742729E-02}        & \num{3.312725e-01} & \num{1.039822e+00}    & \num{4.341471e+00} & \num{2.027886e+02} \\
\runid & \texttt{na16} & \num{3.552267700434E-02} & \num{5.471500e-01} & \num{8.356300e-01} & \num{2.638943e-01} & \num{9.294462e-02} & \num{4.654541611671E-02}        & \num{3.731179e-01} & \num{1.030337e+00}    & \num{2.031620e+01} & \num{2.084002e+02} \\
\midrule
       & mean          & \num{2.886175e-02}       & \num{5.181288e-01} & \num{8.376830e-01} & \num{2.542700e-01} & \num{9.649327e-02} & \num{4.073809e-02}              & \num{4.137990e-01} & \num{1.027303e+00}    & \num{1.106347e+01} & \num{2.425679e+02} \\
\bottomrule
\end{tabular}
\end{table}

\begin{table}
\caption{Registration quality as a function of the regularization parameters $(\sigma_u, \sigma_v)$ for the DDEM algorithm. We consider the \emph{diffeomorphic update rule} with forces based on the \emph{gradient of the deformed template image} (default method; left table) and \emph{symmetrized forces} (right table). We limit this study to the datasets {\tt na02} (template image) and {\tt na01} (reference image) of the \nirep{} repository. We use three resolution levels with 15, 10, and 5 iterations on the individual levels, respectively. We report values for the relative change of the residual, the Dice coefficient, and the $\min$ and $\max$ of the determinant of the deformation gradient (computed via the considered \demons{} implementation). We highlight the best registration results (diffeomorphic deformation map and highest Dice score) in red.\label{t:ddem:sweep}}
\resetrunid\centering\scriptsize
\begin{tabular}{lcclLll}\toprule
       & \mcol{2}{parameters}    & \mcol{2}{\null}                               & \mcol{2}{$\det\igrad\vect{y}^{-1}$}           \\
\midrule
run    & $\sigma_u$ & $\sigma_d$ & mismatch              & dice                  & $\min$                & $\max$                \\
\midrule
\runid & 1.0        & 0.0        & \num{2.326337062E-02} & \num{8.912739868E-01} & $< 0$                 &                       \\
\runid & 2.0        &            & \num{4.707762226E-02} & \num{8.676730083E-01} & $< 0$                 &                       \\
\runid & 3.0        &            & \num{9.053304046E-02} & \num{8.244120686E-01} & $< 0$                 &                       \\
\runid & 4.0        &            & \num{1.507504433E-01} & \num{7.740749386E-01} & $< 0$                 &                       \\
\runid & 5.0        &            & \num{2.115481496E-01} & \num{7.382035391E-01} & $< 0$                 &                       \\
\rowcolor{red!20}
\runid & 6.0        &            & \num{2.835077345E-01} & \num{7.053460167E-01} & \num{2.309063375E-01} & \num{2.325822830E+00} \\
\runid & 7.0        &            & \num{3.488076925E-01} & \num{6.820082983E-01} & \num{3.964146376E-01} & \num{2.088901997E+00} \\
\runid & 8.0        &            & \num{3.860752881E-01} & \num{6.683613464E-01} & \num{4.650743902E-01} & \num{1.898323417E+00} \\
\runid & 0.0        & 1.0        & \num{5.586471409E-02} & \num{8.557338397E-01} & $< 0$                 &                       \\
\runid &            & 2.0        & \num{1.270513088E-01} & \num{7.947549498E-01} & $< 0$                 &                       \\
\rowcolor{red!20}
\runid &            & 3.0        & \num{2.525332272E-01} & \num{7.251251079E-01} & \num{2.050652057E-01} & \num{2.311142206E+00} \\
\runid &            & 4.0        & \num{3.732355833E-01} & \num{6.813779115E-01} & \num{5.042325258E-01} & \num{1.755817890E+00} \\
\runid &            & 5.0        & \num{4.460618794E-01} & \num{6.611711863E-01} & \num{6.121471524E-01} & \num{1.551860809E+00} \\
\runid &            & 6.0        & \num{5.222468972E-01} & \num{6.413094117E-01} & \num{7.003012896E-01} & \num{1.366043925E+00} \\
\runid &            & 7.0        & \num{5.846424103E-01} & \num{6.288886970E-01} & \num{7.586294413E-01} & \num{1.267131329E+00} \\
\runid &            & 8.0        & \num{6.248054504E-01} & \num{6.206149829E-01} & \num{7.877182364E-01} & \num{1.218394876E+00} \\
\bottomrule
\end{tabular}
\qquad\resetrunid\centering\scriptsize
\begin{tabular}{lcclLlll}\toprule
       & \mcol{2}{parameters}    & \mcol{2}{\null}                               & \mcol{2}{$\det\igrad\vect{y}^{-1}$}           \\
\midrule
run    & $\sigma_u$ & $\sigma_d$ & mismatch              & dice                  & $\min$                & $\max$                \\
\midrule
\runid & 1.0        & 0.0        & \num{2.947510965E-02} & \num{8.888131275E-01} & $< 0$                 &                       \\
\runid & 2.0        &            & \num{4.659114033E-02} & \num{8.664689294E-01} & $< 0$                 &                       \\
\runid & 3.0        &            & \num{8.699354529E-02} & \num{8.260966270E-01} & $< 0$                 &                       \\
\runid & 4.0        &            & \num{1.422277242E-01} & \num{7.771980392E-01} & $< 0$                 &                       \\
\runid & 5.0        &            & \num{1.944655180E-01} & \num{7.425547341E-01} & $< 0$                 &                       \\
\runid & 6.0        &            & \num{2.585209012E-01} & \num{7.106354627E-01} & $< 0$                 &                       \\
\rowcolor{red!20}
\runid & 7.0        &            & \num{3.113600612E-01} & \num{6.871142373E-01} & \num{3.200788796E-01} & \num{2.208411217E+00} \\
\runid & 8.0        &            & \num{3.590679169E-01} & \num{6.727600361E-01} & \num{4.076102376E-01} & \num{1.998745799E+00} \\
\runid & 0.0        & 1.0        & \num{5.305717513E-02} & \num{8.595240677E-01} & $< 0$                 &                       \\
\runid &            & 2.0        & \num{1.094698012E-01} & \num{8.054290753E-01} & $< 0$                 &                       \\
\rowcolor{red!20}
\runid &            & 3.0        & \num{2.117441297E-01} & \num{7.385951313E-01} & \num{3.005822562E-02} & \num{2.953015566E+00} \\
\runid &            & 4.0        & \num{3.197921813E-01} & \num{6.921746293E-01} & \num{3.478576839E-01} & \num{1.858527780E+00} \\
\runid &            & 5.0        & \num{3.967760503E-01} & \num{6.695722238E-01} & \num{5.314045548E-01} & \num{1.600968242E+00} \\
\runid &            & 6.0        & \num{4.741586745E-01} & \num{6.490836917E-01} & \num{6.404433250E-01} & \num{1.417105079E+00} \\
\runid &            & 7.0        & \num{5.373765826E-01} & \num{6.354826318E-01} & \num{7.091712952E-01} & \num{1.294634461E+00} \\
\runid &            & 8.0        & \num{5.777886510E-01} & \num{6.280042608E-01} & \num{7.440091968E-01} & \num{1.245173931E+00} \\
\bottomrule
\end{tabular}
\end{table}

\begin{table}
\caption{Registration quality as a function of the regularization parameters ($\sigma_u$, $\sigma_v$) for the LDDDEM algorithm. We consider the \emph{log-domain update rule} with forces based on the \emph{gradient of the deformed template image} (LDDDEM; left table) and the \emph{symmetric log-domain update rule} with \emph{symmetrized forces} (SLDDDEM; default method; right table). We limit this study to the datasets {\tt na02} (template image) and {\tt na01} (reference image) of the \nirep{} repository. We use three resolution levels with 15, 10, and 5 iterations on the individual levels, respectively. We report values for the relative change of the residual, the Dice coefficient, and the $\min$ and $\max$ of the determinant of the deformation gradient. The best registration (diffeomorphic and highest Dice score) is highlighted in red.\label{t:lddem:sweep}}
\resetrunid\centering\scriptsize
\begin{tabular}{lcclLll}\toprule
       & \mcol{2}{}    & \mcol{2}{mismatch}                            & \mcol{2}{$\det\igrad\vect{y}^{-1}$}                     \\
\midrule
run    & $\sigma_u$ & $\sigma_v$ & residual              & dice                  & $\min$                & $\max$                \\
\midrule
\runid & 1.0        & 0.0        & \num{5.280553550E-02} & \num{8.592480912E-01} & $< 0$                                         \\
\runid & 2.0        &            & \num{7.529053837E-02} & \num{8.344186267E-01} & $< 0$                                         \\
\runid & 3.0        &            & \num{1.162667722E-01} & \num{7.969439762E-01} & $< 0$                                         \\
\runid & 4.0        &            & \num{1.162667722E-01} & \num{7.582460126E-01} & $< 0$                                         \\
\runid & 5.0        &            & \num{2.302291691E-01} & \num{7.291434119E-01} & $< 0$                                         \\
\rowcolor{red!20}
\runid & 6.0        &            & \num{2.951605320E-01} & \num{7.008194050E-01} & \num{2.545172870E-01} & \num{2.316173792E+00} \\
\runid & 7.0        &            & \num{3.538450003E-01} & \num{6.794690534E-01} & \num{4.116433859E-01} & \num{2.088604689E+00} \\
\runid & 8.0        &            & \num{3.893532753E-01} & \num{6.666936656E-01} & \num{4.750167131E-01} & \num{1.901278734E+00} \\
\runid & 0.0        & 1.0        & \num{6.912566721E-02} & \num{8.405829037E-01} & $< 0$                                         \\
\runid &            & 2.0        & \num{1.324395835E-01} & \num{7.912953235E-01} & $< 0$                                         \\
\runid &            & 3.0        & \num{2.346027344E-01} & \num{7.292713844E-01} & $< 0$                                         \\
\rowcolor{red!20}
\runid &            & 4.0        & \num{3.575688004E-01} & \num{6.856127695E-01} & \num{4.370177388E-01} & \num{1.764521003E+00} \\
\runid &            & 5.0        & \num{4.312382936E-01} & \num{6.648984510E-01} & \num{5.815097094E-01} & \num{1.561049223E+00} \\
\runid &            & 6.0        & \num{5.057483315E-01} & \num{6.449780368E-01} & \num{6.878822446E-01} & \num{1.381839514E+00} \\
\runid &            & 7.0        & \num{5.665974617E-01} & \num{6.317722049E-01} & \num{7.460328341E-01} & \num{1.293170452E+00} \\
\runid &            & 8.0        & \num{6.073175073E-01} & \num{6.240371933E-01} & \num{7.769927979E-01} & \num{1.241199493E+00} \\
\bottomrule
\end{tabular}
\qquad\resetrunid\centering\scriptsize
\begin{tabular}{lcclLll}\toprule
       & \mcol{2}{}    & \mcol{2}{mismatch}                            & \mcol{2}{$\det\igrad\vect{y}^{-1}$}                     \\
\midrule
run    & $\sigma_u$ & $\sigma_v$ & residual              & dice                  & $\min$                & $\max$                \\
\midrule
\runid & 1.0        & 0.0        & \num{7.080791146E-02} & \num{8.523078638E-01} & $<0$                                          \\
\runid & 2.0        &            & \num{7.844036072E-02} & \num{8.356583216E-01} & $<0$                                          \\
\runid & 3.0        &            & \num{1.125537157E-01} & \num{8.003196139E-01} & $<0$                                          \\
\runid & 4.0        &            & \num{1.638007462E-01} & \num{7.620517053E-01} & $<0$                                          \\
\runid & 5.0        &            & \num{2.101594359E-01} & \num{7.338559193E-01} & $<0$                                          \\
\rowcolor{red!20}
\runid & 6.0        &            & \num{2.689733505E-01} & \num{7.055499502E-01} & \num{7.773168385E-02} & \num{2.332587004E+00} \\
\runid & 7.0        &            & \num{3.266469836E-01} & \num{6.838564571E-01} & \num{3.401958644E-01} & \num{2.055663347E+00} \\
\runid & 8.0        &            & \num{3.688737452E-01} & \num{6.704843177E-01} & \num{4.105187953E-01} & \num{1.880491972E+00} \\
\runid & 0.0        & 1.0        & \num{7.470855862E-02} & \num{8.356075482E-01} & $<0$                                          \\
\runid &            & 2.0        & \num{1.294154227E-01} & \num{7.903013996E-01} & $<0$                                          \\
\rowcolor{red!20}
\runid &            & 3.0        & \num{2.212747037E-01} & \num{7.334714177E-01} & \num{1.548299566E-03} & \num{2.637705326E+00} \\
\runid &            & 4.0        & \num{3.252030015E-01} & \num{6.914403234E-01} & \num{3.490672112E-01} & \num{1.780172467E+00} \\
\runid &            & 5.0        & \num{3.960536718E-01} & \num{6.702251884E-01} & \num{5.088196993E-01} & \num{1.556589603E+00} \\
\runid &            & 6.0        & \num{4.720625579E-01} & \num{6.497214296E-01} & \num{6.488980651E-01} & \num{1.389436603E+00} \\
\runid &            & 7.0        & \num{5.338659883E-01} & \num{6.363664751E-01} & \num{7.125169635E-01} & \num{1.302288532E+00} \\
\runid &            & 8.0        & \num{5.710264444E-01} & \num{6.290682123E-01} & \num{7.459799647E-01} & \num{1.254114509E+00} \\
\bottomrule
\end{tabular}
\end{table}

\begin{table}
\caption{Registration quality as a function of the regularization parameters $\sigma_d$ for the SDDEM algorithm (left) and $\sigma_v$ for the SLDDDEM algorithm (right). We set $\sigma_u = 1$ for both \demons{} variants. These two approaches gave us the best results based on the experiments reported in \tabref{t:ddem:sweep} and \tabref{t:lddem:sweep}. We limit this study to the datasets {\tt na02} (template image) and {\tt na01} (reference image) of the \nirep{} data. We use a multi-resolution approach with 15, 10, and 5 iterations per level (default setting) as a baseline. We increase the number of iterations per level by a factor of 2, 5, 10, and 100. We report values for the Dice coefficient and the $\min$ and $\max$ of the determinant of the deformation gradient.\label{t:dem-iterations}}
\resetrunid\centering\scriptsize
\begin{tabular}{lrclLll}\toprule
       & \mcol{2}{}    & \mcol{2}{mismatch}                              & \mcol{2}{$\det\igrad\vect{y}^{-1}$}             \\
\midrule
run    & iter & $\sigma_d$ & residual              & dice                  & $\min$                & $\max$                \\
\midrule
\runid &    1 & 2.5        & \num{1.602468044E-01} & \num{7.697818121E-01} & $<0$                                          \\
\runid &    2 &            & \num{1.487383097E-01} & \num{7.763855298E-01} & $<0$                                          \\
\runid &    5 &            & \num{1.412891746E-01} & \num{7.803777750E-01} & $<0$                                          \\
\runid &   10 &            & \num{1.396494508E-01} & \num{7.817192478E-01} & $<0$                                          \\
\runid &  100 &            & \num{1.368549168E-01} & \num{7.839196345E-01} & $<0$                                          \\
\runid &    1 & 3.0        & \num{2.117441297E-01} & \num{7.385951313E-01} & \num{3.005822562E-02} & \num{2.953015566E+00} \\
\runid &    2 &            & \num{1.982779354E-01} & \num{7.456242809E-01} & $<0$                                          \\
\runid &    5 &            & \num{1.911524534E-01} & \num{7.508293603E-01} & \num{2.373533696E-02} & \num{3.162640572E+00} \\
\runid &   10 &            & \num{1.891631484E-01} & \num{7.524607674E-01} & \num{8.167906851E-02} & \num{3.234217167E+00} \\
\runid &  100 &            & \num{1.860409975E-01} & \num{7.537887901E-01} & \num{5.574842542E-02} & \num{3.400515318E+00} \\
\runid &    1 & 3.5        & \num{2.701704800E-01} & \num{7.129391932E-01} & \num{2.394281179E-01} & \num{2.227959156E+00} \\
\runid &    2 &            & \num{2.524770498E-01} & \num{7.203035725E-01} & \num{2.079314888E-01} & \num{2.398765326E+00} \\
\runid &    5 &            & \num{2.376044840E-01} & \num{7.249019865E-01} & \num{2.020636797E-01} & \num{2.455271482E+00} \\
\runid &   10 &            & \num{2.377780080E-01} & \num{7.264445709E-01} & \num{2.288740277E-01} & \num{2.491136551E+00} \\
\runid &  100 &            & \num{2.332255840E-01} & \num{7.281488570E-01} & \num{1.942364275E-01} & \num{2.563106775E+00} \\
\bottomrule
\end{tabular}
\qquad\resetrunid\centering\scriptsize
\begin{tabular}{lrclLll}\toprule
       & \mcol{2}{}    & \mcol{2}{mismatch}                              & \mcol{2}{$\det\igrad\vect{y}^{-1}$}             \\
\midrule
run    & iter          & $\sigma_v$ & residual              & dice                  & $\min$                & $\max$       \\
\midrule
\runid &    1 & 2.5        & \num{1.710193157E-01} & \num{7.602645513E-01} & $<0$                                          \\
\runid &    2 &            & \num{1.688379198E-01} & \num{7.660483517E-01} & $<0$                                          \\
\runid &    5 &            & \num{1.618115455E-01} & \num{7.694946951E-01} & $<0$                                          \\
\runid &   10 &            & \num{1.599550247E-01} & \num{7.707331051E-01} & $<0$                                          \\
\runid &  100 &            & \num{1.649025977E-01} & \num{7.700247098E-01} & $<0$                                          \\
\runid &    1 & 3.0        & \num{2.212747037E-01} & \num{7.334714177E-01} & \num{1.548299566E-03} & \num{2.637705326E+00} \\
\runid &    2 &            & \num{2.059686780E-01} & \num{7.436279758E-01} & \num{1.458691340E-02} & \num{2.747811556E+00} \\
\runid &    5 &            & \num{2.059686780E-01} & \num{7.436279758E-01} & \num{1.458691340E-02} & \num{2.747811556E+00} \\
\runid &   10 &            & \num{2.046824694E-01} & \num{7.447072963E-01} & \num{1.217605397E-01} & \num{2.800973415E+00} \\
\runid &  100 &            & \num{2.063782662E-01} & \num{7.447238786E-01} & \num{1.293125749E-01} & \num{2.919600010E+00} \\
\runid &    1 & 3.5        & \num{2.728112042E-01} & \num{7.103996835E-01} & \num{2.342625856E-01} & \num{2.076637030E+00} \\
\runid &    2 &            & \num{2.590842843E-01} & \num{7.168705380E-01} & \num{1.946153939E-01} & \num{2.165000677E+00} \\
\runid &    5 &            & \num{2.533491254E-01} & \num{7.205800657E-01} & \num{2.325878441E-01} & \num{2.244817019E+00} \\
\runid &   10 &            & \num{2.485840917E-01} & \num{7.220047818E-01} & \num{2.952381968E-01} & \num{2.277714014E+00} \\
\runid &  100 &            & \num{2.481485307E-01} & \num{7.229721018E-01} & \num{2.775656879E-01} & \num{2.322824478E+00} \\
\bottomrule
\end{tabular}
\end{table}

\begin{table}
\caption{Registration quality for the diffeomorphic \demons{} algorithm for the entire \nirep{} data. We report results for SDDEM (\emph{diffeomorphic update rule}; force: \emph{symmetrized}) for varying regularization parameters $(\sigma_u,\sigma_d)$. We report values for the Dice coefficient, and the $\min$ and $\max$ values of the determinant of the deformation gradient. The bottom row provides the mean values across all 15 runs for each individual method/setting, respectively. The runs were executed with a multi-resultion approach with 150, 100, and 50 iterations per level.\label{t:demons-all-nirep-data}}
\resetrunid\centering\scriptsize
\begin{tabular}{clLlllLlllLlllLll}\toprule
\mcol{3}{}                                     & \mcol{2}{$\det\igrad\vect{y}^{-1}$}           & \mcol{2}{}                     & \mcol{2}{$\det\igrad\vect{y}^{-1}$}           & \mcol{2}{}                     & \mcol{2}{$\det\igrad\vect{y}^{-1}$}           \\\midrule
data          & run    & dice                  & $\min$                & $\max$                & run    & dice                  & $\min$                & $\max$                & run    & dice                  & $\min$                & $\max$                \\\midrule
              & \mcol{4}{SDDEM($0,3$)}                                                         & \mcol{4}{SDDEM($0,3.5$)}                                                       & \mcol{4}{SDDEM($0,1$)}                                                         \\\midrule
\texttt{na02} & \runid & \num{7.524607674E-01} & \num{8.167906851E-02} & \num{3.234217167E+00} & \runid & \num{7.264445709E-01} & \num{2.288740277E-01} & \num{2.491136551E+00} & \runid & \num{8.607238668E-01} & $<0$                  & \num{3.838982391E+01} \\
\texttt{na03} & \runid & \num{7.264901779E-01} & $<0$                  & \num{2.911321640E+00} & \runid & \num{7.058524654E-01} & \num{1.775113046E-01} & \num{2.418567896E+00} & \runid & \num{8.045950690E-01} & $<0$                  & \num{2.848730659E+01} \\
\texttt{na04} & \runid & \num{7.520353422E-01} & $<0$                  & \num{3.362325668E+00} & \runid & \num{7.278378233E-01} & \num{1.533267945E-01} & \num{2.824498177E+00} & \runid & \num{8.525523138E-01} & $<0$                  & \num{3.622441101E+01} \\
\texttt{na05} & \runid & \num{7.448940663E-01} & $<0$                  & \num{2.852404594E+00} & \runid & \num{7.179882911E-01} & \num{1.142086461E-01} & \num{2.440938711E+00} & \runid & \num{8.552234706e-01} & $<0$                  & \num{2.021568871E+01} \\
\texttt{na06} & \runid & \num{7.581988161E-01} & \num{9.624088183E-03} & \num{3.629817724E+00} & \runid & \num{7.338353626E-01} & \num{1.332834810E-01} & \num{2.792528868E+00} & \runid & \num{8.654651783E-01} & $<0$                  & \num{2.264057732E+01} \\
\texttt{na07} & \runid & \num{7.573275130E-01} & $<0$                  & \num{3.790106058E+00} & \runid & \num{7.318410853E-01} & \num{1.311681867E-01} & \num{2.800957203E+00} & \runid & \num{8.617095407E-01} & $<0$                  & \num{3.364328384E+01} \\
\texttt{na08} & \runid & \num{7.558222534E-01} & $<0$                  & \num{3.059933186E+00} & \runid & \num{7.307236822E-01} & $<0$                  & \num{2.399970531E+00} & \runid & \num{8.648908419E-01} & $<0$                  & \num{2.369104195E+01} \\
\texttt{na09} & \runid & \num{7.390791963E-01} & $<0$                  & \num{3.436971664E+00} & \runid & \num{7.160904978E-01} & \num{6.185429171E-02} & \num{2.640543938E+00} & \runid & \num{8.263383450E-01} & $<0$                  & \num{3.046219254E+01} \\
\texttt{na10} & \runid & \num{7.515453527E-01} & $<0$                  & \num{3.464956045E+00} & \runid & \num{7.292254417E-01} & \num{2.118150890E-01} & \num{2.577309132E+00} & \runid & \num{8.445947494E-01} & $<0$                  & \num{2.355178642E+01} \\
\texttt{na11} & \runid & \num{7.525696982E-01} & $<0$                  & \num{3.075066328E+00} & \runid & \num{7.292198609E-01} & \num{1.507869363E-01} & \num{2.339290142E+00} & \runid & \num{8.491450883E-01} & $<0$                  & \num{2.470713234E+01} \\
\texttt{na12} & \runid & \num{7.378139167E-01} & $<0$                  & \num{3.654606342E+00} & \runid & \num{7.164448116E-01} & \num{1.057251617E-01} & \num{2.803232431E+00} & \runid & \num{8.215823101E-01} & $<0$                  & \num{3.250938797E+01} \\
\texttt{na13} & \runid & \num{7.510133739E-01} & $<0$                  & \num{3.129977465E+00} & \runid & \num{7.289087682E-01} & \num{1.010739729E-01} & \num{2.494370699E+00} & \runid & \num{8.345589066E-01} & $<0$                  & \num{3.215924454E+01} \\
\texttt{na14} & \runid & \num{7.449150251E-01} & $<0$                  & \num{2.820612431E+00} & \runid & \num{7.230210744E-01} & \num{8.664506674E-02} & \num{2.224625587E+00} & \runid & \num{8.287514014E-01} & $<0$                  & \num{2.646317673E+01} \\
\texttt{na15} & \runid & \num{7.341751715E-01} & $<0$                  & \num{2.811891079E+00} & \runid & \num{7.105671221E-01} & \num{1.060965136E-01} & \num{2.203224182E+00} & \runid & \num{8.343307651E-01} & $<0$                  & \num{2.299076080E+01} \\
\texttt{na16} & \runid & \num{7.447690043E-01} & $<0$                  & \num{3.195658922E+00} & \runid & \num{7.228238394E-01} & \num{1.513934731E-01} & \num{2.509691477E+00} & \runid & \num{8.298199119E-01} & $<0$                  & \num{3.381185150E+01}
\\\midrule
\mcol{2}{mean}         & \num{7.468739783E-01} & \num{6.086877112E-03} & \num{3.228657754E+00} &        & \num{7.233883131E-01} & \num{1.275841964E-01} & \num{2.530725702E+00} &        & \num{8.422854506E-01} & 0.00                  & \num{2.866317774e+01}
\\\bottomrule
\end{tabular}
\end{table}

We refer to \secref{s:reg-quality} of the main manuscript for a discussion of these results.

\section{Optimality Conditions}
\label{s:optimality-cond-dev}

As we have mentioned in \secref{s:methods} of the main manuscript, we consider an optimize-then-discretize approach. We have seen in~\eqref{e:lagrangian} that the Lagrangian of the optimization problem in~\eqref{e:varopt} is given by (for simplicity, we consider the compressible formulation and neglect the boundary conditions)
\begin{align}
\fun{L}[m,\lambda,\vect{v}] \defeq &
  \half{1} \iom{(m(t=1) - m_R)^2}
+ \half{\beta} \big\langle\dop{A}\vect{v},\vect{v}\big\rangle_{L^2(\Omega)^d}
+ \int_0^1
\langle
\p_t m + \igrad m \cdot \vect{v},\lambda
\rangle_{L^2(\Omega)}\d{t}
\\\nonumber
&\quad
+ \langle m(t=0) - m_T,\lambda(t=0)\rangle_{L^2(\Omega)},
\end{align}

\noindent where $\dop{A} = \dop{B}^\ast\dop{B}$ is the self-adjoint regularization operator (e.g., a vector-Laplacian $-\ilap$) and $\langle\,\cdot\,,\,\cdot\,\rangle_{L^2(\Omega)}$ is the standard $L^2$ inner product. Our derivation will be formal only. That is, we assume that all variables and functions meet the regularity requirements to be able to carry out the necessary computations. We know from Lagrange multiplier theory that the first variations of the Lagrangian with respect to all variables have to vanish for an admissible solution of~\eqref{e:varopt}. The variations of $\fun{L}$ with respect to the state, adjoint, and control variables are given by
\[
\begin{aligned}
\fun{L}_{m}[m,\lambda,\vect{v}](\hat{m})
& = \langle\lambda(t=1) - m(t=1) + m_R,\hat{m}(t=1)\rangle_{L^2(\Omega)}
+ \int_0^1\langle-\p_t \lambda - \idiv \lambda\vect{v},\hat{m}\rangle_{L^2(\Omega)}\d{t} \\
\fun{L}_{\lambda}[m,\lambda,\vect{v}](\hat{\lambda})
& = \int_0^1\langle\p_t m + \vect{v}\cdot \igrad m ,\hat{\lambda}\rangle_{L^2(\Omega)}\d{t}
+ \langle m(t=0) - m_T,\hat{\lambda}(t=0)\rangle_{L^2(\Omega)} \\
\fun{L}_{v}[m,\lambda,\vect{v}](\vect{\hat{v}})
&=
\int_0^1\langle\lambda \igrad m,\hat{\vect{v}}\rangle_{L^2(\Omega)^d}\d{t}
+ \beta \big\langle\dop{A}\vect{v},\vect{\hat{v}}\big\rangle_{L^2(\Omega)^d},
\end{aligned}
\]

\noindent respectively. To drive the variation of $\fun{L}$ with respect to the state variable $m$, we apply integration by parts. In particular, collecting all terms that depend on $m$ we have
\[
\begin{aligned}
& \half{1}\langle m(t=1) - m_R, m(t=1) - m_R\rangle_{L^2(\Omega)}
+ \int_0^1\langle\p_t m + \igrad m \cdot \vect{v},\lambda\rangle_{L^2(\Omega)}\d{t}
+ \langle m(t=0),\lambda(t=0)\rangle_{L^2(\Omega)} \\
= \;& \half{1}\langle m(t=1) - m_R, m(t=1) - m_R\rangle_{L^2(\Omega)}
  + \langle m(t=1),\lambda(t=1)\rangle_{L^2(\Omega)} - \langle m(t=0),\lambda(t=0)\rangle_{L^2(\Omega)} \\
& + \int_0^1\langle-\p_t \lambda - \idiv \lambda\vect{v},m\rangle_{L^2(\Omega)}\d{t}
+ \langle m(t=0),\lambda(t=0)\rangle_{L^2(\Omega)}\\
=\; &
\half{1}\langle m(t=1) - m_R, m(t=1) - m_R\rangle_{L^2(\Omega)}
+ \langle m(t=1),\lambda(t=1)\rangle_{L^2(\Omega)}
+ \int_0^1\langle-\p_t \lambda - \idiv \lambda\vect{v},m\rangle_{L^2(\Omega)}\d{t}.
\end{aligned}
\]

\noindent Computing variations with respect to $m$ results in the expression given above. Suppose that $\vect{\phi}^\star \defeq (m^\star,\vect{v}^\star,\lambda^\star)$ are any primal and dual optimal points with zero duality gap. Then, the strong form of the KKT (first order optimality) conditions is given by
\[
\begin{aligned}
                       \p_t m^\star + \vect{v}^\star \cdot \igrad m^\star & = 0        && \quad\text{in } \Omega\times(0,1] \\
                                                            m^\star - m_T & = 0        && \quad\text{in } \Omega\times\{0\} \\
                  -\p_t \lambda^\star - \idiv \lambda^\star\vect{v}^\star & = 0        && \quad\text{in } \Omega\times[0,1) \\
                                            \lambda^\star - m^\star + m_R & = 0        && \quad\text{in } \Omega\times\{1\} \\
\beta \dop{A}\vect{v}^\star + \int_0^1 \lambda^\star \igrad m^\star \d{t} & = \vect{0} && \quad\text{in } \Omega.
\end{aligned}
\]

We apply Newtons method to solve the system given above. We have to compute second-order variations of the Lagrangian $\fun{L}$. Formally, we have
\[
\begin{aligned}
\fun{L}_{mm}[m,\lambda,\vect{v}](\tilde{m},\hat{m})
& = \langle \tilde{m}(t=1),\hat{m}(t=1)\rangle_{L^2(\Omega)} \phantom{\int_0^1}\\
\fun{L}_{m \lambda}[m,\lambda,\vect{v}](\tilde{\lambda},\hat{m})
& =  \langle \tilde{\lambda}(t=1),\hat{m}(t=1)\rangle
  + \int_0^1\langle -\p_t \tilde{\lambda} - \idiv \tilde{\lambda}\vect{v},\hat{m}\rangle_{L^2(\Omega)}\d{t} \\
\fun{L}_{m\vect{v}}[m,\lambda,\vect{v}](\vect{\tilde{v}},\vect{\hat{v}})
& = \int_0^1\langle-\idiv\lambda\vect{\tilde{v}},\hat{m}\rangle_{L^2(\Omega)^d} \d{t} \\
\fun{L}_{\lambda m}[m,\lambda,\vect{v}](\tilde{m},\hat{\lambda})
& = \langle\tilde{\lambda}(t=1),\hat{\lambda}(t=0)\rangle_{L^2(\Omega)}
  + \int_0^1\langle\p_t\tilde{m}+\igrad\tilde{m}\cdot\vect{v},\hat{\lambda}\rangle_{L^2(\Omega)}\d{t}\\
\fun{L}_{\lambda\lambda}[m,\lambda,\vect{v}](\tilde{\lambda},\hat{\lambda})
& = 0 \phantom{\int_0^1}\\
\fun{L}_{\lambda\vect{v}}[m,\lambda,\vect{v}](\vect{\tilde{v}},\vect{\hat{v}})
& = \int_0^1\langle\igrad m\cdot\vect{\tilde{v}},\hat{\lambda}\rangle_{L^2(\Omega)^d} \d{t}\\
\fun{L}_{\vect{v}m}[m,\lambda,\vect{v}](\tilde{m},\vect{\hat{v}})
& = \int_0^1\langle\lambda\igrad\tilde{m},\vect{\hat{v}}\rangle_{L^2(\Omega)^d}\d{t} \\
\fun{L}_{\vect{v}\lambda }[m,\lambda,\vect{v}](\tilde{\lambda},\vect{\hat{v}})
& = \int_0^1\langle\tilde{\lambda}\igrad m,\vect{\hat{v}} \rangle_{L^2(\Omega)^d}\d{t} \\
\fun{L}_{\vect{v}\vect{v}}[m,\lambda,\vect{v}](\vect{\tilde{v}},\vect{\hat{v}})
& = \beta \big\langle\dop{A}\vect{\tilde{v}},\vect{\hat{v}}\big\rangle_{L^2(\Omega)^d} \phantom{\int_0^1}.
\end{aligned}
\]

To illustrate the discretized KKT system, we collect these terms in a matrix. We obtain
\[
\begin{bmatrix*}
\di{H}_{mm} & \di{H}_{mv}         & \di{A}^{\!\T} \\
\di{H}_{vm} & \di{H}_{\text{reg}} & \di{C}^{\!\T} \\
\di{A}        & \di{C}              & \di{0}        \\
\end{bmatrix*}
\begin{bmatrix*}
\tilde{\di{m}}        \\
\tilde{\di{v}}        \\
\tilde{\dig{\lambda}} \\
\end{bmatrix*}
=
-
\begin{bmatrix*}
\di{g}_{m}       \\
\di{g}_{v}       \\
\di{g}_{\lambda} \\
\end{bmatrix*}.
\]

Here, $\tilde{\di{m}}\in\ns{R}^{nn_t}$, $\tilde{\di{v}}\in\ns{R}^{3n}$, and $\tilde{\dig{\lambda}}\in\ns{R}^{nn_t}$ are the search directions, $\di{g}_m\in\ns{R}^{nn_t}$, $\di{g}_v\in\ns{R}^{3n}$, and $\dig{\lambda}\in\ns{R}^{nn_t}$ are the discrete gradients (first variation) of the Lagrangian with respect the state, control and adjoint variable, respectively. Further, $\di{H}_{mm}\in\ns{R}^{nn_t,nn_t}$, $\di{H}_{vm}\in\ns{R}^{nn_t,3n}$, $\di{H}_{mv}\in\ns{R}^{3n,nn_t}$, and $\di{H}_{\text{reg}}\in\ns{R}^{3n,3n}$ are components of the Hessian matrix of the Lagrangian functional, and $\di{A}\in\ns{R}^{nn_t,nn_t}$ and $\di{C}\in\ns{R}^{nn_t,3n}$ are the Jacobian of the state equation with respect to the state and control variables, respectively. More precisely, $\di{H}_{mm}$, $\di{H}_{mv}$, $\di{H}_{vm}$, and $\di{H}_{\text{reg}}$ correspond to $\fun{L}_{mm}$, $\fun{L}_{mv}$, $\fun{L}_{vm}$, and $\fun{L}_{vv}$, respectively, and $\di{A}$ and $\di{C}$ and their transposes correspond to $\fun{L}_{\lambda m}$, $\fun{L}_{vm}$, $\fun{L}_{m\lambda}$, and $\fun{L}_{v\lambda }$, respectively. Under the assumption that $\di{m}$ and $\dig{\lambda}$ fulfill the state and adjoint equations exactly, we have $\di{g}_m = \di{g}_\lambda = \vect{0}$. We use this assumption to eliminate the incremental state and adjoint variables from the KKT system. We obtain $\tilde{\di{m}} = -\di{A}^{-1}\di{C}\tilde{\di{v}}$ and $\tilde{\dig{\lambda}} = -\di{A}^{\!-\T}(\di{H}_{mm} \tilde{\di{m}} + \di{H}_{mv}\tilde{\di{v}})$. The reduced space linear system for the Newton step in the control variable is given by $\di{H}\tilde{\di{v}} = -\di{g}_v$, where $\di{H} \defeq \di{H}_{\text{reg}} + \di{C}^{\!\T}\di{A}^{\!-\T}(\di{H}_{mm}\di{A}^{-1}\di{C} - \di{H}_{mv}) - \di{H}_{vm}\di{A}^{-1}\di{C}$. Notice that the reduced Hessian involves inverses of the state and adjoint operators. This makes $\di{H}$ a dense matrix that is often too large to be computed and/or stored. However, we can define an expression for the Hessian matvec (application of $\di{H}$ to a vector) that involves solving linear systems with the matrices $\mat{A}$ and $\mat{A}^\T$. This approach corresponds to the reduced space method described in the main manuscript.

If we consider the incompressibility constraint $\idiv \vect{v}$ the KKT conditions are given by
\[
\begin{aligned}
                                        \p_t m^\star + \igrad m^\star \cdot \vect{v}^\star & = 0        && \quad\text{in } \Omega\times(0,1] \\
                                                                             m^\star - m_T & = 0        && \quad\text{in } \Omega\times\{0\} \\
                                                                      \idiv \vect{v}^\star & = 0        && \quad\text{in } \Omega            \\
                                   -\p_t \lambda^\star - \idiv \lambda^\star\vect{v}^\star & = 0        && \quad\text{in } \Omega\times[0,1) \\
                                                             \lambda^\star - m^\star + m_R & = 0        && \quad\text{in } \Omega\times\{1\} \\
\beta \dop{A}\vect{v}^\star + \igrad p^\star + \int_0^1 \lambda^\star \igrad m^\star \d{t} & = \vect{0} && \quad\text{in } \Omega.
\end{aligned}
\]

\noindent with primal and dual optimal points $\vect{\phi}^\star \defeq (m^\star,\vect{v}^\star,\lambda^\star,p^\star)$, where $p^\star$ is the optimal Lagrange multiplier (dual variable) for the incompressibility constraint $\idiv\vect{v} = 0$. This incompressibility constraint can be eliminated from the KKT conditions. Suppose that $\dop{A} = -\ilap$. Computing the divergence of the control equation yields
\[
   \idiv (-\ilap\vect{v}^\star + \igrad p^\star + \int_0^1 \lambda^\star \igrad m^\star \d{t})
 =-\idiv\ilap\vect{v}^\star + \ilap p^\star + \idiv \int_0^1 \lambda^\star \igrad m^\star \d{t} = \vect{0}.
\]

\noindent Using the optimality condition $\idiv \vect{v}^\star  = 0$ and the identity $\ilap\vect{v}^\star = \igrad(\idiv \vect{v}^\star) - \icurl(\icurl \vect{v}^\star)$  we obtain
\[
 p^\star = \ilap^{-1} \idiv \int_0^1 \lambda^\star \igrad m^\star \d{t}.
\]

\noindent Inserting $p^\star$ into the control equation yields
\[
\beta \dop{A}\vect{v}^\star + \igrad \ilap^{-1} \idiv \int_0^1 \lambda^\star \igrad m^\star \d{t} + \int_0^1 \lambda^\star \igrad m^\star \d{t} = \vect{0}.
\]

After elimination of the pressure $p^\star$, the strong form of the KKT conditions are given by
\[
\begin{aligned}
                                        \p_t m^\star + \igrad m^\star \cdot \vect{v}^\star & = 0        && \quad\text{in } \Omega\times(0,1] \\
                                                                             m^\star - m_T & = 0        && \quad\text{in } \Omega\times\{0\} \\
                                   -\p_t \lambda^\star - \idiv \lambda^\star\vect{v}^\star & = 0        && \quad\text{in } \Omega\times[0,1) \\
                                                             \lambda^\star - m^\star + m_R & = 0        && \quad\text{in } \Omega\times\{1\} \\
\beta \dop{A}\vect{v}^\star + \igrad \ilap^{-1} \idiv \int_0^1 \lambda^\star \igrad m^\star \d{t} + \int_0^1 \lambda^\star \igrad m^\star \d{t} & = \vect{0} && \quad\text{in } \Omega.
\end{aligned}
\]

The same algorithm discussed above can be used to solve for an optimal $\vect{v}^\star$. The derivations necessary to setup the KKT system are along the same lines as without incompressibility constraint.

\end{document}